\documentclass[11pt,a4paper,reqno,final]{article}

\newif\iffinal
\finaltrue

\usepackage{amsmath}
\usepackage{amsthm}
\usepackage{enumerate}
\usepackage{amssymb}

\usepackage{enumerate}
\usepackage{verbatim}
\usepackage{url}
\usepackage[latin1]{inputenc}

\usepackage{caption}
\usepackage{geometry}
\usepackage{graphicx,color}
\usepackage{fullpage}

\usepackage{tikz}
\usetikzlibrary{calc}
\usepackage{xifthen}
\pgfmathsetseed{1234}
\usetikzlibrary{decorations.pathreplacing}

\usepackage[obeyDraft]{todonotes}

\def\clap#1{\hbox to 0pt{\hss#1\hss}}

  \newcounter{constant}
  \newcommand{\nc}[1]{\refstepcounter{constant}\label{#1}}
  \newcommand{\uc}[1]{c_{\ref{#1}}}
\usepackage{array}
\newcolumntype{e}{>{\displaystyle}r @{\,} >{\displaystyle}c @{\,} >{\displaystyle}l}

\def\arraypar#1{\parbox[c]{\textwidth - 2cm}{\centering #1}}

\usepackage{environ}
\NewEnviron{display}{
\begin{equation}\begin{array}{c}
\arraypar{\BODY}
\end{array}\end{equation}
}

\numberwithin{equation}{section}

\newcommand{\mcup}{{\textstyle \bigcup\limits}}
\newcommand{\mcap}{{\textstyle \bigcap\limits}}

\usepackage[notcite,notref,color]{showkeys}
\definecolor{labelkey}{gray}{.80}

\usepackage[colorlinks]{hyperref}
\usepackage[figure,table]{hypcap} 
\hypersetup{final}
\usepackage{xcolor}
\hypersetup{
    colorlinks,
    linkcolor={blue!50!black},
    citecolor={blue!50!black},
    urlcolor={blue!80!black}
}

\def\Z{{\mathbb Z}}

\def\R{{\mathbb R}}

\def\Fc{{\mathcal F}}

\renewcommand*\d{\mathop{}\!\mathrm{d}}

\def\be{\begin{equation}}
\def\ee{\end{equation}}
\def\bea{\begin{equation*}}
\def\eea{\end{equation*}}
\def\bal{\begin{aligned}}
\def\eal{\end{aligned}}
\def\eps{\varepsilon}

\DeclareMathOperator{\cross}{Cross}
\DeclareMathOperator{\cir}{Circ}
\DeclareMathOperator{\arm}{Arm}

\def\b{\tfrac{5}{4}}
\def\bb{\tfrac{4}{3}}
\def\bOverTwo{\tfrac{2}{3}}

\def\x{r}

\def\xprime{{r'}}

\def\y{h}

\def\Pr{{\mathbb P}}
\DeclareMathOperator{\E}{{\mathbb E}}
\DeclareMathOperator{\Var}{Var}
\DeclareMathOperator{\Cov}{Cov}

\DeclareMathOperator{\ind}{{\bf 1}}

\DeclareMathOperator{\Inf}{Inf}
\DeclareMathOperator{\supp}{supp}

\DeclareMathOperator{\interior}{int}
\DeclareMathOperator{\closure}{cl}

\newtheorem{thm}{Theorem}[section]
\newtheorem{lma}[thm]{Lemma}
\newtheorem{cor}[thm]{Corollary}
\newtheorem{prop}[thm]{Proposition}

\newtheorem{definition}[thm]{Definition}

\theoremstyle{remark}
\newtheorem{remark}{Remark}
\newtheorem{preex}[thm]{Example}

\theoremstyle{definition}
\newtheorem*{acknow}{Acknowledgements}

\def\dual{densely dotted}

\newcommand{\boxRectangle}[3][solid]{\vcenter{\hbox{\;\tikz[scale=1.5]{
    \draw (0, 0) rectangle (6ex, 4ex);
    \draw [#1] (0, 2ex) .. controls (2.6ex, 3ex) and (3.3ex, 1ex) .. (6ex, 2ex);
    \draw (3ex, .2ex) node[anchor=north] {$\scriptstyle #2$};
    \draw (5.7ex, 2ex) node[anchor=west] {$\scriptstyle #3$};
  }}}
}

\newcommand{\boxThreeFourUp}[2]{\vcenter{\hbox{\tikz[scale=1.5]{
    \draw (0, 0) rectangle (3ex, 4ex);
    \draw (2.8ex, 2ex) -- (3.2ex, 2ex);
    \draw (0, 2ex) .. controls (1ex, 3ex) and (2ex, 1ex) .. (3ex, 3ex);
    \draw (1.5ex, .2ex) node[anchor=north] {$\scriptstyle #1$};
    \draw (0ex, 2ex) node[anchor=east] {$\scriptstyle #2$};
  }\;}}
}

\newcommand{\boxThreeFourDown}[2]{\vcenter{\hbox{\tikz[scale=1.5]{
    \draw (0, 0) rectangle (3ex, 4ex);
    \draw (2.8ex, 2ex) -- (3.2ex, 2ex);
    \draw (0, 2ex) .. controls (1ex, 3ex) and (2ex, 1ex) .. (3ex, 1ex);
    \draw (1.5ex, .2ex) node[anchor=north] {$\scriptstyle #1$};
    \draw (0ex, 2ex) node[anchor=east] {$\scriptstyle #2$};
  }\;}}
}

\newcommand{\boxThreeFourMiddle}[4][solid]{\vcenter{\hbox{\;\tikz[scale=1.5]{
    \draw (0, 0) rectangle (3ex, 4ex);
    \draw (2.8ex, 2ex) -- (3.2ex, 2ex) (2.8ex, 2.8ex) -- (3.2ex, 2.8ex);
    \draw [#1](0, 2ex) .. controls (1ex, 3.5ex) and (2ex, 1ex) .. (3ex, 2.6ex);
    \draw (1.5ex, .2ex) node[anchor=north] {$\scriptstyle #2$};
    \draw (.4ex, 2ex) node[anchor=east] {$\scriptstyle #3$};
    \draw (2.8ex, 2.4ex) node[anchor=west] {$\scriptstyle #4$};
  }}}
}

\newcommand{\boxThreeFourSmall}[4][solid]{\vcenter{\hbox{\;\tikz[scale=1.5]{
    \draw (0, 0) rectangle (3ex, 4ex);
    \draw (2.8ex, 2.2ex) -- (3.2ex, 2.2ex) (2.8ex, 2.5ex) -- (3.2ex, 2.5ex);
    \draw [#1](0, 2ex) .. controls (1ex, 3.5ex) and (2ex, 2ex) .. (3ex, 2.4ex);
    \draw (1.5ex, .2ex) node[anchor=north] {$\scriptstyle #2$};
    \draw (.4ex, 2ex) node[anchor=east] {$\scriptstyle #3$};
    \draw (2.8ex, 2.4ex) node[anchor=west] {$\scriptstyle #4$};
  }}}
}

\newcommand{\boxFourThree}[3][solid]{\vcenter{\hbox{\tikz[scale=1.5]{
    \draw (0,0) rectangle (4ex, 3ex);
    \draw [#1] (0, 1.5ex) .. controls (1.3ex, 2.5ex) and (2.6ex, .5ex) .. (4ex, 1.5ex);
    \draw (2ex, .2ex) node[anchor=north] {$\scriptstyle #2$};
    \draw (3.6ex, 1.5ex) node[anchor=west] {$\scriptstyle #3$};
  }}}
}

\newcommand{\boxHighTwoCross}[2]{\vcenter{\hbox{\tikz[scale=1.5]{
    \draw (0, 0) rectangle (2ex, 6ex);
    \draw (1.8ex, 3ex) -- (2.2ex, 3ex);
    \draw [\dual] (0, 5ex) .. controls (.6ex, 6ex) and (1.4ex, 4ex) .. (2ex, 5ex);
    \draw (0, 3.5ex) .. controls (.6ex, 5ex) and (1.4ex, 2.5ex) .. (2ex, 4ex);
    \draw (1ex, .2ex) node[anchor=north] {$\scriptstyle #1$};
    \draw (0ex, 3ex) node[anchor=east] {$\scriptstyle #2$};
  }\;}}
}

\newcommand{\boxHighTwoMiddle}[4][color=white]{\vcenter{\hbox{\tikz[scale=1.5]{
    \draw (0, 2ex) .. controls (1.4ex, 1ex) and (.6ex, 3ex) .. (2ex, 3.3ex);
    \draw [#1] (0, 4ex) .. controls (1.4ex, 5ex) and (.6ex, 3.2ex) .. (2ex, 3.7ex);
    \draw (0, 0) rectangle (2ex, 6ex);
    \draw (1.8ex, 2ex) -- (2.2ex, 2ex) (1.8ex, 3ex) -- (2.2ex, 3ex) (1.8ex, 4ex) -- (2.2ex, 4ex);
    \draw [decorate,decoration={brace,amplitude=.5ex,raise=.4ex},yshift=0pt] (2.2ex, 4ex) -- (2.2ex, 2ex);
    \draw (1ex, .2ex) node[anchor=north] {$\scriptstyle #2$};
    \draw (0ex, 3ex) node[anchor=east] {$\scriptstyle #3$};
    \draw (2.4ex, 3ex) node[anchor=west] {$\scriptstyle #4$};
  }}}
}

\newcommand{\boxHighSingle}[3]{\vcenter{\hbox{\tikz[scale=1.5]{
    \draw (0, 2ex) .. controls (1.4ex, 1ex) and (.6ex, 3.3ex) .. (2ex, 3.5ex);
    \draw (0, 0) rectangle (2ex, 6ex);
    \draw (1.8ex, 3.3ex) -- (2.2ex, 3.3ex) (1.8ex, 3.8ex) -- (2.2ex, 3.8ex);
    \draw (1ex, .2ex) node[anchor=north] {$\scriptstyle #1$};
    \draw (0ex, 3ex) node[anchor=east] {$\scriptstyle #2$};
    \draw (2ex, 3.5ex) node[anchor=west] {$\scriptstyle #3$};
  }}}
}

\newcommand{\boxFork}[4][solid]{\vcenter{\hbox{\tikz[scale=1.5]{
    \draw (0, 0) rectangle (4ex, 4ex);
    \draw (-.2ex, 1.5ex) -- (.2ex, 1.5ex) (-.2ex, 2.5ex) -- (.2ex, 2.5ex);
    \draw (3.8ex, 1.5ex) -- (4.2ex, 1.5ex) (3.8ex, 2.5ex) -- (4.2ex, 2.5ex);
    \draw[#1] (0, 3ex) .. controls (.5ex, 3ex) and (1ex, 2.5ex).. (1ex, 2ex) .. controls (1ex, 1.5ex) and (.5ex, 1ex) .. (0ex, 1ex);
    \draw[#1] (1ex, 2ex) .. controls (1.6ex, 3ex) and (2.4ex, 1ex) .. (3ex, 2ex);
    \draw[#1] (4ex, 3ex) .. controls (3.5ex, 3ex) and (3ex, 2.5ex).. (3ex, 2ex) .. controls (3ex, 1.5ex) and (3.5ex, 1ex) .. (4ex, 1ex);
    \draw [decorate,decoration={brace,amplitude=.7ex,raise=.4ex},yshift=0pt] (0, 0) -- (0, 4ex);
    \draw (2ex, .2ex) node[anchor=north] {$\scriptstyle #2$};
    \draw (-.5ex, 2ex) node[anchor=east] {$\scriptstyle #3$};
    \draw (3.8ex, 2ex) node[anchor=west] {$\scriptstyle #4$};
  }}}
}

\newcommand{\boxSquareTickLabel}[3]{\vcenter{\hbox{\tikz[scale=1.5]{
    \draw (0, 0) rectangle (4ex, 4ex);
    \draw (3.8ex, 2ex) -- (4.2ex, 2ex);
    \draw (0, 1.7ex) .. controls (1.3ex, 4ex) and (2.6ex, 1.4ex) .. (4ex, 3.2ex);
    \draw (2ex, .2ex) node[anchor=north] {$\scriptstyle #1$};
    \draw (.4ex, 2ex) node[anchor=east] {$\scriptstyle #2$};
    \draw [decorate,decoration={brace,amplitude=.4ex,raise=.4ex},yshift=0pt] (4ex, 4ex) -- (4ex, 2ex);
    \draw (4.3ex, 3ex) node[anchor=west] {$\scriptstyle #3$};
  }\;}}
}

\newcommand{\boxSquareTwoTicks}[3]{\vcenter{\hbox{\tikz[scale=1.5]{
    \draw (0, 0) rectangle (4ex, 4ex);
    \draw (3.7ex, 1.3ex) -- (4.3ex, 1.3ex) (3.8ex, 2ex) -- (4.2ex, 2ex) (3.7ex, 2.7ex) -- (4.3ex, 2.7ex);
    \draw (0, 2ex) .. controls (1.3ex, 1.4ex) and (2.6ex, 4ex) .. (4ex, 3.4ex);
    \draw [decorate,decoration={brace,amplitude=.3ex,raise=.4ex},yshift=0pt] (4.2ex, 2.7ex) -- (4.2ex, 1.3ex);
    \draw (.4ex, 2ex) node[anchor=east] {$\scriptstyle #1$};
    \draw (2ex, .2ex) node[anchor=north] {$\scriptstyle #2$};
    \draw (4.4ex, 2ex) node[anchor=west] {$\scriptstyle #3$};
  }}}
}

\newcommand{\boxSquareTwoCrossings}[3]{\vcenter{\hbox{\tikz[scale=1.5]{
    \draw (0, 0) rectangle (4ex, 4ex);
    \draw (3.7ex, 1.3ex) -- (4.3ex, 1.3ex);
    \draw (3.8ex, 2ex) -- (4.2ex, 2ex);
    \draw (3.7ex, 2.7ex) -- (4.3ex, 2.7ex);
    \draw (0, 2ex) .. controls (1.3ex, 1ex) and (2.6ex, 2ex) .. (4ex, 2.2ex);
    \draw [\dual] (2ex, 4ex) .. controls (2ex, 3ex) and (2.3ex, 3ex) .. (4ex, 2.4ex);
    \draw [decorate,decoration={brace,amplitude=.3ex,raise=.4ex},yshift=0pt] (4.2ex, 2.7ex) -- (4.2ex, 1.3ex);
    \draw (.4ex, 2ex) node[anchor=east] {$\scriptstyle #1$};
    \draw (2ex, .2ex) node[anchor=north] {$\scriptstyle #2$};
    \draw (4.4ex, 2ex) node[anchor=west] {$\scriptstyle #3$};
  }}}
}

\newcommand{\boxSquareCorridor}[5]{\vcenter{\hbox{\tikz[scale=1.5]{
    \draw [color=gray!50!white](4.5ex, 2ex) -- (6ex, 2ex) (4.5ex, 4ex) -- (6ex, 4ex);
    \draw (0, 0) rectangle (6ex, 6ex);
    \draw (4.5ex, 0) -- (4.5ex, 6ex);
    \draw (5.7ex, 2.3ex) -- (6.3ex, 2.3ex) (5.8ex, 3ex) -- (6.2ex, 3ex) (5.7ex, 3.7ex) -- (6.3ex, 3.7ex);
    \draw (4.5ex, 1.5ex) .. controls (5ex, 1ex) and (5ex, 3.5ex) .. (6ex, 3.3ex);
    \draw [\dual] (5ex, 6ex) .. controls (6ex, 5ex) and (4ex, 3.5ex) .. (6ex, 3.5ex);
    \draw [decorate,decoration={brace,amplitude=.5ex,raise=.5ex},yshift=0pt] (6.2ex, 3.7ex) -- (6.2ex, 2.3ex);
    \draw [decorate,decoration={brace,amplitude=.5ex,raise=.4ex},yshift=0pt] (6ex, 6ex) -- (6ex, 4ex);
    \draw [decorate,decoration={brace,amplitude=.3ex,raise=.2ex},yshift=0pt] (4.5ex, 6ex) -- (6ex, 6ex);
    \draw (.4ex, 3ex) node[anchor=east] {$\scriptstyle #1$};
    \draw (3ex, .2ex) node[anchor=north] {$\scriptstyle #2$};
    \draw (6.5ex, 3ex) node[anchor=west] {$\scriptstyle #3$};
    \draw (5.2ex, 6ex) node[anchor=south] {$\scriptstyle #4$};
    \draw (6.3ex, 5ex) node[anchor=west] {$\scriptstyle #5$};
  }}}
}

\newcommand{\boxUseFork}[4]{\vcenter{\hbox{\tikz[scale=1.5]{
    \draw (0, 0) rectangle (4ex, 4ex);
    \draw (-4ex, -3.5ex) rectangle (4ex, 4.5ex) (0, -3.5ex) rectangle (8ex, 4.5ex);
    \draw (-.2ex, 1.5ex) -- (.2ex, 1.5ex) (-.2ex, 2.5ex) -- (.2ex, 2.5ex);
    \draw (3.8ex, 1.5ex) -- (4.2ex, 1.5ex) (3.8ex, 2.5ex) -- (4.2ex, 2.5ex);
    \draw (0, 3ex) .. controls (.5ex, 3ex) and (1ex, 2.5ex).. (1ex, 2ex) .. controls (1ex, 1.5ex) and (.5ex, 1ex) .. (0ex, 1ex);
    \draw (1ex, 2ex) .. controls (1.6ex, 3ex) and (2.4ex, 1ex) .. (3ex, 2ex);
    \draw (4ex, 3ex) .. controls (3.5ex, 3ex) and (3ex, 2.5ex).. (3ex, 2ex) .. controls (3ex, 1.5ex) and (3.5ex, 1ex) .. (4ex, 1ex);
    \draw (-4ex, -1ex) .. controls (0ex, -2ex) and (0ex, 1ex).. (4ex, 2ex);
    \draw (0ex, 2ex) .. controls (3ex, 4ex) and (5ex, 3ex).. (8ex, 4ex);
    \draw [decorate,decoration={brace,amplitude=.5ex,raise=.2ex},yshift=0pt] (4ex, -3.5ex) -- (0ex, -3.5ex);
    \draw [decorate,decoration={brace,amplitude=.7ex,raise=.2ex},yshift=0pt] (-4ex, -3.5ex) -- (-4ex, 4.5ex);
    \draw (-4.3ex, .5ex) node[anchor=east] {$\scriptstyle #1$};
    \draw (2ex, -3.7ex) node[anchor=north] {$\scriptstyle #2$};
    \draw [fill, color=white](4.4ex, 1.1ex) rectangle (6.2ex, 2.7ex);
    \draw [decorate,decoration={brace,amplitude=.2ex,raise=.2ex},yshift=0pt] (4.2ex, 2.5ex) -- (4.2ex, 1.5ex);
    \draw (4.2ex, 2ex) node[anchor=west] {$\scriptstyle #3$};
  }}}
}

\newcommand{\boxAnnulus}[2]{\vcenter{\hbox{\tikz[scale=1.5]{
    \draw (5ex, 3.5ex) rectangle (7ex, 5.5ex);
    \draw (3ex, 1.5ex) rectangle (9ex, 7.5ex);
    \draw [rounded corners=.25em]
      (4ex, 6.5ex) .. controls (6ex, 6.5ex + .4ex*rand) ..
      (8ex, 6.5ex) .. controls (8ex + .4ex*rand, 4.5ex) ..
      (8ex, 2.5ex) .. controls (6ex, 2.5ex + .4ex*rand) ..
      (4ex, 2.5ex) -- (4ex + .4ex*rand, 4.5ex) -- cycle;
    \draw (9.4ex, 3.5ex) -- (9.4ex, 5.5ex);
    \draw (9.2ex, 3.5ex) -- (9.6ex, 3.5ex) (9.2ex, 5.5ex) -- (9.6ex, 5.5ex);
    \draw (9.4ex, 4.5ex) node[anchor=west] {$\scriptstyle #1$};
    \draw (2.4ex, 1.5ex) -- (2.4ex, 7.5ex);
    \draw (2.2ex, 1.5ex) -- (2.6ex, 1.5ex) (2.2ex, 7.5ex) -- (2.6ex, 7.5ex);
    \draw (2ex, 4.5ex) node[anchor=east] {$\scriptstyle #2$};
  }}}
}

\newcommand{\boxRings}[5][solid]{\vcenter{\hbox{\tikz[scale=1.5]{
    \draw (0, 0) rectangle (6ex, 8ex) (6ex, 0) rectangle (12ex, 8ex);
    \draw (5.8ex, 4.2ex) -- (6.2ex, 4.2ex) (5.8ex, 4.8ex) -- (6.2ex, 4.8ex);
    \draw [#1] (0, 3ex) .. controls (2ex, 2ex) and (4ex, 5ex) .. (6ex, 4.65ex);
    \draw [#1] (6ex, 4.35ex) .. controls (8ex, 4ex) and (8ex, 7ex) .. (12ex, 6ex);
    \draw [#1,decorate,decoration={brace,amplitude=.7ex,raise=.2ex},yshift=0pt] (12ex, 0ex) -- (0, 0);
    \draw (5ex, 3.5ex) rectangle (7ex, 5.5ex);
    \draw (3ex, 1.5ex) rectangle (9ex, 7.5ex);
    \draw (3ex, 1ex) -- (9ex, 1ex);
    \draw (3ex, .8ex) -- (3ex, 1.2ex) (9ex, .8ex) -- (9ex, 1.2ex);
    \draw [fill, color=white] (4ex, .2ex) rectangle (8ex, 1.3ex);
    \draw [#1,rounded corners=.25em]
      (4ex, 6.5ex) .. controls (6ex, 6.5ex + .4ex*rand) ..
      (8ex, 6.5ex) .. controls (8ex + .4ex*rand, 4.5ex) ..
      (8ex, 2.5ex) .. controls (6ex, 2.5ex + .4ex*rand) ..
      (4ex, 2.5ex) -- (4ex + .4ex*rand, 4.5ex) -- cycle;
    \draw (6ex, -.4ex) node[anchor=north] {$\scriptstyle #2$};
    \draw (.4ex, 4ex) node[anchor=east] {$\scriptstyle #3$};
    \draw (9.4ex, 3.5ex) -- (9.4ex, 5.5ex);
    \draw (9.2ex, 3.5ex) -- (9.6ex, 3.5ex) (9.2ex, 5.5ex) -- (9.6ex, 5.5ex);
    \draw (9ex, 4.5ex) node[anchor=west] {$\scriptstyle #4$};
    \draw (6ex, 1.7ex) node[anchor=north] {$\scriptstyle #5$};
  }\;\;}}
}

\newcommand{\lro}[1][]{\overset{\tiny \mathcal{O}}\longleftrightarrow}
\newcommand{\lrv}[1][]{\overset{ \mathcal{V}}\longleftrightarrow}
\newcommand{\lambdacritical}[0]{{\lambda_{c}^\star}}
\newcommand{\lambdaoccupied}[0]{{\lambda_{c}}}

\begin{document}

\title{Sharpness of the phase transition for continuum\\percolation in $\mathbb{R}^2$}
\date{\today}

\author{Daniel Ahlberg\thanks{Instituto de Matemática Pura e Aplicada, Estrada Dona Castorina 110, 22460-320 Rio de Janeiro}
\thanks{Department of Mathematics, Uppsala University, SE-75106 Uppsala}
  \and
  Vincent Tassion\thanks{Université de Genève, 2-4 rue du Lièvre, 1211 Genève}
  \and
  Augusto Teixeira\footnotemark[1]}

\maketitle

\begin{abstract}
  We study the phase transition of random radii Poisson Boolean percolation: Around each point of a planar Poisson point process, we draw a disc of random radius, independently for each point.
  The behavior of this process is well understood when the radii are uniformly bounded from above.
  In this article, we investigate this process for unbounded (and possibly heavy tailed) radii distributions.
  Under mild assumptions on the radius distribution, we show that both the vacant and occupied sets undergo a phase transition at the same critical parameter $\lambda_c$.
  Moreover,
  \begin{itemize}
  \item For $\lambda < \lambda_c$, the vacant set has a unique unbounded connected component and we give precise bounds on the one-arm probability for the occupied set, depending on the radius distribution.
  \item At criticality, we establish the box-crossing property, implying that no unbounded component can be found, neither in the occupied nor the vacant sets.
    We provide a polynomial decay for the probability of the one-arm events, under sharp conditions on the distribution of the radius.
  \item For $\lambda > \lambda_c$, the occupied set has a unique unbounded component and we prove that the one-arm probability for the vacant decays exponentially fast.
  \end{itemize}
  The techniques we develop in this article can be applied to other models such as the Poisson Voronoi and confetti percolation.
\end{abstract}

\renewcommand\footnotemark{}
\renewcommand\footnoterule{}
\let\thefootnote\relax

\bigskip

\begin{center}
  \begin{tikzpicture}[scale=1.4]
    \clip (-.6, -.6) rectangle (4.6, 3);
    \begin{scope}
      \iffinal
      \foreach \z in {1,...,1250}
      { \pgfmathrandominteger{\x}{0}{400};
        \pgfmathrandominteger{\y}{0}{300};
        \pgfmathrandominteger{\c}{3}{7};
        \draw[color=gray,fill=gray,opacity=.3] ({.02*\x-1}, {.02*\y-1}) circle (.004 * \c * \c);
      }
      \fi
    \end{scope}
  \end{tikzpicture}
\end{center}

\bigskip


\newpage

\tableofcontents

\bigskip

\section{Introduction}

Percolation is the branch of probability theory focused on the study of the geometry and connectivity properties of random media.
Since its foundation in the 1950s, with the work of Broadbent and Hammersley~\cite{broham57}, the area reached new heights over the decades to come, specially in two dimensions: During the 1980s following Kesten's determination of the critical threshold~\cite{kesten80}, and at the turn of the century with Schramm's introduction of Schramm-Loewner evolution~\cite{schramm00} and Smirnov's proof of Cardy's formula~\cite{smirnov01}.
This progress has been well documented in a range of books on the subject; see for instance \cite{kesten1982percolation}, \cite{grimmett99}, \cite{bolrio06} and~\cite{garste}.

Bernoulli percolation on a symmetric planar lattice, e.g.\@ on the square or triangular lattices, is a cornerstone within percolation theory.
Loved for its simple yet challenging structure, Bernoulli percolation has become quintessential in the study of phase transitions and other phenomena emanating from statistical mechanics and mathematical physics.

In the site-percolation version of this model, each vertex of the lattice is independently declared `open' with probability $p$ and `closed' with probability $1-p$.
A random graph is obtained from the initial lattice by removing the closed vertices and the connected components of this graph are called \emph{clusters}.
As the parameter increases the model undergoes a sharp phase transition.
More precisely, in this planar case, it is well-known that there exists a critical parameter $p_c$, strictly in between zero and one, such that
\begin{enumerate}[\quad (i)]
\item for $p<p_c$, the probability of observing an open path from $0$ to distance $n$
  decays exponentially fast in $n$;
\item at $p=p_c$, there is no infinite open connected component, and the
  probability of an open path from $0$ to distance $n$ decays polynomially fast
  in $n$; and
\item for $p>p_c$, there exists a unique infinite open cluster.
\end{enumerate}
The phase transition is said to be `sharp' because of the abrupt change in the decay of the connection probabilities, from exponentially small in the subcritical regime to uniformly bounded in supercritical. Proofs of the above properties can be found in~\cite{grimmett99} and~\cite{bolrio06}.

In this paper we investigate the sharpness of the phase transition for Poisson
Boolean percolation in $\R^2$, establishing results analogous to the ones
described above for Bernoulli percolation. Poisson Boolean percolation is an
archetype for percolation in the continuum, and shares many features of
Bernoulli percolation while posing significant additional challenges. Apart from
being continuous rather than discrete, these challenges come from its
asymmetrical nature (the `open' and `closed' set have different properties) and
long-range dependencies.

One particular strength of the present paper is its generality, and its structure has been oriented with this in mind.
Although our results are presented for Poisson Boolean percolation, they extend to many other percolation processes in $\R^2$, such as Poisson Voronoi and confetti percolation, as will be described below.
Together with an accompanying paper~\cite{ahltastei}, our results give a precise description of the phase transition for Poisson Boolean percolation.
There, very specific properties of Poisson Boolean percolation will be exploited, as opposed to the robust methods developed here.

\bigskip

An important ingredient in establishing properties (i)--(iii) for Bernoulli percolation is what we call the dual process.
It can be defined by looking at the closed vertices on a modified graph, called the matching or dual graph.
This dual graph has a critical parameter $p_c^\star$ which was proved to satisfy the duality relation
\begin{equation}
  \label{eq:2}
  p_c + p_c^\star = 1.
\end{equation}
This relation is especially useful in the study of self-similar processes, for which the dual process at $p$ coincides with the primal process at $1-p$. In this case it can be linked with the equality $p_c = p_c^\star$, yielding $p_c = 1/2$, see more examples in Section~\ref{s:other_models}.

Kesten's~\cite{kesten80} original proof of the sharpness of the phase transition for Bernoulli percolation is based on the analysis of the crossing probabilities for rectangles; a rectangle is said to be crossed if there is a path from left to right, made out of open vertices.
The proof involves three essential ingredients:
\begin{description}
\item \emph{Finite-size criterion} -- if for some $n$, the probability to cross
  a $n$ by $3n$ rectangle in the short direction is smaller than some small
  constant $\theta>0$, then the two-point connection probability decays
  exponentially fast. On the other side, if for some $n$, the probability to
  cross a $3n$ by $n$ rectangle in the long direction is larger than
  $1-\theta>0$, then there exists an infinite cluster almost surely, see
  \cite{Rus81}, \cite{kesten1982percolation} and \cite{aizchachafrorus83}.
\item \emph{Russo-Seymour-Welsh theory} -- relates crossing probabilities for rectangle with different aspect ratios, see \cite{russo1978note} and \cite{seymour1978percolation}.
\item \emph{Threshold phenomenon} -- given a positive constant $c > 0$ and a large rectangle $K$, there is only a small interval of parameters $p$ for which the crossing probability of $K$ remains between $c$ and $1 - c$ this interval is called the critical window.
  The first proof of this sharp threshold is due to Kesten in \cite{kesten80} and is based on a geometric interpretation of the derivative of the crossing probability for a rectangle.
  A better understanding for threshold phenomena has since been obtained in work like \cite{russo82}, \cite{kahkallin88}, \cite{talagrand94} and elsewhere.
\end{description}

The strategy described above is inherently two-dimensional.
As such, the finite-size criterion extends well beyond the product structure of the Bernoulli percolation measure, and require only a minimal assumption on the decay of long-range dependencies.
The original Russo-Seymour-Welsh techniques and sharp threshold results rely in a much stronger sense on the product structure of the Bernoulli measure, and do not extend easily to more general percolation models. These will thus be the two foremost challenges in for the present study, in which we study the sharpness of the phase transition for continuum percolation in $\R^2$, which may present arbitrarily slow decay of dependencies.

\paragraph{The phase transition for Poisson Boolean percolation}
\label{sec:phase-trans-poiss}

Poisson Boolean percolation was introduced by Gilbert in~\cite{gilbert61}.
In this model we start with a Poisson point process on $\mathbb{R}^2$
with intensity parameter $\lambda > 0$.
We then independently associate to each of
these points a disc of a randomly chosen radius, according to a fixed
probability measure $\mu$ on $\R_+$. The set $\mathcal{O} \subseteq \mathbb{R}^2$ of points which are covered by at least one of the above
discs is called \emph{the occupied set}, while its
complement $\mathcal{V}:=\R^2\setminus\mathcal{O}$ is referred to as \emph{the vacant set}. We denote by $\Pr_\lambda$ the measure associated with this construction, and defer a more formal definition
to Section~\ref{s:preliminaries}.

Denote by $[0 \lro \partial B_r]$ the event that there exists an occupied path from $0$ to distance $r$.
We write $[0 \lro \infty]$ if $[0 \lro \partial B_r]$ holds for every $r \ge 1$.
Analogously we define $[0 \lrv \partial B_r]$ and $[0 \lrv \infty]$ for the corresponding events for the vacant set. Finally, define the critical parameters
\begin{equation}
  \label{eq:1}
  \begin{array}{ll}
    \lambdaoccupied := \sup \big\{ \lambda \geq 0 \; : \; \mathbb P_\lambda\big[0 \lro \infty\big] = 0 \big\},&\\
    \lambdacritical := \inf \big\{ \lambda \geq 0 \; : \; \mathbb P_\lambda\big[0 \lrv \infty\big] = 0 \big\}.&
  \end{array}
\end{equation}
We aim to describe the percolative properties of the model at and around these critical values.

That $\lambda_c$ and $\lambda_c^\star$ are finite may be obtained via a comparison with Bernoulli percolation on $\Z^2$, where a standard Peierls argument may be used.
To show, on the other hand, that $\lambda_c$ and $\lambda_c^\star$ are strictly positive is a different matter, and requires a condition on the radii distribution in order to be true.
The most fundamental condition in the study of Poisson Boolean percolation in $\R^2$ is that of finite second moment on the radius distribution:
\begin{equation}
  \label{eq:7}
  \int_0^\infty x^2 \mu(\d x) < \infty.
\end{equation}
It was observed by Hall~\cite{hall85} that~\eqref{eq:7} is necessary in order to avoid the entire plane to be almost surely covered, regardless of the intensity (as long as positive) of the Poisson point process.
Gou\'er\'e~\cite{gouere08} further showed that this condition is also sufficient for $\lambda_c$ to be strictly positive. Lower bounds on $\lambda_c^\star$ may be obtained from their comparison with $\lambda_c$. Under the assumption of bounded support on $\mu$ it is known that $\lambda_c^\star=\lambda_c$; see~\cite{roy90}. At the critical point $\lambda_c$ it is further known that there is almost surely no unbounded occupied component in the case of $\mu$ having bounded support~\cite{roy90}, whereas the analogous statement for an unbounded vacant component is only known to hold in the case of unit radii, see~\cite{alexander96}.

Our goal with the present and forthcoming paper~\cite{ahltastei} is to extend these results to hold under the condition~\eqref{eq:7}.
We will in this paper make no further assumption of the above results, as they will be easy consequences of the techniques we develop.
We merely acknowledge the observation that~\eqref{eq:7} is necessary in order for the model to present non-trivial behavior.

Following Kesten's lead, we shall in this paper focus on the study of crossings of rectangles and build a theory around them. Quantitative estimates on the rate of decay of connection probabilities, in the different regimes, will be obtained as consequences of this. The first result we state will thus relate crossing probabilities to the critical parameter $\lambda_c$.

Let us first define, for every $\x, h \ge 1$, the event $\cross(\x, h)$ that there exists an occupied path inside the rectangle $[0, \x] \times [0, h]$ from the left side to the right side. That is,
\begin{equation}
  \cross(\x, h) = \boxRectangle{\x}{\y}.
\end{equation}
We also write $\cross^\star(\x, h)$ for the existence of a vacant crossing of the same box.

\bigskip

The results of this paper are based on the following theorem.

\nc{c:thm_sharp_cross}
\begin{thm}
  \label{thm:sharpnessCrossing}
  Assume the second moment condition~\eqref{eq:7}. Then $\lambda_c$ is strictly between zero and infinity, and
  \begin{enumerate}[(i)]
  \item\label{item:3} for all $\lambda > \lambdaoccupied$ and all $\kappa > 0$,
    we have
    \begin{equation}
      \lim_{r \to \infty} \mathbb P_\lambda \big[ \cross(\kappa r, r) \big] = 1.
    \end{equation}
  \item\label{item:2} for $\lambda = \lambdaoccupied$ and all $\kappa > 0$,
    there exists a constant $\uc{c:thm_sharp_cross} = \uc{c:thm_sharp_cross}(\kappa)>0$ such that
    \begin{equation}
      \uc{c:thm_sharp_cross} < \mathbb P_{\lambdaoccupied} \big[ \cross(\kappa r, r) \big] < 1 - \uc{c:thm_sharp_cross}, \text{ for every $r \ge 1$.}
    \end{equation}
  \item\label{item:1} for all $\lambda < \lambdaoccupied$ and all $\kappa > 0$,
    \begin{equation}
      \lim_{r \to \infty} \mathbb{P}_\lambda\big[\cross(r,\kappa r)\big] = 0.
    \end{equation}
  \end{enumerate}
  Moreover, at $\lambda_c$ there is almost surely no unbounded cluster of either kind.
\end{thm}

\begin{remark}
  Our proof gives quantitative bounds for the rate of convergence in parts~\emph{(i)} and~\emph{(iii)}.
\end{remark}

The next two theorems are intended to illustrate what can be obtained using the techniques of the present work.
They give an overview of the percolative behavior of the vacant and occupied sets in Poisson Boolean percolation.
Several of these characteristics resemble the known features of Bernoulli percolation, and these two theorems are general: Similar theorems can be proved also for other percolation processes using the methods of this paper (see Section~\ref{s:other_models}).
However, their general nature calls for a slightly stronger moment condition than the one in~\eqref{eq:7}.
These two alternative hypotheses are
\begin{gather}
  \label{e:log_condition}
  \int_0^\infty x^2\log x\, \mu(\d x) < \infty, \text{ and}\\
  \label{e:alpha_condition}
  \text{for some $\alpha > 0$,} \quad \int_0^\infty x^{2 + \alpha} \,\mu(\d x) < \infty.
\end{gather}

The three hypotheses \eqref{eq:7}, \eqref{e:log_condition} and
\eqref{e:alpha_condition} will be useful because they imply some decorrelation
inequalities discussed in Section~\ref{def:rho}. The natural second moment
assumption \eqref{eq:7} already implies some spatial decorrelation properties,
and the hypotheses \eqref{e:log_condition} and \eqref{e:alpha_condition} imply
quantitative bounds on the spatial decorrelation. The hypothesis \eqref{eq:7} will
be sufficient for most of the paper, but the hypotheses~\eqref{e:log_condition}
and \eqref{e:alpha_condition} will be useful to apply some renormalization
methods presented in Section~\ref{s:finite-size}.

For the vacant set we have the following.

\nc{c:poly_decay}
\nc{c:exp_decay}
\begin{thm}
  \label{thm:sharpnessVacant}
  Assume the $2+\log$ moment condition~\eqref{e:log_condition}.
  Then $\lambda_c^\star=\lambda_c$ and is thus strictly between zero and infinity.
  Moreover,
  \begin{enumerate}[(i)]
  \item\label{item:5} for $\lambda = \lambdacritical$ there exists a constant $\uc{c:poly_decay} > 0$ such that
    \begin{equation}
      \mathbb{P}_{\lambda_c^\star}\big[ 0 \lrv \partial B_r \big] \leq \frac1{\uc{c:poly_decay}} r^{-\uc{c:poly_decay}}.
    \end{equation}
  \item\label{item:4} for all $\lambda > \lambdacritical$ there exists a constant $\uc{c:exp_decay} = \uc{c:exp_decay}(\lambda) > 0$ such that
    \begin{equation}
      \mathbb{P}_\lambda\big[ 0 \lrv \partial B_r \big] \leq \frac1{\uc{c:exp_decay}}
      \exp\{-\uc{c:exp_decay}r\}.
    \end{equation}
\end{enumerate}
\end{thm}

\begin{remark}
  Condition~\eqref{e:log_condition} in Theorem~\ref{thm:sharpnessVacant} is not sharp.
  In the forthcoming paper \cite{ahltastei}, we show that the second moment condition suffices to prove all of the above.
  The fact that $\lambdacritical = \lambdaoccupied$ is analogous to $p_c + p_c^\star = 1$ for Bernoulli percolation.
\end{remark}

There are important distinctions between the vacant and the occupied sets in
regard to their percolation properties. An important difference comes from the
fact that parts~\eqref{item:5} and \eqref{item:4} of
Theorem~\ref{thm:sharpnessVacant} do not hold in general for the occupied set,
since the existence of long occupied connections can be triggered by a single
disk of large radius. For example, if one chooses the radius distribution $\mu$
such that $\mu[x, \infty) = x^{-2}(\log x)^{-2}$ for $x$ large enough, then the
hypothesis \eqref{eq:7} is satisfied but for every $\lambda > 0$ and every $r$
large enough
  \begin{equation}
    \Pr_\lambda\big[0 \lro \partial B_r\big] \,\ge\, c \cdot \lambda \int_r^\infty x \,\mu[x,\infty)\d
    x \,\ge\, c \cdot \lambda/\log r.
  \end{equation}
  Therefore, in this case the one-arm probability $\Pr_\lambda[0 \lro \partial B_r]$ cannot decay polynomially fast at criticality, just because this choice of $\mu$ implies that the probability that $0$ is contained in a ball of radius larger than $r$ does not decay polynomially.  Nevertheless, under the stronger moment assumption \eqref{e:alpha_condition}, one can prove the polynomial decay of the arm exponent as shown in the next result.

\begin{thm}
  \label{thm:sharpnessOccupied}
  Assume that the $2+\alpha$ moment condition~\eqref{e:alpha_condition} holds
  for some $\alpha>0$. Then,
  \begin{enumerate}[(i)]
  \item\label{item:5_o} \nc{c:poly_o_decay} for $\lambda = \lambdaoccupied$, there exists a constant $\uc{c:poly_o_decay} = \uc{c:poly_o_decay}(\alpha) > 0$ such that
    \begin{equation}
      \mathbb{P}_{\lambda_c}\big[ 0 \lro \partial B_r \big] \leq \frac1{\uc{c:poly_o_decay}} r^{-\uc{c:poly_o_decay}}.
    \end{equation}
  \item\label{item:4_o} for all $\lambda < \lambdaoccupied$, there exists a constant $\uc{c:poly_decay}=\uc{c:poly_decay}(\lambda) > 0$ such that
    \begin{equation}
      \mathbb{P}_\lambda\big[ 0 \lro \partial B_r \big] \leq \uc{c:poly_decay} r^{-\alpha}.
    \end{equation}
\end{enumerate}
\end{thm}

\begin{remark}
The exponent in part~\emph{(ii)} of Theorem~\ref{thm:sharpnessOccupied} cannot be improved in general. To see this, consider the radii distribution for which $\mu[x,\infty)=x^{-(2+\alpha)}(\log x)^{-2}$ for large $x$. This distribution satisfies~\eqref{e:alpha_condition}, but
$$
\Pr_\lambda\big[ 0 \lro \partial B_r \big]\,\ge\,\lambda r^2\mu[2r,\infty)\,=\,\lambda (2r)^{-\alpha}(\log 2r)^{-2}
$$
for large $r$.
\end{remark}

We will now describe the main steps in the proof of the above results, which roughly speaking correspond to the three steps described for the case of Bernoulli percolation above.
These are \emph{finite-size criterion}, \emph{Russo-Seymour-Welsh theory} and a \emph{threshold phenomenon}.
We describe each of these steps below and emphasize that their corresponding proofs can be read independently of one another.

\paragraph{Finite-size criterion}

In Section~\ref{s:finite-size} we prove that the crossing probabilities converge to 1 as soon as they become close enough to 1.
Roughly speaking we prove that there exists $\theta>0$ and $r_0(\lambda)$ such that the following are equivalent
\begin{gather}
  \mathbb{P}_\lambda \big[\cross(3\x, \x)\big] > 1 - \theta \text{ for some $\x \geq r_0$,}\\
  \lim_{r \to \infty} \mathbb{P}_\lambda \big[\cross(3\x, \x)\big] = 1,
\end{gather}
and the analogous statement holds for vacant crossings, see Proposition~\ref{prop:finite-size}.
This motivates us to introduce
\begin{gather}
  \label{eq:11}
  \lambda_0 := \sup\big\{\lambda \ge 0 : \lim_{\x \to \infty} \Pr_\lambda \big[ \cross(\x, 3\x) \big] = 0\big\},\\
  \label{eq:9}
  \lambda_1 := \inf\big\{\lambda \ge 0 : \lim_{\x \to \infty} \Pr_\lambda \big[ \cross(3\x, \x) \big] = 1\big\}.
\end{gather}
The first important consequence of the finite-size criterion is that
$$
0\,<\,\lambda_0\,\le\,\lambda_c\,\le\, \lambda_1\, < \,\infty.
$$
These results are proved in Section~\ref{s:finite-size}.

Up to now, it is not clear that $\lambda_0$ is actually equal to $\lambda_1$.
We call $[\lambda_0, \lambda_1]$ the \emph{critical regime} and we can see that, for every $\lambda \in [\lambda_0, \lambda_1]$,
\begin{equation}
  \label{e:state_easy_cross}
  \inf_{\x \ge 1} \Pr_\lambda \big[ \cross(\x, 3\x) \big] > 0 \quad \text{and} \quad \inf_{\x \ge 1} \Pr_\lambda \big[ \cross^\star(\x, 3\x) \big] > 0.
\end{equation}
These statements say that in the critical regime, rectangles have non-degenerate probabilities of being crossed by both the vacant and the occupied set.
However, one should notice that these crossings happen in the easy direction of the rectangles.
In order to obtain a more precise description of the behavior in the critical regime we require similar statements regarding crossing probabilities in the long direction. This is precisely the purpose of Russo-Seymour-Welsh theory, as we describe next.

\paragraph{Russo-Seymour-Welsh theory}

In this step we show that \eqref{e:state_easy_cross} implies that
\begin{equation}
  \label{e:state_hard_cross}
  \inf_{\x \ge 1} \Pr_\lambda \big[ \cross(3 \x, \x) \big] > 0 \quad \text{and} \quad \inf_{\x \ge 1} \Pr_\lambda \big[ \cross^\star(3 \x, \x) \big] > 0.
\end{equation}
Note that the only difference between the above and \eqref{e:state_easy_cross}
is that now the rectangles are crossed in the long direction. This `box-crossing' property for the critical regime rules out the existence of unbounded clusters and provides bounds on arm probabilities. These results are
proved in Section~\ref{s:critical}.

\begin{remark}
  RSW bounds were obtained separately for the occupied
  and vacant regimes during the 1990s by Roy~\cite{roy90} and
  Alexander~\cite{alexander96}. However, these proofs are technical, and pose
  pose heavy restrictions on the radii distribution. Recently,
  Tassion~\cite{tassion14+} found an argument allowing for greater generality.
  His argument, presented in the setting of Voronoi percolation but extends
  verbatim, shows that
  \begin{equation}
   \label{eq:4} \liminf_{\x\to\infty}\Pr_\lambda\big[\cross(r,r)\big]>0\quad\Rightarrow\quad\liminf_{\x\to\infty}\Pr_\lambda\big[\cross(3\x,\x)\big]>0.
  \end{equation}
  While this is all that is needed in `symmetric' percolation models, such as Poisson Voronoi percolation, it is indispensable for `non-symmetric' models, which we consider here, that crossings of rectangles the easy way imply crossings also the hard way.
  This extension is not straightforward.
\end{remark}

\paragraph{Sharp thresholds}

The final step of the proof is to show that $\lambda_0 = \lambda_1$.
The main strategy here is to show that the occupied crossing probabilities grow very fast from zero to one as we increase the density $\lambda$ of the system.

There is a solid theory that predicts the occurrence of threshold phenomena in the setting of Boolean function on the discrete cube $\{0, 1\}^n$ equipped with a product structure.
Poisson point processes is the natural analogue to the discrete Bernoulli measure.
It is therefore natural to believe that these techniques should carry over also to the Poissonian continuum setting.
This is indeed the case, although doing so is not straightforward task.
We will here follow an approach introduced in~\cite{ahlbrogrimor14}.

This part of the argument involves the analysis of Boolean functions and is the least general part of our argument, since it involves representing the crossing probabilities in our Poisson point process as a certain function of discrete Bernoulli random variables. 

This construction and the proof that $\lambda_0 = \lambda_1$ can be found in Section~\ref{s:sharp}.

\paragraph{Other examples}

We have chosen Poisson Boolean percolation as the main example to illustrate the techniques of this work.
However, the techniques that we develop apply in greater generality.
We emphasize this fact in Section~\ref{s:other_models}, where we apply these techniques to Poisson Voronoi and confetti percolation. We direct the reader to that section for a statement of the results we obtain for these models, and give here just a short description of the scope of our techniques.

The finite-size criterion described above only uses the fact that the law of our percolation process is invariant under translations and right-angle rotations, see \eqref{e:invariance}, and that it satisfies a mixing condition, see \eqref{e:rho_to_zero}.

The Russo-Seymour-Welsh part of our argument, we solely use the invariance of $\mathcal{O}$ and $\mathcal{V}$ under translations, right-angle rotations and reflections in coordinate axes, see~\eqref{e:invariance}, the mixing property~\eqref{e:rho_to_zero}, the FKG inequality~\eqref{eq:3} and a certain continuity property of the crossing probabilities stated in Lemma~\ref{l:cross_continuous}.

The threshold argument is the least general and uses, in addition to the assumptions made above, the fact that process is based on a Poisson point process.

\paragraph{Some previous works on the model}

The reference book \cite{meeroy96} provides a general exposition of continuum percolation.
The case of uniformly bounded radius distribution has been extensively studied in  \cite{MR782099}, \cite{MR783056}, \cite{MR914958}, \cite{roy90}, \cite{MR1273055}, \cite{alexander96}.
These works have already established much of our results in this restricted setting.

In Lemma~A.2 of \cite{MR3225978}, it is proved that for sufficiently small $\lambda > 0$, $\mathbb{P}_\lambda[B_r \leftrightarrow B_{2r}^c]$ goes to zero with $r$ under the second moment condition \eqref{eq:7}.
This implies that $\lambda_0 > 0$ under the second moment condition.
In \cite{zbMATH04197125}, Roy studied a Poisson soup of bounded sticks in the plane, proving a Russo-Seymour-Welsh theorem and establishing that $\lambdacritical = \lambdaoccupied$.
Uniqueness of the unbounded components (for both vacant and occupied regions) has been established in~\cite{zbMATH00681409}.

\paragraph{Organization of the paper}

In Section~\ref{s:preliminaries} we provide some notation and preliminary results that are needed throughout the text.
Sections~\ref{s:finite-size}, \ref{s:critical} and~\ref{s:sharp} present the three steps of the proof that were described in the introduction, namely: the finite-size criterion, Russo-Seymour-Welsh result and the sharp threshold step.
These sections can be read independently of one another.
Section~\ref{s:continuity} proves the continuity of the critical parameter with respect to the law $\mu$ of the radius distribution.
The main results of the paper are then proved in Section~\ref{s:main_results}.
We treat other models and provide some open questions in Sections~\ref{s:other_models} and \ref{s:open}.

\begin{acknow}
  This work began during a visit of V.T. to IMPA, that he thanks for support and hospitality.
  We thank the Centre Intradisciplinaire Bernoulli (CIB) and Stardû for hosting the authors.
  D.A. was during the course of this project financed by grant 637-2013-7302 from the Swedish Research Council.
  A.T. is grateful to CNPq for its financial contribution to this work through
  the grants 306348/2012-8, 478577/2012-5 and 309356/2015-6 and FAPERJ through grant number 202.231/2015.
  V.T. acknowledges support from the Swiss NSF.
\end{acknow}

\section{Notation and preliminary results}
\label{s:preliminaries}

Throughout the text we let $c$ denote positive constants which may depend on the radius distribution and may change from line to line.
However, numbered constants such as $c_0, c_1, \dots$ refer to their first appearance in the text.
We will further write $B^\infty(r) = [-r, r]^2$ for the ball centered at the origin in the supremum norm and let $B(x,r)$, on the other hand, denote the closed Euclidean ball with center $x$ and radius $r$.
When $x$ is the origin we omit it from the above notation.

Rather informally, a realization of Poisson Boolean percolation is obtained by decorating the points in a Poisson point process in $\R^2$ by Euclidean discs with independent radii sampled from some distribution.
There are various (equivalent) ways of making this description formal. We will start this section by describing one way which will be suitable for our purposes.

\subsection{Definition of the process}

To define our process, let us first introduce a Poisson point process on the following space of point measures
\begin{equation}
  \label{e:Omega}
  \Omega = \Big\{ \omega = \sum_i \delta_{(x_i,z_i)}: (x_i,z_i) \in \mathbb{R}^2 \times \mathbb{R}_+ \text{ and } \omega \big( K \times \mathbb{R}_+ \big) < \infty \text{ for all $K$ compact} \Big\}.
\end{equation}
We endow this space with the $\sigma$-algebra $\mathcal{M}$ generated by the evaluation maps $A \mapsto \omega(A)$, for $A \in \mathcal{B}(\mathbb{R}^2 \times \mathbb{R}_+)$, the Borel sets on $\R^2 \times \mathbb{R}_+$.

We next fix an intensity parameter $\lambda\ge0$ and some probability measure $\mu$ on $\mathbb{R}_+$ which will give the radius distribution of our discs.
We can now define on $(\Omega, \mathcal{M})$ a Poisson point process with intensity $\lambda \cdot dx \, \mu(dz)$, i.e., Lebesgue measure on $\mathbb{R}^2$ product with $\mu$ and multiplied by $\lambda \geq 0$.
The law of this process is denoted by $\mathbb{P}_\lambda$ throughout the text and we complete the $\sigma$-algebra $\mathcal{M}$ with respect to $\mathbb{P}_\lambda$.

For each point $(x,z) \in \mathbb{R}^2 \times \mathbb{R}_+$ in the support of the measure $\omega \in \Omega$, we associate the disc $B(x, z)$ and the \emph{occupied} region of the plane is consequently given by
\begin{equation}
  \label{e:def_D}
  \mathcal{O} := \bigcup_{(x,z) \in \supp(\omega)} B(x, z),
\end{equation}
while the \emph{vacant} set is given by its complement $\mathcal{V} := \mathbb{R}^2 \setminus \mathcal{O}$. (Below, we will often identify $\omega$ with its support in order to ease the notation.)

From the definition of the model, it is trivial to conclude that
\begin{display}
  \label{e:invariance}
  the law of $\mathcal{O}$ is invariant under translations,\\ right-angle
  rotations, and reflection in coordinate axes.
\end{display}
Although the law is invariant under arbitrary rotations, and reflections in axes with other orientation, we will only use the above statement throughout the text.

There is a natural partial ordering of elements in $\Omega$, namely, $\omega\leq\omega'$ if $\omega(A) \leq \omega'(A)$ for all $A \in \mathcal{B}(\mathbb{R}^2 \times \mathbb{R}_+)$.
An event $A\in\mathcal{M}$ is said to be \emph{increasing} if $\omega\in A$ implies that $\omega'\in A$ for all $\omega\leq\omega'$.
It is \emph{decreasing} if its complement is increasing.

A useful property of increasing events is that they are positively correlated.
The following proposition, known as the FKG inequality, was proved by Roy in his
doctorate thesis; see also~\cite[Theorem~2.2]{meeroy96}: If $A_1$ and $A_2$ are
increasing events, then
\begin{equation}
  \label{eq:3}
  \Pr_\lambda(A_1\cap A_2)\,\geq\,\Pr_\lambda(A_1)\Pr_\lambda(A_2).
\end{equation}
The above also holds when $A_1$ and $A_2$ are decreasing events.

We will also use the following standard consequence of the FKG inequality,
referred to as the square-root trick. Let $A_1,\ldots,A_k$ be $k$ increasing events (or
$k$ decreasing events), then
\begin{equation}
  \label{eq:32}
  \max_{1\le i \le k}\Pr_\lambda(A_i)\,\ge\, 1-\big [1-\Pr_\lambda(A_1\cup\ldots\cup A_k)\big]^{1/k}.
\end{equation}

\subsection{Crossing events}

Throughout the text we will often deal with crossing events of various types. In particular, we will be interested in the following general definition.
Given subsets $A_1, A_2$ and $\mathcal{C}$ of $\mathbb{R}^2$ let
\begin{equation}
  \label{e:path_in_C}
  A_1 \overset{\mathcal{C}}\longleftrightarrow A_2 := \text{there is a path in $\mathcal{C}$ connecting $A_1$ to $A_2$}.
\end{equation}
We are now in position to give a formal definition of the crossing event $\cross(\x,h)$, for $\x, h > 0$, as
\begin{equation}
  \label{e:def_cross}
  \cross(\x, h) := \big[ A_1 {\overset{\mathcal{O} \cap K}\longleftrightarrow} A_2 \big],
\end{equation}
where $K$ denotes the box $[0,\x] \times [0,h]$ and $A_1 = \{0\} \times [0,h]$, $A_2 = \{\x\} \times [0,h]$ stand for its left and right sides.
We also define the event corresponding to a vacant crossing as $\cross^\star(\x,h) := [A_1 {\overset{\mathcal{V} \cap K}\longleftrightarrow} A_2 ]$.
It is rather straightforward to verify that events of this type are measurable.
Notice that occupied crossing events are increasing, and vacant crossing events are decreasing.

\subsection{Decay of spatial correlations}

The second moment condition given in \eqref{eq:7} is sufficient to imply a spatial decorrelation, in the sense of the function $\rho_\lambda$ introduced in Definition~\ref{def:rho} below.
First set $Y_v = \ind_{\{v \in \mathcal{V}\}}$ for $v \in \mathbb{R}^2$.
\begin{definition}
  \label{def:rho}
  Given $0 < r, s < \infty$, let
  \begin{equation}
    \rho_\lambda (r, s) := \sup_{f_1, f_2} \Cov\big(f_1(\mathcal O), f_2(\mathcal O)\big),
  \end{equation}
  where the above suppremum is taken over all functions
  $f_1,f_2:\mathcal P(\mathbb R^2)\to [-1,1]  $ such that $f_1(\mathcal O) \in \sigma(Y_v; v
  \in B^\infty(r))$ and $f_2(\mathcal O) \in \sigma(Y_v; v \not \in B^\infty(r+s))$.
\end{definition}

The function $\rho_\lambda$ has a nice geometric interpretation.
Namely, one can observe that $\rho_\lambda (r, r+s)$ is directly related to the probability that there exists one big occupied disk crossing the annulus $B^\infty(r+s)\setminus B^\infty(r)$.
By making this observation rigorous, and computing this crossing probability (see Lemma~\ref{l:decouple}), we prove the following upper bound.

\nc{c:decorrelation}
\begin{prop}
  \label{c:decouplef1f2}
  For any $\lambda > 0$ and $r, s \ge 1$, we have
  \begin{equation}
    \rho_\lambda(r, s)\, \leq\, \uc{c:decorrelation} \lambda \left(1+\frac{r}{s}\right)^2\int_{s/2}^\infty x^2\mu(\d x).
  \end{equation}
  In particular, the second moment condition~\eqref{eq:7} implies that
  \begin{equation}
    \label{e:rho_to_zero}
    \text{for any $\kappa > 0$, }\lim_{\x \to \infty} \rho_\lambda(\x, \kappa \x) = 0.
  \end{equation}
  and the $2 + \alpha$ moment condition in \eqref{e:alpha_condition} implies that for every $\lambda > 0$ and $\eps>0$
  \begin{equation}
    \label{e:ralpha_decay}
    \sup\limits_{\x \geq 1, \kappa \geq \eps} r^\alpha \rho_\lambda(\x, \kappa\x) < \infty.
  \end{equation}
\end{prop}

For $K \subseteq \mathbb{R}^2$, we write
\begin{equation}
  \mathcal{O}_K = \bigcup_{(x,z) \in\omega:\, x \in K} B(x, z).
\end{equation}
The next lemma says that deleting balls from $\mathcal{O}$ which are far from a given box do not alter the configuration of $\mathcal{O}$ inside that box.

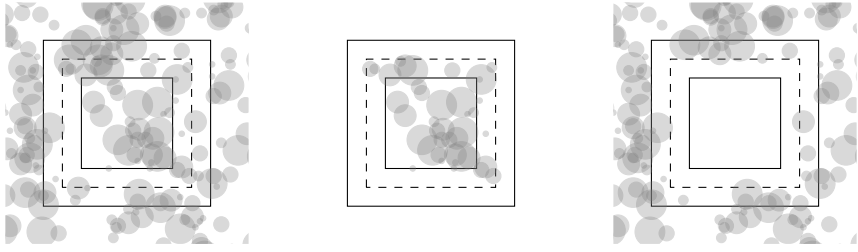
\begin{figure}[ht]
  \label{f:D_B}
  \centering
  \begin{tikzpicture}[scale=1]
    \clip (-.6, -.6) rectangle (2.6, 2.6) (3.4, -.6) rectangle (6.6, 2.6) (7.4, -.6) rectangle (10.6, 2.6);
    \foreach \w in {0,...,2} {
      \draw (0.4 + 4*\w, 0.4) rectangle (1.6 + 4*\w, 1.6);
      \draw[dashed] (0.15 + 4*\w, 0.15) rectangle (1.85 + 4*\w, 1.85);
      \draw (-.1 + 4*\w, -.1) rectangle (2.1 + 4*\w, 2.1); }
    \begin{scope}
      \iffinal
      \foreach \z in {1,...,250}
      { \pgfmathrandominteger{\x}{0}{200};
        \pgfmathrandominteger{\y}{0}{200};
        \pgfmathrandominteger{\c}{3}{7};
          \draw[color=gray,fill=gray,opacity=.3] ({.02*\x-1}, {.02*\y-1}) circle (.004 * \c * \c);
        \ifthenelse{\cnttest{\x}>{57}\AND\cnttest{\x}<{142}\AND\cnttest{\y}>{57}\AND\cnttest{\y}<{142}}{
          \draw[color=gray,fill=gray,opacity=.3] ({.02*\x+3}, {.02*\y-1}) circle (.004 * \c * \c);}{}
        \ifthenelse{\cnttest{\x}<{57}\OR\cnttest{\x}>{142}\OR\cnttest{\y}<{57}\OR\cnttest{\y}>{142}}{
          \draw[color=gray,fill=gray,opacity=.3] ({.02*\x+7}, {.02*\y-1}) circle (.004 * \c * \c);}{}
      }
      \fi
    \end{scope}
  \end{tikzpicture}
  \caption{An illustration of the sets $\mathcal{O}$, $\mathcal{O}_{B^\infty(r + s/2)}$ and $\mathcal{O}_{B^\infty(r + s/2)^c}$ respectively. Note that in this particular realization we have $\mathcal{O}_{B^\infty(r + s/2)} \cap B^\infty(r) = \mathcal{O} \cap B^\infty(r)$ and also that $\mathcal{O}_{B^\infty(r + s/2)^c} \cap B^\infty(r + s)^c = \mathcal{O} \cap B^\infty(r + s)^c$ as in Lemma~\ref{l:decouple}.}
\end{figure}

\nc{c:separate_boxes}
\begin{lma}
  \label{l:decouple}
  For all $r, s > 0$, writing $K_1 = B^\infty(r)$ and $K_2 = B^\infty(r + s)$, we have
  \begin{equation}
    \label{e:decouple1}
    \mathbb{P}_\lambda \big( \mathcal{O}_{K_2^c} \cap K_1 \neq \varnothing \big)\leq 8\lambda\left(1+\frac{r+1}{s}\right)^2\int_s^\infty x^2\mu(\d x).
  \end{equation}
  Also,
  \begin{equation}
    \label{e:decouple2}
    \mathbb{P}_\lambda \big( \mathcal{O}_{K_1} \cap K_2^c \neq \varnothing \big)\leq 8\lambda\left(\frac{r}{s}\right)^2\int_s^\infty x^2\mu(\d x).
  \end{equation}
\end{lma}
\begin{remark}
  The above lemma implies in particular that, whenever $\int_0^\infty x^2\mu(\d x) < \infty$,
  \begin{display}
    \label{e:finite_in_compact}
    if $K \subseteq \mathbb{R}^2$ is compact, then $\mathbb{P}_\lambda$-almost surely there\\
    are only finitely many balls in $\mathcal{O}$ that intersect $K$.
  \end{display}
\end{remark}

\begin{proof}
We first derive the second bound. Since $K_1$ has area $8r^2$, an immediate estimate of the expected number of points in $K_1$ with radius at least $s$ gives that
  \begin{equation*}
      \mathbb{P}_\lambda \big( \mathcal{O}_{K_1} \cap K_2^c \neq \varnothing \big) \,\leq\, 8r^2\lambda\int_s^\infty\mu(\d x) \,\leq\, 8\lambda\left(\frac{r}{s}\right)^2\int_s^\infty x^2\mu(\d x),
  \end{equation*}
  hence establishing~\eqref{e:decouple2}.

  For the first bound, we split the complement of $K_2 =B^\infty(r+s)$ into the disjoint annuli $A_i = B^\infty(r + s + i + 1) \setminus B^\infty(r + s + i)$ and write $K_2^c = \bigcup_{i\ge0} A_i$.
  Using the fact that the area of $A_i$ is bounded by $8(r+s+i+1)$, we can now estimate
  \begin{equation*}
    \begin{split}
      \mathbb{P}_\lambda \big( \mathcal{O}_{K_2^c} \cap K_1 \neq \varnothing \big)& \,\leq\, \sum_{i\ge0} \mathbb{P}_\lambda \big( \mathcal{O}_{A_i} \cap K_1 \neq \varnothing \big)\\
      & \leq\, 8\lambda\sum_{i\ge0} (r + s + i+1)\int_{s+i}^\infty\mu(\d x).
    \end{split}
  \end{equation*}
  Exchanging the order of summation yields the bound
  $$
  8\lambda\int_s^\infty\sum_{i=0}^{x-s}(r+s+i+1)\mu(\d x) \,\leq\, 8\lambda\int_s^\infty (r+x+1)^2\mu(\d x),
  $$
  from which~\eqref{e:decouple1} follows.
\end{proof}

Let us now prove Proposition~\ref{c:decouplef1f2}.

\begin{proof}[Proof of Proposition~\ref{c:decouplef1f2}]
  Recall that $\lVert f_1 \rVert_\infty, \lVert f_2 \rVert_\infty \leq 1$. Let
  $K=B^\infty(r+s/2)$. By the Poissonian character of the law
  $\mathbb{P}_\lambda$, the two random variables $ f_1(\mathcal{O}_K)$ and
  $f_2(\mathcal{O}_{K^c})$ are independent. Therefore,
  \begin{equation}
    \begin{split}
      \Cov(f_1(\mathcal O), f_2(\mathcal O) & \,\leq\, \mathbb{P}_\lambda
      \big(f_1(\mathcal O) \neq f_1(\mathcal O_K)\big) + \mathbb{P}_\lambda \big(f_2(\mathcal
      O)\neq f_2(\mathcal O_{K^c})\big)\\
      & \leq\,     \mathbb{P}_\lambda \big( \mathcal{O}_{K^c} \cap B^\infty(r)\neq\varnothing \big) + \mathbb{P}_\lambda \big( \mathcal{O}_K \cap B^\infty(r+s)^c\neq\varnothing \big),
    \end{split}
  \end{equation}
  and the proof then follows from Lemma~\ref{l:decouple}.
  For the last conclusion \eqref{e:ralpha_decay}, one may use Markov's inequality.
\end{proof}

\subsection{An alternative notion of decoupling}

Note that Definition~\ref{def:rho} is symmetric for $\mathcal O$ and $\mathcal
V$, in the sense that the roles of the occupied and vacant sets and are
interchangeable. In what follows, we will introduce another notion of decoupling
that will not be symmetric for $\mathcal{O}$ and $\mathcal{V}$. This distinction
will be very relevant to explain the different behavior of these sets, see
Remark~\ref{r:D_is_bad}, and explains why we obtain Theorem~\ref{thm:sharpnessCrossing} under~\eqref{eq:7} but require~\eqref{e:log_condition} for Theorem~\ref{thm:sharpnessVacant} with the techniques used here.

To motivate this new notion, observe from Definition~\ref{def:rho} that
\begin{equation}
  \label{e:rho_bounds_f1f2}
  \mathbb{E}_\lambda\big(f_1(\mathcal O)  f_2(\mathcal O)\big) \,\leq\,
  \mathbb{E}_\lambda\big(f_1(\mathcal O)\big) \mathbb{E}_\lambda\big(f_2(\mathcal O)\big) + \rho(r, s),
\end{equation}
for every functions $f_1,f_2:\mathcal P(\mathbb R^2)\to [0,1]$ such that   $f_1(\mathcal O) \in
\sigma(Y_v; v \in B^\infty(r))$ and $f_2(\mathcal O) \in \sigma(Y_v; v \not \in B^\infty(r + s))$.
This will be used in several parts of the text.

However, some results will be stronger using the following modified version of the above.
\begin{definition}
  \label{d:rho_bar}
  Let $\bar{\rho}_\lambda:\mathbb{R}_+^2 \to \mathbb{R}_+$ be defined as the smallest value such that
  \begin{equation}
    \label{e:rho_bar_bounds_f1f2}
    \mathbb{E}_\lambda\big(f_1(\mathcal O) f_2(\mathcal O)\big) \,\leq\,
    \mathbb{E}_\lambda\big(f_1(\mathcal O)\big) \Big( \mathbb{E}_\lambda\big(f_2(\mathcal O)\big) + \bar{\rho}_\lambda(r, s) \Big),
  \end{equation}
  for all \emph{decreasing} functions $f_1, f_2:\mathcal P(\mathbb R^2)\to
  [0,1]$ satisfying $f_1(\mathcal O) \in \sigma(Y_v; v \in B^\infty(r))$,
  $f_2(\mathcal O) \in \sigma(Y_v; v \not \in B^\infty(r + s))$. (A function
  $f:\mathcal P(\mathbb R^2) \to[0,1]$ is said to be decreasing if $f(A)\ge f(B)$ whenever $A\subset B$.)
\end{definition}

\begin{remark}
  Note that the error $\bar{\rho}_\lambda$ in the above definition is being multiplied by $\mathbb{E}_\lambda(f_1)$.
  This represents an important improvement in some cases, as will become clear later for instance when we compare \eqref{eq:19} with \eqref{eq:16}.
  Although the error term in Definition~\ref{d:rho_bar} is smaller than the one in \eqref{e:rho_bounds_f1f2}, the bound \eqref{e:rho_bar_bounds_f1f2} can only be used when $f_1$ and $f_2$ are decreasing functions.
  This restriction was intentionally introduced in the definition of $\bar{\rho}_\lambda$ and it is necessary in our proof that $\bar{\rho}_\lambda$ vanishes, see Proposition~\ref{p:bar_rho} below.
  This restriction reflects the distinct behavior that the vacant and occupied sets display at criticality, see Remark~\ref{r:D_is_bad}.
\end{remark}

We now prove an analogue of Proposition~\ref{c:decouplef1f2}.
\begin{prop}
  \label{p:bar_rho}
  For any $\lambda > 0$ and $r, s \geq 1$ we have
  \begin{equation}
    \bar{\rho}_\lambda(r, s) \,\leq\, \uc{c:decorrelation} \lambda \left(1+\frac{r}{s}\right)^2\int_{s/2}^\infty x^2\mu(\d x).
  \end{equation}
  In particular, if $\int_0^\infty x^2 \mu (\d x) < \infty$, then
  \begin{equation}
    \label{e:rho_bar_to_zero}
    \text{for any } \kappa > 0,\, \lim_{\x \to \infty} \bar{\rho}_\lambda(\x, \kappa \x) = 0.
  \end{equation}
\end{prop}

\begin{proof}
  We write $K_1 = B^\infty(r)$ and $K_2 = B^\infty(r + s)^c$ and define
  \begin{equation}
    \check{\mathcal{O}} = \bigcup_{\substack{(x,z) \in \supp(\omega); \\ B(x, z) \cap K_1 = \varnothing}} B(x, z).
  \end{equation}
  Note that this set is contained in $\mathcal{O}$ and independent of $f_1(\mathcal{O})$.

  Using the fact that $f_2$ is decreasing, we can write
  \begin{equation*}
    \begin{split}
      \mathbb{E}_\lambda \big( f_1(\mathcal{O}) f_2(\mathcal{O}) \big) & \,\leq\, \mathbb{E}_\lambda \big( f_1(\mathcal{O}) f_2(\check{\mathcal{O}}) \big) \,\leq\, \mathbb{E}_\lambda \big( f_1(\mathcal{O}) \big) \mathbb{E}_\lambda \big( f_2(\check{\mathcal{O}}) \big)\\
      & \leq\, \mathbb{E}_\lambda \big( f_1(\mathcal{O}) \big) \Big( \mathbb{E}_\lambda \big( f_2(\mathcal{O}) \big) + \mathbb{P}_\lambda \big( \check{\mathcal{O}} \cap K_2 \neq \mathcal{O} \cap K_2 \big) \Big),
    \end{split}
  \end{equation*}
  where we have used the fact that $f_1,f_2\in[0,1]$.
  In view of the definition of $\bar{\rho}_\lambda$, all we need to do is to bound the
  probability above,
  \begin{equation}
    \begin{split}
      \mathbb{P}_\lambda (\check{\mathcal{O}} & \cap K_2 \neq \mathcal{O} \cap K_2) \,\leq\, \mathbb{P}_\lambda \Big(
      \begin{array}{c}
        \text{there is $(x,z) \in \supp(\omega)$ such that}\\
        \text{$B^\infty(x,z)$ touches $K_1$ and $K_2$}
      \end{array}
      \Big)\\
      & \leq\, \mathbb{P}_\lambda \big(\mathcal{O}_{B^\infty(r + s/2)} \cap K_2 \neq \varnothing \big) + \mathbb{P}_\lambda \big( \mathcal{O}_{B^\infty(r + s/2)^2} \cap K_1 \neq \varnothing \big).
    \end{split}
  \end{equation}
  The proof now ends in the same way as the proof of Proposition~\ref{c:decouplef1f2}.
\end{proof}

\subsection{Local continuity}

Another rather straightforward consequence of the Poissonian nature of the model is the following `local continuity' property.
In essence this property says that locally the topological properties of the set $\mathcal{O}$ are unlikely to be affected by slight change in the radii of the discs.

Given a set $A\subseteq\R^2$ we define its $\eps$-interior and $\eps$-closure, for $\eps>0$, as follows:
\begin{equation*}
\begin{aligned}
  \interior(A,\eps)&:=\{x\in\R^2:B(x,\eps)\subseteq A\},\\
  \closure(A,\eps)&:=\{x\in\R^2:B(x,\eps)\cap A\neq\varnothing\}.
\end{aligned}
\end{equation*}

We omit the proof of the following proposition.

\begin{prop}\label{prop:perturbation}
Let $K \subseteq \R^2$ be compact and convex and $A_1, A_2 \subseteq K$ closed.
For every $\lambda > 0$
\begin{equation}
\lim_{\eps\to0}\Pr_\lambda \Big( A_1 {\overset{\closure(\mathcal{O},\eps) \cap K}\longleftrightarrow} A_2\text{ occurs but not }A_1 {\overset{\interior(\mathcal{O},\eps) \cap K}\longleftrightarrow} A_2 \Big) \,=\,0.
\end{equation}
\end{prop}

A consequence of the above proposition that will be used later is the following lemma.
Consider a variation of the standard crossing event defined as
\begin{equation}
  \label{e:cross_below_y}
  \boxSquareTickLabel{\x}{h}{y} := \big[ A_1 {\overset{\mathcal{O} \cap K}\longleftrightarrow} A_2 \big],
\end{equation}
where $K$ is again the box $[0,\x] \times [0,h]$, while $A_1 = \{0\} \times [0,h]$ and $A_2 = \{\x\} \times [y, h]$.

\begin{lma}
  \label{l:cross_continuous}
  Given $\x, h$ and $\lambda > 0$, the function $F:[0,h] \to [0,1]$ given by
  \begin{equation}
    F(y) = \mathbb{P}_\lambda \Bigg[ \smash{\boxSquareTickLabel{\x}{h}{y}} \Bigg] \text{ is continuous.}
  \end{equation}
\end{lma}

\section{Finite size criterion for percolation}
\label{s:finite-size}

The first tool in our study of the phase transition for Poisson Boolean percolation in $\R^2$ will be the following result that bootstraps the probability of crossing a rectangle: from close to one to converging fast to one.
This will be fundamental in this paper; as a first indication of its importance we show how it relates to the existence of an unbounded component, and the non-triviality of the threshold parameters $\lambda_0$ and $\lambda_1$.

\begin{prop}\label{prop:finite-size}
Assume the second moment condition~\eqref{eq:7}. Then, there exists a constant $\theta = \theta(\mu) > 0$ and an increasing function $r_0:[0,\infty)\to[0,\infty)$ such that the following are equivalent:
  \begin{enumerate}[\quad (i)]
  \item\label{item:8} There exists $\x\ge r_0(\lambda) $ such that
    $\Pr_\lambda(\cross(3\x,\x))>1-\theta$.
  \item\label{item:7} $\displaystyle\lim_{\x\to\infty}\Pr_\lambda(\cross(3\x,\x))=1$.
\end{enumerate}
Under the stronger condition~\eqref{e:log_condition}, if either of~(i) or~(ii) holds, then
\begin{enumerate}[\quad (i)]
\setcounter{enumi}{2}
  \item $\sum_{k\ge0}\big[1-\Pr_\lambda(\cross(3^{k+1},3^k))\big]<\infty$.
  \end{enumerate}
These statements remains valid for $\cross$ replaced by $\cross^\star$.
\end{prop}

\begin{proof}
Fix $\theta=1/100$ and let $C$ be a sufficiently large constant so that, for all $r\ge1$,
$$
\rho_\lambda(5r,r)\le g(r) := C \lambda \int_{r/2}^\infty x^2\mu(\d x).
$$
The existence of the constant follows from Proposition~\ref{c:decouplef1f2}.
Set $r_0:=\min\{r\ge1:g(r)\le\theta/2\}$.
Since $g(r)$ increases with $\lambda$ it is clear that $r_0$ does too.
Since~\emph{(ii)} trivially implies~\emph{(i)}, it will suffice to prove that~\emph{(i)} implies~\emph{(ii)} and, under the stronger assumption~\eqref{e:log_condition}, that~\emph{(ii)} implies~\emph{(iii)}.

Let $p(\x):=1-\Pr_\lambda(\cross(3\x,\x))$ and let $B_k:=\cross(9\x,\x)$.
Note that a $9\x\times\x$-rectangle may be tiled with seven $3\x\times\x$-rectangles, four positioned horizontally and three vertically, in such a way that if each is crossed in the `hard' direction, then the $9\x\times\x$-rectangle is crossed horizontally.
Consequently, using the union bound, we obtain that
$$
\Pr_\lambda(B_\x)\,\ge\,1-7\Pr_\lambda(\neg\cross(3\x,\x))\,=\,1-7p(\x).
$$
Let $B_\x'$ denote the translate of $B_\x$ along the vector $(0,2\x)$. Since the occurrence of either of $B_\x$ or $B_\x'$ implies the occurrence of $\cross(9\x,3\x)$, he have
$$
p(3\x)\,\le\,\Pr_\lambda(\neg B_\x \cap \neg B_\x')\,\le\,\Pr_\lambda(\neg B_\x)^2+\rho_\lambda(5r,r).
$$
Hence, by definition of $g$, we obtain for every $r\ge1$ the bound
\begin{equation}\label{eq:contraction}
p(3r)\le49p(r)^2+g(r).
\end{equation}

Now, assume there exists $r\ge r_0$ so that $p(r)<\theta$. Then, via iterated use of~\eqref{eq:contraction}, we find that $p(3^kr)<\theta$ for all $k\ge0$.
Consequently, further use of~\eqref{eq:contraction} gives, for $\ell=1,2,\ldots,k$
\begin{equation}\label{eq:it_contraction}
p(3^kr)\,\le\,\frac{1}{2^\ell}p(3^{k-\ell}r)+\sum_{j=1}^\ell\frac{1}{2^{j-1}}g(3^{k-j}r)\,\le\,\frac{1}{2^\ell}+2g(3^{k-\ell}r).
\end{equation}
Hence, sending $\ell$ to infinity with $k$ shows that $\lim_{k\to\infty}p(3^kr)=0$.
Since for every $\x'\in[\x,3\x]$ we have $p(\x')\le\Pr_\lambda(\neg B_\x)\le7p(\x)$, it follows that $\lim_{r\to\infty}p(r)=0$, so~\emph{(ii)} holds.

Finally, we assume that~\emph{(ii)} holds and pick $k_0$ so that $3^{k_0}\ge r_0$ and $p(3^{k_0})<\theta$. Summing over $k$ in~\eqref{eq:it_contraction}, with $r=3^{k_0}$ and $\ell=\lfloor k/2\rfloor$, leads to
$$
\sum_{k\ge1}p(3^{k+k_0})\,\le\,\sum_{k\ge1}\frac{1}{2^{(k-1)/2}}+2\sum_{k\ge1}g(3^{k/2+k_0})\,\le\,4+2\int_0^\infty g(3^{y/2})\,\d y.
$$
Since $g$ is itself an integral we may use Fubini's theorem to reverse the order of integration, and obtain
$$
\sum_{k\ge1}p(3^{k+k_0})\,\le\, 4+C\lambda\int_0^\infty x^2\log x\,\mu(\d x),
$$
for a possibly larger constant $C$, which is finite under the assumption~\eqref{e:log_condition}, and~\emph{(iii)} follows.
\end{proof}

\begin{remark}
Note that the proof of Proposition~\ref{prop:finite-size} assumes very little about the underlying percolation model. Indeed, the argument applies to any percolation model on $\R^2$ satisfying the invariance assumption~\eqref{e:invariance} and decay of correlations corresponding to conditions~\eqref{eq:7} or~\eqref{e:log_condition}.
\end{remark}

The following more quantitative statement sheds some light on the role of spatial correlations on the rate of convergence in Proposition~\ref{prop:finite-size}.
Given a function $f:[1, \infty) \to (0, \infty)$ we say that it is \emph{regularly varying} if for any $a \in (0, \infty)$
\begin{equation}\label{eq:reg_var}
\frac{f(a\x)}{f(\x)} \text{ converges to a non-zero limit.}
\end{equation}
The following proposition says, in particular, that if the spatial dependence decays polynomially fast, then the convergence in Proposition~\ref{prop:finite-size} occurs (at least) at the same polynomial rate.

\begin{prop}\label{prop:finite_quant}
Assume the finite moment condition~\eqref{eq:7} and let $f:[1,\infty)\to(0,\infty)$ be a regularly varying function satisfying $\lim_{\x\to\infty}f(\x)=0$ and $\rho_\lambda(5r,r)\le f(r)$ for all $r \ge 1$.
Then, if $\lim_{\x\to\infty}\Pr_\lambda(\cross(3\x,\x))=1$, then there exists $c<\infty$ such that
$$
\Pr_\lambda(\cross(3\x,\x))\ge1-cf(\x)\quad\text{for all }r\ge1.
$$
The statement remains true if $\cross$ is replaced by $\cross^\star$.
\end{prop}

\begin{proof}
With $p(\x):=1-\Pr_\lambda(\cross^\star(3\x,\x))$, we may repeat the proof of Proposition~\ref{prop:finite-size} to obtain
\begin{equation}\label{eq:f-contraction}
p(3\x)\le49p(\x)^2+f(\x).
\end{equation}
Since $f$ is regularly varying there exists a constant $c' \in (0, 50)$ such that $f(3\x)/f(\x) \ge c'$ for all $r\ge1$. Let $\theta'=c'/100$ and pick $r'$ such that $p(r)\le\theta'$ and $f(r)\le\theta'/2$ for all $r\ge r'$. Let $c=\max\{2/c',1/f(r')\}$. We claim that
$$
p(3^kr')\le cf(3^kr')\quad\text{for all }k\ge0.
$$
The case $k = 0$ is immediate, and the remaining cases follows from~\eqref{eq:f-contraction} via an induction step. This proves the statement of the proposition along exponentially growing sub-sequences from which one can easily derive the full statement.
\end{proof}

Propositions~\ref{prop:finite-size} and~\ref{prop:finite_quant} show how an observation in a finite region can be sufficient to draw conclusions regarding crossing probabilities over arbitrarily large regions.
The next well-known result connects the finite observation to the existence of unbounded components.

\begin{prop}\label{prop:finite-infinite}
The condition $\sum_{k\ge0}\big[1-\Pr_\lambda(\cross(3^{k+1},3^k))\big]<\infty$ implies that the probability, at $\lambda$, that the origin is contained in an unbounded occupied component is positive.
The same holds when $\cross$ and `occupied' are replaced by $\cross^\star$ and `vacant', respectively.
\end{prop}

\begin{proof}
Denote by $R_k$ the rectangle $[0,3^{k+1}]\times[0,3^k]$, and by $R_k'$ the rectangle $[0,3^k]\times[0,3^{k+1}]$, resulting from a right-angle rotation of $R_k$.
Consider the sequence, indexed by $k\ge0$, alternating between $R_k$ or $R_k'$, depending on whether $k$ is odd or even.
The above condition implies that the expected number of these rectangles that fail to be crossed in the `hard' direction is finite.
Hence, by Borel-Cantelli, it follows that all but finitely many contain a crossing in the `hard' direction almost surely.
However, all but finitely many of these crossings must intersect as a consequence of how the rectangles are placed.
Thus, with probability one, there is an unbounded occupied component.
Due to invariance with respect to translations it follows that the origin has positive probability to be contained in this component.
\end{proof}

We close this section by showing that the (finite-size) phase transition is non-trivial.

\begin{cor}\label{cor:nontriviality}
  Assume the second moment condition~\eqref{eq:7}.
  Then,
  $$
  0<\lambda_0\le\lambda_c\le\lambda_1<\infty.
  $$
  Under the stronger assumption~\eqref{e:log_condition} we also have $\lambda_0\le\lambda_c^\star\le\lambda_1$.
\end{cor}

\begin{proof}
  Assume that~\eqref{eq:7} holds.
  We first show that $\lambda_0 > 0$.
  Let $\theta>0$ and $r_0:[0,\infty)\to[0,\infty)$ be as in Proposition~\ref{prop:finite-size}.
  Clearly $\Pr_0(\cross(\x, 3\x))=0$ for all $\x\ge1$, in particular for $\x=r_0(1)$.
  That the crossing probability varies continuously with respect to the intensity parameter is an easy consequence of~\eqref{e:decouple1} of Lemma~\ref{l:decouple} and the fact that it is unlikely to add any point at all to a bounded region for small enough changes in the parameter.
  Hence, we may pick $\eps \in (0,1)$ such that $\Pr_\eps(\cross(\x,3\x)) < \theta$, implying that
  \begin{equation*}
    \lim_{\x\to\infty} \Pr_\eps(\cross(\x,3\x)) = 0,
  \end{equation*}
  and thus that $\lambda_0 \ge \eps > 0$.

  The fact that $\lambda_1 < \infty$ follows from the fact that this holds for the model with constant radius, together with a stochastic domination argument.

  We next show that $\lambda_c\ge\lambda_0$ and $\lambda_c^\star\le\lambda_1$.
  First, consider $\lambda<\lambda_0$. Then, since a path connecting the origin to $\partial B(r)$ has to cross one of four rectangles of dimension $r/3\times r$ in the easy direction, the union bound and the definition of $\lambda_0$ gives that
  $$
  \lim_{\x\to\infty}\Pr_\lambda(0\lro\partial B(r))\,\le\,\lim_{\x\to\infty}4\Pr_\lambda(\cross(r/3,r))\,=\,0.
  $$
  Thus, $\lambda\le\lambda_c$, which shows that $\lambda_0\le\lambda_c$. An analogous argument shows that $\lambda_c^\star\le\lambda_1$,

  To prove that $\lambda_c^\star\ge\lambda_0$, we make the stronger assumption that~\eqref{e:log_condition} holds and fix $\lambda<\lambda_0$, in which case we have $\lim_{\x\to\infty}\Pr_\lambda(\cross^\star(3\x,\x))=1$.
  Then, Proposition~\ref{prop:finite-size} shows that
\begin{equation}\label{eq:finite_sum}
\sum_{k\ge0}\big[1-\Pr_\lambda(\cross^\star(3^{k+1},3^k))\big]<\infty,
\end{equation}
which by Proposition~\ref{prop:finite-infinite} implies the almost sure existence of an unbounded vacant component, and thus that $\lambda\le\lambda_c^\star$. Hence $\lambda_c^\star\ge\lambda_0$.

An analogous argument shows, under the additional condition~\eqref{e:log_condition}, that $\lambda_c\le\lambda_1$. However, we claim that this additional assumption is unnecessary. Our goal will therefore be to show that~\eqref{eq:finite_sum} holds, with $\cross^\star$ replaced by $\cross$, also under the weaker second moment condition~\eqref{eq:7}.
As before, let $p(r)=1-\Pr_\lambda(\cross(3\x,\x))$. Repeating the first steps of the proof of Proposition~\ref{prop:finite-size}, observing that the complements of occupied crossing events are decreasing, we obtain the following analogue to~\eqref{eq:contraction}
$$
p(3r)\,\le\,7p(r)\big(7p(r)+\bar\rho_\lambda(5r,r)\big).
$$
Since both $p(r)$ and $\bar\rho_\lambda(5r,r)$ tend to zero as $r\to\infty$, we may find $r_0$ such that $7p(r)+\bar\rho_\lambda(5r,r)$ is bounded by $\frac{1}{14}$ for all $r\ge r_0$. So, for some $k_0\ge1$, we obtain that
$$
p(3^{k+k_0})\,\le\,\frac{1}{2}p(3^{k-1+k_0}r_0)\,\le\,\frac{1}{2^k}p(3^{k_0}),
$$
which is summable. Hence~\eqref{eq:finite_sum} holds for $\cross$ also under the weaker condition~\eqref{eq:7}, and $\lambda_c\le\lambda_1$ follows.
\end{proof}

\section{Russo-Seymour-Welsh theory}
\label{s:critical}
In the study of planar percolation, Russo-Seymour-Welsh (RSW) techniques play a
central role and have numerous consequences. The original proof for Bernoulli
percolation (\cite{russo1978note,seymour1978percolation}) is strongly based on planarity and the independence structure of the Bernoulli percolation measure, and does not extend easily to other contexts.
Considerable technicalities had to be overcome even in the extension to Poisson Boolean percolation in $\R^2$ with fixed radii, see \cite{roy90,alexander96}, and in the case of (unbounded) random radii such a result has until this point not been obtained.
In the last few years some new arguments have been developed to prove RSW-results for
dependent percolation models, e.g.\@ Voronoi percolation
\cite{bolrio06-2,tassion14+,ahlberg2015quenched} and the random-cluster model
\cite{beffara2012self,duminil2015continuity}.

In this section we develop an RSW-theory applicable for Poisson Boolean percolation.
Our method of proof will be greatly inspired by that of~\cite{tassion14+}. However, due to the
asymmetry between the vacant and occupied regions, we need a stronger version of the result proved in \cite{tassion14+}. In that paper, the RSW statement bounds the
crossing probabilities for rectangles in the long direction, assuming a bound on
the crossing probabilities for squares, see \eqref{eq:4}. Here we assume only
that rectangles are crossed in the easier direction.
Our proof will be rather general and applies in settings far beyond Poisson Boolean percolation; see Remark~\ref{rem:scope} below.

\begin{thm}[RSW-Theorem]
  \label{thm:rsw}
  Assume the finite second moment condition stated in \eqref{eq:7}.
  Then, for any $\lambda > 0$, if for some $\kappa>0$ we have
  \begin{equation}
    \label{e:easy_crossing}
    \inf_{\x \geq 1} \mathbb{P}_\lambda\big(\cross(\kappa \x, \x) \big) > 0,
  \end{equation}
  then the same is true for all $\kappa>0$.
  The same holds for $\cross$ replaced by $\cross^\star$.
\end{thm}

\begin{remark}
  \label{rem:scope}
  Throughout this section, we are going to establish the above result and some of its consequences for Poisson Boolean percolation with finite second moments.
  However, let us emphasize that the same proof works for several different types of percolation measures on the plane. More precisely, the only properties of the random set $\mathcal{O}$ that we use in this section are:
  \begin{equation}
    \label{e:natural}
    \begin{array}{c}
      \text{the translation, reflection and rotational symmetries in \eqref{e:invariance},}\\
      \text{the FKG inequality \eqref{eq:3},}\\
      \text{the decay of correlations stated in \eqref{e:rho_to_zero} and}\\
      \text{the continuity of the crossing probabilities stated in Lemma~\ref{l:cross_continuous}}.
    \end{array}
  \end{equation}
\end{remark}

\subsection{Consequences of Theorem~\ref{thm:rsw}}
\label{sec:stat-theor-cons}

The Russo-Seymour-Welsh Theorem stated above, which will be proved in Section~\ref{sec:proof-theorem}, has several important consequences that we develop next. These consequences concern `box-crossing' and `one-arm' probabilities in the critical regime $[\lambda_0,\lambda_1]$. However, we first state a `high-probability' version of Theorem~\ref{thm:rsw}. All proofs of the consequences listed here will be given in Section~\ref{sec:proof-cons} below.

\begin{cor}\label{cor:high-probab}
  Assume the finite second moment condition~\eqref{eq:7}.
  Given $\lambda > 0$, if for some $\kappa>0$ we have
  \begin{equation*}
    \lim_{\x\to\infty}\mathbb{P}_\lambda \big(\cross(\kappa \x, \x) \big) = 1,
  \end{equation*}
  then it is true for all $\kappa>0$.
  The same holds for $\cross$ replaced by $\cross^\star$.
\end{cor}

Recall the definition of $\lambda_0$ and $\lambda_1$ in~\eqref{eq:11} and~\eqref{eq:9}.
Outside of the interval $[\lambda_0,\lambda_1]$ the probability of crossing a fixed-ratio rectangle will converge rather rapidly to either zero or one as the side length of the rectangle increases.
This is merely a matter of definition and the FKG-inequality.
One of the main consequences of the RSW Theorem is that within the interval $[\lambda_0,\lambda_1]$ the probability of crossing a fixed-ratio rectangle remains rather balanced.
This is sometimes referred to as the `box-crossing' property, which will also be used to provide bounds on the one-arm probabilities.

For $0 < r < r'$, we define
\begin{equation}
  \label{e:def_arm}
  \arm(r, r') := \big[ B^\infty(r) \overset{\mathcal{O}}{\longleftrightarrow} \partial B^\infty(r') \big].
\end{equation}
As before we write $\arm^\star$ for the above event, where $\mathcal{O}$ is replaced by $\mathcal{V}$.
Given $\x' > \x \geq 1$ we also introduce
\begin{equation}
  \label{e:def_annulus}
  \cir(\x, \xprime) := \boxAnnulus{2\x}{2\x'} := \arm^\star(r, r')^c
\end{equation}
Finally we let $\cir^\star(\x,\xprime) := \arm(r, r')^c$.

\nc{c:q-bound}
\begin{cor}\label{cor:box-crossing}
  Assume the second moment condition in~\eqref{eq:7}.
  Then, for every $\kappa > 0$, there exists $\uc{c:q-bound} = \uc{c:q-bound}(\kappa) > 0$ such that for every $\lambda \in [\lambda_0, \lambda_1]$ we have
  \begin{equation}
    \label{item:9}
    \Pr_\lambda(\cross(\kappa\x, \x)) \in (\uc{c:q-bound}, 1 - \uc{c:q-bound})\quad \text{for every $\x \ge 1$,}
  \end{equation}
  and also
  \begin{equation}
    \label{e:cir_non_trivial}
    \mathbb{P}_\lambda (\arm(\x, 2\x)) \in (c', 1-c')\quad \text{for every $\x \geq 1$,}
  \end{equation}
  where $c'=\uc{c:q-bound}(4)^4$.
  As before, the above result also holds for the vacant set, i.e.\ with $\cross$ and $\arm$ replaced by $\cross^\star$ and $\arm^\star$ respectively.
\end{cor}

Another consequence of Theorem~\ref{thm:rsw} concerns the decay of arm probabilities
at criticality. Unlike previous results of this section, the bounds on the
arm events distinguish between vacant and occupied sets. The one-arm probability
always decays polynomially fast for the vacant set. For the occupied set, we can
prove polynomial decay of the one-arm probability only under the stronger
assumption \eqref{e:alpha_condition} and we explain why this restriction is
necessary in Remark~\ref{r:D_is_bad}.

\begin{cor}[Bounds on the arm-events]
  \label{cor:arm_events}
  Assume the second moment conditions in~\eqref{eq:7}.
  \begin{enumerate}[(i)]
  \item\label{item:10} There exists a function $f: (0,1) \to (0, 1)$ such that $\lim_{x \to 0}
    f(x) = 0$ and for every $\lambda \in [\lambda_0, \lambda_1]$ we have
    \begin{equation}
      \label{eq:arm_decay}
      \Pr_\lambda(\arm(r,r')) \le f \Big( \frac{r}{r'} \Big).
    \end{equation}

  \item\label{item:17} If $\arm$ is replaced by $\arm^\star$, then a stronger conclusion holds:
    there exists $c > 0$ such that
    \begin{equation}
      \label{eq:13'}
      \Pr_\lambda ( \arm^\star(r,r') ) \le \frac1c \left(\frac r{r'}\right)^c.
    \end{equation}
  \item Under the stronger assumption~\eqref{e:alpha_condition} the conclusion
    of~\eqref{eq:13'} holds also for $\arm$.
  \end{enumerate}

\end{cor}

\begin{remark}
  \label{r:D_is_bad}
  In Corollary~\ref{cor:arm_events}, under the $2 + \eps$-moment condition \eqref{e:alpha_condition}, we can choose $f(x) = \tfrac1c x^c$ for some constant $c > 0$.
  However, if $\mu([x, \infty))$ is $1/(x^2 \log^2(x))$, then the second moment condition holds but the arm event does not decay polynomially.
  To see this, consider the event that the origin is covered by a occupied disk of radius at least $r$, and observe that the probability of this event does not decay polynomially fast in $r$.
\end{remark}

An immediate consequence of Corollary~\ref{cor:arm_events} is the following observation.

\begin{cor}\label{cor:no_cluster}
  Assume the second moment conditions in~\eqref{eq:7}. For every $\lambda\in[\lambda_0,\lambda_1]$ there is almost surely no unbounded occupied nor vacant cluster. Consequently,
  $$
  \lambda_c^\star\le\lambda_0\le\lambda_1\le\lambda_c.
  $$
\end{cor}

\subsection{Standard inequalities}
\label{sec:stand-ineq}

Before starting the proof of the RSW Theorem and its consequences, let us recall some standard
inequalities on crossing probabilities.

\begin{lma}\label{lem:standardInequalities}
  For every $\lambda>0$, $r>0$, $\kappa>0$ and integer $j\ge1$ we have:
  \begin{enumerate}[(i)]
  \item\label{item:13} $\Pr_\lambda[\cross((1+j\kappa)r,r)]\ge
    \Pr_\lambda[\cross((1+\kappa)r,r)]^{2j-1}$,
  \item\label{item:14} $\Pr_\lambda[\cross(r,(1+\kappa)r)]\ge 1-\left(1-
   \Pr_\lambda[\cross(r,(1+j\kappa r)]\right)^{1/(2j-1)}$,
  \item\label{item:16} $\Pr_\lambda[\cir(r,2r)]\ge
  \Pr_\lambda[\cross(4r,r)]^4$,
  \item\label{item:15} $\Pr_\lambda[\cross(2r,r)]\ge \Pr[\cir(r,2r)]$.
  \end{enumerate}
  The same holds if we replace $\cross$ and $\cir$ by $\cross^\star$ and
  $\cir^\star$, respectively.
\end{lma}

\begin{proof}
  Part~\eqref{item:13} is a straightforward consequence of the FKG-inequality~\eqref{eq:3}: A horizontal crossing of an $(1+j\kappa)\x\times\x$-rectangle can be enforced by $j$ horizontal crossings of $(1+\kappa)\x\times\x$-rectangles and $j-1$ vertical crossings of $\x\times(1+\kappa)\x$-rectangles overlapping one another. Similarly, part~\eqref{item:16} is the consequence of a circuit being possible to construct out of four overlapping rectangle crossings.

  Part~\eqref{item:14} is a direct consequence of~\eqref{item:13} since it may be rewritten as
  $$
  \Pr_\lambda[\cross^\star((1+\kappa)r,r)]\le\Pr_\lambda[\cross^\star((1+j\kappa)r,r)]^{1/(2j-1)}.
  $$

  Inequality~\eqref{item:15} follows
  from the observation that any circuit in the annulus $B^\infty(2r)\setminus
  B^\infty(r)$ must cross the rectangle $[-r,r]\times[r,2r]$ horizontally.
\end{proof}

\subsection{A useful circuit Lemma}
\label{sec:usef-circ-lemm}

For Bernoulli percolation, the bounds on crossing probabilities provide some
bounds on the arm events. This is based on a circuit argument, that uses
independence. In our case, the argument need to be adapted because of spatial
dependencies.

\begin{lma}[Circuit Lemma]
  \label{lma:circuit}
  Assume the second moment conditions in~\eqref{eq:7}. Let $\lambda\in
  [0,\infty)$.
  \begin{enumerate}[(i)]
  \item\label{item:11} For every $c>0$, there exists a function $f=f_{c,\lambda}: (0,1) \to (0,
    1)$ such that $\lim_{x \to 0} f(x) = 0$ and for all $r' \ge 2r \ge 2$,
    \begin{equation}
      \label{eq:13}
      \text{if }\displaystyle\inf_{r\le s\le r'/2} \Pr_\lambda[\cir^\star(s,2s)]\ge c,
      \text{ then }\Pr_\lambda [ \cir^\star(r,r') ] \ge 1- f \left(\tfrac r{r'}\right).
    \end{equation}

  \item\label{item:6} For every $c>0$, there exists a constant $c'=c'(c,\lambda)>0$ such that for all $r' \ge 2r \ge 2$,
    \begin{equation}
      \label{eq:13''}
      \text{if }\displaystyle\inf_{r\le s\le r'/2} \Pr_\lambda[\cir(s,2s)]\ge
      c,
     \text{ then }\Pr_\lambda [ \cir(r,r') ] \ge 1-\tfrac1{c'}\left(\tfrac
        r{r'}\right)^{c'}.
    \end{equation}

  \item\label{item:12} Under the $2+\eps$ moment
    condition~\eqref{e:alpha_condition}, Item~$(\ref{item:6})$ holds also for
    $\cir$ replaced by $\cir^\star$.
  \end{enumerate}
\end{lma}

\begin{proof}
  We begin with the proof of Item~(\ref{item:11}). It suffices to prove the
  statement for $r'\ge 16r$. Let $r \ge 1$ and $r'\ge 16r$ be given, and assume
  that
  \begin{equation}
    \inf_{r\le s\le r'/2} \Pr_\lambda[\cir^\star(s,2s)]\ge c>0.\label{eq:6}
  \end{equation}
  Let $g(s) := \sup_{s'\ge s} \rho_\lambda(s',2s')$.
  By Proposition~\ref{c:decouplef1f2}, we have
  $\lim_{s\to\infty}g(s) = 0$.

  Set $t = \sqrt{rr'}$ and $\ell_i=4^it$. Note that if $\cir^\star(r,r')$ fails
  to occur, then $A_i=\cir^\star(\ell_i,2\ell_i)$ cannot occur for no
  $i=0,1,\ldots,k-1$ where $k=\lfloor\frac12\log_4(r'/r)\rfloor$ ($k\ge1$ since
  we assumed $r'\ge 16r$). The hypothesis~\eqref{eq:6} implies that
  $\Pr_\lambda(A_i^c)\le1-c<1$. Therefore,
    \begin{align}
    \label{eq:19}
    \Pr_\lambda \Big[ \mcap_{i = 0}^{k-1} A_i^c \Big]
    \,&\le\, \Pr_\lambda[A_{k-1}] \Pr_\lambda\Big[ \mcap_{i = 0}^{k-2} A_i^c \Big] + \rho_\lambda (2\ell_{k-2},\ell_{k-2})\\
    &\le\, (1-c) \Pr_\lambda \Big[ \mcap_{i = 0}^{k-2} A_i^c \Big] + g(t),
  \end{align}
  which by induction gives that
  \begin{equation}
    \label{eq:20}
    \Pr_\lambda \Big[ \mcap_{i = 0}^{k-1} A_i^c \Big] \,\le\, (1-c)^k + \frac{1}{c} g(t).
  \end{equation}
  By the choice of $k$, and the fact that $g(t)\le g(t/r) = g(\sqrt{r'/r})$, we
  finally obtain
  \begin{equation}
    \label{eq:21}
    \Pr_\lambda \big[ \cir^\star(r,r')^c \big] \,\le\,\Pr_\lambda \Big[ \mcap_{i = 0}^{k-1} A_i^c \Big] \,\leq\, \frac{1}{1-c}\Big( \frac{r'}r \Big)^{\frac12\log_4(1-c)}+ \frac{1}{c}g\left(\sqrt{\tfrac {r'}r}\right).
  \end{equation}
  Hence,~\eqref{eq:arm_decay} holds with $f(x)=\frac{1}{1-c}x^\alpha+\frac{1}{c}g(x^{-1/2})$ for any sufficiently small constant $\alpha>0$.
  The stronger moment condition~\eqref{e:alpha_condition} implies~\eqref{e:ralpha_decay}, which readily gives that $f$ may be chosen of the form $\frac{1}{c'}x^{c'}$ for some constant $c'=c'(c,\lambda)>0$, which also proves Item~(\ref{item:12}).

  We now turn to the proof of Item~(\ref{item:6}) and reuse the above notation.
  Set $B_i=\cir(\ell_i,2\ell_i)$.
  As above, by~\eqref{e:cir_non_trivial} it follows that $\Pr_\lambda\big(B_i^c\big)\le 1-c$ for every $i=0,1,\ldots,k-1$.
  By Proposition~\ref{p:bar_rho} we may fix $r_0$ large enough so that
  \begin{equation}
    \label{e:rho_bar_c2}
    \bar{\rho}_\lambda(2r,r) \leq c/2\quad\text{for every }r\ge r_0.
  \end{equation}
  Observe that the events $B_i$ are decreasing. Hence, for $r\ge r_0$, $r'\ge16r$ and $t=\sqrt{rr'}$, we obtain
  \begin{equation}
    \label{eq:16}
    \begin{aligned}
      \Pr_\lambda \big( \cir(r, r')^c \big) \,& \leq\, \Pr_\lambda \Big( \mcap_{i = 0}^{k-1}B_i^c \Big)\\
      & \leq \,\Pr_\lambda \Big( \mcap_{i = 1}^{k-1} B_i^c \Big) \Big( \mathbb{P}_\lambda (B_0^c) + \bar{\rho}_\lambda(2\ell_0, \ell_0) \Big)\\
      & \leq \,\prod_{i = 0}^{k-1} \Big( \mathbb{P}_\lambda \big( B_i^c \big) + \bar{\rho}_\lambda(2\ell_i, \ell_i) \Big)\\
      & \leq \,(1-c/2)^{\frac12\log_4(r'/r)-1}.
    \end{aligned}
  \end{equation}
  By allowing for a larger constant on the right-hand side we may replace the assumption $r\ge r_0$ with $r\ge1$. This concludes the proof of \eqref{eq:13}.
\end{proof}

\subsection{Proof of the consequences of Theorem~\ref{thm:rsw}}
\label{sec:proof-cons}

We now assume the validity of Theorem~\ref{thm:rsw} and prove Corollaries~\ref{cor:high-probab}-\ref{cor:no_cluster}.

\begin{proof}[Proof of Corollary~\ref{cor:high-probab}]
We will show that if $\lim_{\x\to\infty}\Pr(\cross(\kappa\x,\x))=1$ for some $\kappa>0$, then $\lim_{\x\to\infty}\Pr(\cross(2\kappa\x,\x))=1$, from which the statement follows via iteration.

Partition the right-hand side of the rectangle $R=[0,\kappa\x]\times[0,\x]$ into
$n$ equal intervals (of length $\x/n$), and denote the event that there is a horizontal crossing of $R$ into the $k$-th of these intervals by $A_k$. Using the square-root trick~\eqref{eq:32}, we find that
$$
\max_{1\le k\le n}\Pr_\lambda(A_k)\ge1-[1-\Pr_\lambda(\cross(\kappa\x,\x))]^{1/n},
$$
which by assumption tends to $1$ as $\x\to\infty$.

By assumption, Theorem~\ref{thm:rsw} shows that $\inf_{\x\ge1}\Pr_\lambda(\cross(4\x,\x))>0$. Hence, the standard inequalities of Lemma~\ref{lem:standardInequalities} show that $\inf_{\x\ge1}\Pr_\lambda(\cir(\x,2\x))>0$. Hence, by Lemma~\ref{lma:circuit} we conclude that for large enough $\x$,
$$
\Pr_\lambda(\cir( \x/n,\kappa \x/2))\ge1-f(\kappa n/2),
$$
for some function $f$ such that $f(x)\to0$ as $x\to\infty$.

Finally, we combine the above estimates to obtain that for every $\eps>0$
$$
\Pr_\lambda\big(\cross(2\kappa\x,\x)\big)\,\ge\,\Pr_\lambda\Bigg[
\boxRings{2\kappa\x}{\x}{\tfrac{\x}{n}}{\kappa\x}
\Bigg]\,\ge\,\Pr_\lambda(A_k)^2\,\Pr_\lambda\big(\cir(\x/n,\kappa\x/2)\big)\,>\,1-\eps,
$$
for some choice of $k$ and $n$, and all large $\x$, as required.
\end{proof}

\begin{proof}[Proof of Corollary~\ref{cor:box-crossing}]
  Take $\theta > 0$ and $r_0 = r_0(\lambda_1)$ so that Proposition~\ref{prop:finite-size} is in force.
  Then,
  \begin{equation}
    \Pr_{\lambda_1}(\cross(3\x, \x))\le1-\theta \text{ for all $\x \geq \x_0$.}
  \end{equation}
  Otherwise, by continuity, we could find $\eps > 0$ and $\x \ge r_0$ such that $\Pr_{\lambda_1-\eps}(\cross(3\x,\x))>1-\theta$, which by Proposition~\ref{prop:finite-size} would imply that the probability indeed converges to 1; a contradiction to the definition of $\lambda_1$.
  Analogously, we have $\Pr_{\lambda_0}(\cross(\x, 3\x)) \ge \theta$ for all large $\x$.
  Using Theorem~\ref{thm:rsw} and its dual version this establishes \eqref{item:9}.

  Equation \eqref{e:cir_non_trivial} follows from the fact that the complement
  of $\arm(\x,3\x)$ is $\cir^\star(\x,3\x)$ and the standard inequalities of
  Lemma~\ref{lem:standardInequalities}.
\end{proof}

\begin{proof}[Proof of Corollary~\ref{cor:arm_events}]
This is an immediate consequence of Corollary~\ref{cor:box-crossing} and Lemma~\ref{lma:circuit}, as $\Pr_\lambda(\arm(r,r'))=1-\Pr_\lambda(\cir^\star(r,r'))$.
\end{proof}

\begin{proof}[Proof of Corollary~\ref{cor:no_cluster}]
  The bounds $\lambda_c\ge \lambda_1$ and $\lambda_c^\star\le \lambda_0$ follow respectively from Item~(\ref{item:10}) and Item~(\ref{item:17}) of Corollary~\ref{cor:arm_events}.
\end{proof}

\subsection{Proof of Theorem~\ref{thm:rsw}}
\label{sec:proof-theorem}

\nc{c:easy_cross}
In this section we prove Theorem~\ref{thm:rsw}.
We assume throughout the proof the second moment condition~\eqref{eq:7} and that $\lambda>0$ is fixed. By assumption there exists $\kappa\in(0,\tfrac13]$ such that $\inf_{\x\ge1}\Pr_\lambda(\cross(\kappa\x,\x))>0$. Clearly, there is no restriction assuming the upper bound on $\kappa$. Hence, we note that either also $\inf_{\x\ge1}\Pr_\lambda(\cross^\star(\kappa\x,\x))>0$, or we have
$$
\sup_{\x\ge1}\Pr_\lambda(\cross(3\x,\x))=1.
$$
The latter implies, via Proposition~\ref{prop:finite-size}, that $\lim_{\x\to\infty}\Pr_\lambda(\cross(3\x,\x))=1$, in which case there is nothing left to prove. So, we may without loss of generality assume that for some $\kappa>0$
\begin{equation*}
  \inf_{\x \geq 1} \mathbb{P}_\lambda
    \big(\cross(\kappa \x, \x) \big) > 0\quad\text{and}\quad\inf_{\x \geq 1} \mathbb{P}_\lambda
    \big( \cross^\star(\kappa \x, \x) \big) > 0.
\end{equation*}
Moreover, by the standard inequalities of Lemma~\ref{lem:standardInequalities} (part~\eqref{item:14} and its dual version) this is equivalent to assuming the existence of a constant $\uc{c:easy_cross}>0$ such that
\begin{equation}
  \label{e:easy_cross}
  \inf_{\x \geq 1} \mathbb{P}_\lambda \big( \cross \big(\x, \b \x \big) \big) \geq \uc{c:easy_cross} \quad \text{and} \quad \inf_{\x \geq 1} \mathbb{P}_\lambda \big( \cross^\star \big(\x, \b \x \big) \big) \geq \uc{c:easy_cross}.
\end{equation}
In what follows we will explore several consequences of the above assumption.
This will be done in a series of lemmas that will culminate in the proof of Theorem~\ref{thm:rsw}.

Throughout this section we will need to introduce a range of events that are closely related to the crossing event defined in \eqref{e:def_cross}. We will typically depict these events graphically, since we believe the proof is better explained this way. However, we emphasize that it is always possible to give a formal definition.
For instance, recall the definition of $[A \overset{\mathcal{C}}{\leftrightarrow} B]$ in~\eqref{e:path_in_C}.
We start with the definition of the event of having a vacant crossing above an occupied crossing, both reaching the right-hand side on the upper half:
\begin{equation}
  \boxHighTwoCross{\x}{h} := \bigcup_{q \in [0, h] \cap \mathbb{Q}} \Big[ \big( \{0\} \times [0,q] {\overset{\mathcal{O} \cap K}\longleftrightarrow} A \big) \cap \big( \{0\} \times [q, h] {\overset{\mathcal{V} \cap K}\longleftrightarrow} A \big) \Big],
\end{equation}
where $K = [0, \x] \times [0,h]$ and $A = \{\x\} \times [h/2, h]$.

\nc{c:prob_dual}
\nc{c:box_height}
\begin{lma}
  \label{l:long_rectangle}
  There exist $r_0\ge 1$, $0 < \uc{c:prob_dual} < 1/8$ and $\uc{c:box_height} \ge
  23$ such that
  \begin{equation}
    \label{e:high_two_cross}
    \inf_{\x \geq r_0} \mathbb{P}_\lambda \Bigg[ \boxHighTwoCross{\x/\uc{c:box_height}}{\tfrac{\x}{3}} \Bigg] \geq \uc{c:prob_dual}
  \end{equation}
  and
  \begin{equation}
    \label{e:cross_long}
    \sup_{\x \geq r_0} \mathbb{P}_\lambda \big( \cross(\tfrac{\x}{3}, \tfrac{\x}{\uc{c:box_height}}) \big) \leq \frac{1}{16} \quad \text{and} \quad \sup_{\x \geq r_0} \mathbb{P}_\lambda \big( \cross^\star(\tfrac{\x}{3},\tfrac{\x}{\uc{c:box_height}}) \big) \leq \frac{1}{16}.
  \end{equation}
\end{lma}

\begin{proof}
  We first argue for the existence of $r_0\ge1$ and $\uc{c:box_height}>18$ for which~\eqref{e:cross_long} holds.
  Given $s\ge1$, set $A_0:=\cross(\tfrac{5}{4}s,s)$ and let $A_i$ denote the translate of $A_i$ along the vector $(2is,0)$.
  Translation invariance, duality and \eqref{e:easy_cross} imply that for every $i$,
  \begin{equation}
    \label{eq:5}
    \Pr_\lambda(A_i)= 1- \Pr_\lambda \big( \cross^\star \big(s, \b s \big) \big)\le  1-\uc{c:easy_cross}.
  \end{equation}
  If for some $k\ge1$ $\cross(2ks, s)$ occurs, then $A_i$ occur for all $i=0,1,\ldots,k-1$. Therefore,
  \begin{align}
    \label{eq:25}
    \Pr_\lambda\big( \cross(2ks,s)\big)\,\le\, \Pr_\lambda\Big(\bigcap_{i=0}^{k-1} A_i \Big)
    \,\le\, \Pr_\lambda\Big( \bigcap_{i=0}^{k-2}A_i \Big)
    (1-\uc{c:easy_cross}) +\rho_\lambda(s, \tfrac34 s).
  \end{align}
  Via induction in $k$ we obtain that
  \begin{equation}
    \label{eq:26}
    \Pr_\lambda\big( \cross(2ks,s)\big)\le
    (1-\uc{c:easy_cross})^k+\frac1{\uc{c:easy_cross}}
    \rho_\lambda(s,\tfrac34 s).
  \end{equation}
  For $k\ge4$ and $s_0\ge1$ large enough this expression is bounded by $1/16$ for all large $s\ge s_0$.
  Hence, the choice $r_0=6ks_0$ and $\uc{c:box_height}=6k\ge23$ guarantees that
  \begin{equation}
    \label{eq:27}
    \Pr_\lambda\big( \cross(\tfrac{\x}{3},
    \tfrac{\x}{\uc{c:box_height}})\big)\le \frac1{16}.
  \end{equation}
  for every $r\ge r_0$. With exactly the same proof, we can show that
  \eqref{eq:27} also holds with $\cross$ replaced by $\cross^\star$. Hence
  \eqref{e:cross_long} holds for this choice of $r_0$ and $\uc{c:box_height}$.

  To conclude the proof, we will show that any choice of $\uc{c:box_height}\ge 23$ implies
  \begin{equation}
    \label{eq:28}
    \inf_{r \ge r_0}  \mathbb{P}_\lambda \Bigg[ \boxHighTwoCross{\x/\uc{c:box_height}}{\tfrac{\x}3} \Bigg] > 0.
  \end{equation}
  To prove this, observe that
  \begin{equation}
    \begin{split}
      \mathbb{P}_\lambda \Bigg[
      \boxHighTwoCross{\x/\uc{c:box_height}}{\tfrac{\x}3} \Bigg] \,& \geq\,
      \Pr_\lambda \big( \cross(\tfrac{\x}{\uc{c:box_height}}, \tfrac{\x}{18})
      \big) \Pr_\lambda \big( \cross^\star(\tfrac{\x}{\uc{c:box_height}},
      \tfrac{\x}{18}) \big)-\rho_\lambda(\tfrac{r}{18}, \tfrac{r}{18})\\
      & \ge\, \uc{c:easy_cross}^2- \rho_\lambda(\tfrac{r}{18},
      \tfrac{r}{18}).
    \end{split}
  \end{equation}
  By \eqref{e:rho_to_zero}, this shows that the probability on the left hand side
  of \eqref{eq:28} is larger than $ \uc{c:easy_cross}^2/2$ for $r$ large enough, as required.
\end{proof}

The next event we will be interested in is
\begin{equation}
  \boxHighTwoMiddle[\dual]{\x}{h}{2a} := \bigcup_{q \in [0,h] \cap \mathbb{Q}} \Big[ \big( \{0\} \times [0,q] {\overset{\mathcal{O} \cap K}\longleftrightarrow} A \big) \cap \big( \{0\} \times [q, h] {\overset{\mathcal{V} \cap K}\longleftrightarrow} A \big) \Big],
\end{equation}
where $K = [0,\x] \times [0,h]$ and $A = \{r\} \times [h/2, h/2 + a]$.

\begin{definition}
  For every $r\ge r_0$, let $0 \le \alpha_\x\le r/2$ be such that
  \begin{equation}
    \label{e:alpha}
    \mathbb{P}_\lambda \Bigg[ \boxHighTwoMiddle[\dual]{\x/\uc{c:box_height}}{\x}{2 \alpha_\x} \Bigg] = \frac{\uc{c:prob_dual}}2.
  \end{equation}
\end{definition}

Let us briefly show that
\begin{equation}
  \label{e:alpha_height}
  \alpha_\x \leq \frac{\x}6.
\end{equation}
By Lemma~\ref{l:cross_continuous}, the quantity
\begin{equation}
  \label{eq:29}
   f(a)=\mathbb{P}_\lambda \Bigg[ \boxHighTwoMiddle[\dual]{\x/\uc{c:box_height}}{\x}{2a} \Bigg]
\end{equation}
is continuous  in $a$, and by Lemma~\ref{l:long_rectangle}, it satisfies $f(0)=0$ and $f(\tfrac r6)\ge
\uc{c:prob_dual}$. Hence, there exists $0\le \alpha_\x\le r/6$ such that $f(\alpha_\x)=\tfrac{\uc{c:prob_dual}}2.$

The next event is what we call the \emph{fork} and it is crucial in what follows
\begin{equation}
  \boxFork{\x}{h}{2a} := \bigg[
  \begin{array}{c}
  \big( \{0\} \times [0,h/2 - a] {\overset{\mathcal{O} \cap K}\longleftrightarrow} \{\x\} \times [h/2 + a, h] \big)\\
  \mcap \big( \{0\} \times [h/2 + a, h] {\overset{\mathcal{O} \cap K}\longleftrightarrow} \{\x\} \times [0,h/2 - a] \big)
  \end{array}
  \bigg],
\end{equation}
and analogously for the dual (where the picture is traced using dashed curves).

\begin{lma}
  \label{l:fork}
  Assuming \eqref{e:easy_cross}, for every $\x \geq r_0$, one of the following holds
  \begin{equation}
    \label{e:fork}
      \mathbb{P}_\lambda\bigg[ \boxFork{\x}{\x}{2 \alpha_\x} \bigg] \geq 2^{-13} \quad \text{or} \quad
      \mathbb{P}_\lambda\bigg[ \boxFork[\dual]{\x}{\x}{2 \alpha_\x} \bigg] \geq 2^{-13}.
  \end{equation}
\end{lma}

\begin{proof}
  By duality and invariance under right-angle's rotations, see \eqref{e:invariance},
  \begin{equation}
    \mathbb{P}_\lambda \big( \cross(\x, \x) \big) + \mathbb{P}_\lambda \big( \cross^\star(\x, \x) \big) = 1
  \end{equation}
  The two cases in the statement of the lemma will correspond to which of the above probabilities is at least one half.
  Therefore, we assume without loss of generality that
  \begin{equation}
    \label{e:cross_square}
    \mathbb{P}_\lambda \big( \cross(\x, \x) \big) \geq \frac 12.
  \end{equation}
  Using symmetry and then a simple union bound, we get
  \begin{equation}
    \frac 14 \leq \mathbb{P}_\lambda\bigg[ \boxSquareTickLabel{\x}{\x}{\tfrac{\x}2} \bigg] \leq \mathbb{P}_\lambda\bigg[ \boxSquareTwoTicks{\x}{\x}{2 \alpha_\x} \bigg] + \mathbb{P}_\lambda\bigg[ \boxSquareTwoCrossings{\x}{\x}{2 \alpha_\x} \bigg],
  \end{equation}
  where the event appearing in the last term of the above sum can be rigorously
  defined as \[ \boxSquareTwoCrossings{\x}{\x}{2 \alpha_\x}=
\mcup_{q \in [0, \x] \cap \mathbb{Q}} \big( \{0\} \times [0,q] {\overset{\mathcal{O} \cap K}\longleftrightarrow} A_2 \big) \cap \big( A_1(q) {\overset{\mathcal{V} \cap K}\longleftrightarrow} A_2 \big),\] with $K = [0, \x]^2$, $A_1(q) = \{0\} \times [q, \x] \cup [0, \x] \times \{\x\}$ and $A_2 = \{\x\} \times [\x/2, \x/2 + \alpha_\x].$

  Recalling that $\uc{c:box_height} > 1$ and \eqref{e:alpha_height}, we bound the second term in the above sum by a sum of two terms.
  They respectively correspond to whether the dual path depicted above stays confined in an $(r/\uc{c:box_height})$-wide rectangle or not.
  More precisely, recalling \eqref{e:alpha_height},
  \begin{equation}
    \begin{split}
      \mathbb{P}_\lambda\bigg[ \boxSquareTwoCrossings{\x}{\x}{2 \alpha_\x} \bigg] \,& \leq\, \mathbb{P}_\lambda\Bigg[ \boxSquareCorridor{\x}{\x}{2 \alpha_\x}{\x/\uc{c:box_height}}{\tfrac{\x}3} \Bigg] + \mathbb{P}_\lambda \Bigg[ \boxHighTwoMiddle[\dual]{\x/\uc{c:box_height}}{\x}{2 \alpha_\x} \Bigg]\\
      &\leq\, \frac{1}{16} + \frac{\uc{c:prob_dual}}{2} \,\leq\, \frac 18,
    \end{split}
  \end{equation}
  where we also have used~\eqref{e:cross_long}, \eqref{e:alpha} and that $\uc{c:prob_dual}<1/8$.
  Therefore
  \begin{equation}
    \mathbb{P}_\lambda\bigg[ \boxSquareTwoTicks{\x}{\x}{\alpha_\x} \bigg] \geq \frac 18.
  \end{equation}

  To finish the bound in \eqref{e:fork}, we are going to apply the FKG
  inequality, see \eqref{eq:3} above. We start by defining the events $A_1, A_2,
  A_3, A_4$ obtained by reflecting the above square along the vertical and
  horizontal axis. Of course the above bound will remain valid by the assumed
  symmetry of the system. Also, using the fact that the law is
  invariant with respect to rotations by right angles, the probability of
  observing a vertical crossing of the above box is at least one half, by
  \eqref{e:cross_square}. It is not difficult to see that the fork event occurs
  as soon as there is a vertical crossing of the box, together with the four
  events $A_1, A_2, A_3, A_4$. Plugging the above together with the
  FKG-inequality leads to the bound in \eqref{e:fork}.
\end{proof}

\nc{c:hit_middle}
\begin{lma}
  \label{l:hit_middle}
  Let $\tilde \alpha_{\uc{c:box_height}\x} := \alpha_{\uc{c:box_height}\x} \wedge \big( \bOverTwo \x \big)$.
  Then, for some $\uc{c:hit_middle} \in (0, \uc{c:easy_cross})$, we have
  \begin{equation}
    \label{e:hit_middle}
    \inf_{\x \geq r_0} \mathbb{P}_\lambda\bigg[ \boxThreeFourMiddle{\x}{\bb \x}{\tilde \alpha_{\uc{c:box_height} \x}} \bigg] \geq \uc{c:hit_middle} \quad \text{and} \quad \inf_{\x \geq r_0} \mathbb{P}_\lambda\bigg[ \boxThreeFourMiddle[\dual]{\x}{\bb \x}{\tilde \alpha_{\uc{c:box_height} \x}} \bigg] \geq \uc{c:hit_middle}.
  \end{equation}
  where the above event is defined as $\{0\} \times [0, 4\x/3] \overset{K}\leftrightarrow A$, with $K = [0,\x] \times [0,4\x/3]$ and $A = \{\x\} \times [2\x/3, 2\x/3 + \tilde{\alpha}_{\uc{c:box_height} \x}]$ (and analogously for the dual).
\end{lma}

\begin{proof}
  We only show the first inequality above, since the dual case is completely analogous.
  For this we split the proof in two cases.
  Either $\alpha_{\uc{c:box_height} \x} \geq \bOverTwo \x$, in which case $\tilde{\alpha}_{\uc{c:box_height} \x} = 2\x/3$ and the result follows from \eqref{e:easy_cross}, together with the vertical reflection symmetry of the system and a union bound.

  We now treat the case $\tilde \alpha_{\uc{c:box_height} \x} =
  \alpha_{\uc{c:box_height} \x} < \bOverTwo \x$. We first recall that
  \eqref{e:alpha} gives
  \begin{equation}
    \mathbb{P}_\lambda \Bigg[ \boxHighTwoMiddle{\x}{\uc{c:box_height}\x}{2 \alpha_{\uc{c:box_height} \x}} \Bigg] \geq \frac{\uc{c:prob_dual}}2.
  \end{equation}
  We now split the interval $\{r\} \times [\uc{c:box_height}/2,
  \uc{c:box_height}/2 + \alpha_{\uc{c:box_height} \x}]$ into eight equal parts.
  Then, using the union bound, there must exist some $h_r \in
  \{\uc{c:box_height}r/2 + i \alpha_{\uc{c:box_height} \x}/8; i = 0, \dots,
  7\}$, such that
  \begin{equation}
    \label{e:cross_to_I}
    \mathbb{P}_\lambda \Bigg[ \boxHighSingle{\x}{\uc{c:box_height}\x}{I_r} \Bigg] \geq \frac{\uc{c:prob_dual}}{16},
  \end{equation}
  where $I_r = \{r\} \times [h_r, h_r + \alpha_{\uc{c:box_height} \x}/8]$.

  Since we are assuming $\alpha_{\uc{c:box_height} \x} < 2\x/3$, the length of $I_r$ is no larger than $\x/12$.
  Thus
  \begin{equation}
    \begin{split}
      \mathbb{P}_\lambda \bigg[ \boxThreeFourMiddle{\x}{\bb \x}{\tilde \alpha_{\uc{c:box_height} \x}} \bigg] & \geq \mathbb{P}_\lambda \Bigg[ \boxThreeFourUp{\x}{\b \x} \mcap \boxThreeFourDown{\x}{\b \x} \mcap \boxHighSingle{\x}{\uc{c:box_height}\x}{I_r} \Bigg]\\
      & \overset{\textnormal{FKG}, \eqref{e:easy_cross}, \eqref{e:cross_to_I}}\geq \Big(\frac{\uc{c:easy_cross}}2\Big)^2 \frac{ \uc{c:prob_dual}}{16}\,:=\, \uc{c:hit_middle}
    \end{split}
  \end{equation}
  Where above we have made a slight abuse of notation in the second term above: The first two events in the triple intersection need to be translated vertically.
  This concludes the proof.
\end{proof}

We now define
\begin{equation}
  p_\x := \mathbb{P}_\lambda\bigg[ \; \boxFourThree{3\x/2}{\tfrac 43 \x} \bigg] \quad \text{and} \quad p_\x^\star := \mathbb{P}_\lambda\bigg[ \; \boxFourThree[\dual]{3\x/2}{\tfrac 43 \x} \bigg].
\end{equation}

\nc{c:empty}
\nc{c:hit_smaller}
\nc{c:rsw}
\begin{lma}
  \label{l:alpha_no_grow}
  If $\alpha_{\x/2} \ge \frac{2}{5 \uc{c:box_height}} \alpha_{\uc{c:box_height}
    \x}$ for some $r \geq 2 r_0$, then $\max\{p_\x, p_\x^\star\} \geq \uc{c:rsw} =
  \uc{c:rsw}(\uc{c:box_height}, \uc{c:hit_middle})$.
\end{lma}

\begin{proof}
  Using Lemma~\ref{l:fork} we can assume without loss of generality that
  \begin{equation}
    \label{e:choose_fork}
    \mathbb{P}_\lambda\bigg[ \boxFork{\x/2}{\x/2}{2 \smash{\alpha_{\x/2}}} \bigg] \geq 2^{-13}
  \end{equation}
  and then prove that $p_\x \geq \uc{c:rsw}$.
  If on the other hand the above bound holds for the dual, an identical proof shows that $p^\star_\x \geq \uc{c:rsw}$.

  Recall that we are assuming $\alpha_{\x/2} \geq \frac{2}{5
    \uc{c:box_height}} \alpha_{\uc{c:box_height} \x} \geq \frac{2}{5
    \uc{c:box_height}} \tilde{\alpha}_{\uc{c:box_height} \x}$. With this, we
  define the intervals $I_j = \{\x\} \times [(2/3) \x + (j-1) \alpha_{\x / 2},
  (2/3) \x + j \alpha_{\x / 2}]$, for $j = 1, \dots, \lceil 2 \uc{c:box_height}
  / 5 \rceil$, which cover the interval $A$ defined below \eqref{e:hit_middle}.
  Therefore, the union bound yields
  \begin{equation}
    \begin{split}
      \max_{j \leq \lceil 2 \uc{c:box_height}/5 \rceil} \mathbb{P}_\lambda & \bigg[ \boxThreeFourSmall{\x}{\bb \x}{I_j} \bigg]\, \geq\, \frac{1}{\lceil 2 \uc{c:box_height} / 5 \rceil} \sum_{j \leq \lceil 2 \uc{c:box_height} / 5 \rceil} \mathbb{P}_\lambda \bigg[ \boxThreeFourSmall{\x}{\bb \x}{I_j} \bigg]\\
      & \geq \,\frac{1}{\lceil 2 \uc{c:box_height} / 5 \rceil} \mathbb{P}_\lambda\bigg[ \boxThreeFourMiddle{\x}{\bb \x}{\tilde \alpha_{\uc{c:box_height} \x}} \bigg] \,\geq \,\frac{\uc{c:hit_middle}}{\lceil 2 \uc{c:box_height} / 5 \rceil} \,=:\, \uc{c:hit_smaller}.
    \end{split}
  \end{equation}

  Let $j_o$ be the index attaining the maximum in the left hand side of the above equation.
  Since the interval $I_{j_o}$ has length $\alpha_{\x / 2}$, it can be covered by the small interval appearing in \eqref{e:choose_fork}, such as illustrated in the picture below
  \begin{equation}
    p_\x \,\geq\, \mathbb{P}_\lambda \Biggl[ \boxUseFork{\tfrac{4}{3}\x}{\x/2}{2 \smash{\alpha_{\x/2}}}{(1/3)\x} \; \Biggr] \overset{\text{FKG}}\geq 2^{-13} \uc{c:hit_smaller}^2 \,=:\, \uc{c:rsw},
  \end{equation}
  finishing the proof of the lemma.
\end{proof}

\begin{remark}
  \label{r:c_rsw_small}
  It is important to observe that $\uc{c:rsw}$ is strictly smaller than $\uc{c:hit_middle}^2$, as one can clearly see from the definitions of  $\uc{c:hit_smaller}$ and $\uc{c:rsw}$.
  This will be used in Lemma~\ref{l:boost_rsw} below.
\end{remark}

\nc{c:tiny_rsw}
\begin{lma}
  \label{l:tiny_rsw}
  There exists a constant $\uc{c:tiny_rsw} > 0$ such that the following holds.
  Suppose that for some $r \geq r_0$, $p_\x \geq \uc{c:rsw}$ (respectively $p_\x^\star \geq
  \uc{c:rsw}$) and that for some $\xprime \geq 30 \x$ we have
  $\alpha_{\uc{c:box_height}\xprime} \leq \x$. Then $p_\xprime \geq
  \uc{c:tiny_rsw}$ (respectively $p_\xprime^\star \geq \uc{c:tiny_rsw}$).
\end{lma}

\begin{proof}
  Using the fact that $p_{\x} \geq \uc{c:rsw}$ and Inequalities~(\ref{item:13})
  and (\ref{item:16}) in Lemma~\ref{lem:standardInequalities}, we can conclude that
  \begin{equation}
    \label{e:build_circuit}
    \mathbb{P} \big( \cir(2 \x, 4 \x) \big) \geq \uc{c:rsw}^{124},
  \end{equation}
  We now recall \eqref{e:hit_middle} and our assumption that $\alpha_{\uc{c:box_height}\xprime} \leq \x$.
  Using symmetry and the FKG-inequality, we obtain
  \begin{equation}
    \label{e:glue_ring}
    p_{\xprime} \,=\, \Pr_\lambda\big(\cross(3 \xprime/2, 4 \xprime/3)\big) \,\geq\, \mathbb{P}_\lambda \Biggl[ \boxRings{2 \xprime}{\frac{4}{3}\xprime}{2 \x}{6 \x} \Biggr] \,\geq\, \uc{c:hit_middle}^2 \uc{c:rsw}^{124} \,=:\, \uc{c:tiny_rsw},
  \end{equation}
  finishing the proof of the lemma.
\end{proof}

We have obtained in the previous lemma a condition for $p_\xprime \geq \uc{c:tiny_rsw}$.
However, the constant is smaller than the original $\uc{c:rsw}$ that lower bounded $p_\x$ in the hypothesis of Lemma~\ref{l:tiny_rsw}.
The purpose of the next lemma is to use the above in order to bootstrap the lower bound of $p_\xprime$ back to $\uc{c:rsw}$.

\nc{c:boost_rsw}
\begin{lma}
  \label{l:boost_rsw}
  There exists a constant $\uc{c:boost_rsw} > 300$ such that the following holds.
  Suppose that for some $r \geq r_0$, $p_\x \geq \uc{c:rsw}$ (respectively $p_\x^\star$) and that for every $\xprime \in [30\x, \uc{c:boost_rsw} \x]$ we have $\alpha_{\uc{c:box_height}\xprime} \leq \x$.
  Then $p_{2 \uc{c:boost_rsw} \xprime} \geq \uc{c:rsw}$ (respectively $p_{2 \uc{c:boost_rsw} \xprime}^\star \geq \uc{c:rsw}$).
\end{lma}

\begin{proof}
  We wish to apply the circuit argument of Lemma~\ref{lma:circuit}. For this
  purpose let us fix a function $f=f_{\lambda,c}$ as in Item (\ref{item:11}) of
  Lemma~\ref{lma:circuit}, corresponding to the constant
  $c:=\uc{c:tiny_rsw}^{124}$. Recall from Remark~\ref{r:c_rsw_small} that
  $\uc{c:rsw} < \uc{c:hit_middle}^2$, so that we can choose $\uc{c:boost_rsw} >
  300$ so that
  \begin{equation}
    \label{e:choose_boost}
    1-f(\tfrac {300}{\uc{c:boost_rsw}}) \geq \frac{\uc{c:rsw}}{\uc{c:hit_middle}^2}.
  \end{equation}

  We only prove the statement for the dual quantities, the same proof also works
  for the primal ones (except that the function $f$ need to be chosen according
  to Item~(\ref{item:6}) in Lemma~\ref{lma:circuit}).

  Assume that $p_\x^\star \geq \uc{c:rsw}$ and that for
  every $\xprime \in [30\x, \uc{c:boost_rsw} \x]$ we have
  $\alpha_{\uc{c:box_height}\xprime} \leq \x$.
  Using Lemma~\ref{l:tiny_rsw}, we know that for every $\x'' \in [30 \x, \uc{c:boost_rsw}r]$, we have $p_{\x''} \geq \uc{c:tiny_rsw}$.
  By the standard inequalities of
  Lemma~\ref{lem:standardInequalities}, this implies in particular
  \begin{equation}
    \label{eq:30}
    \inf_{60r\le r''\le \uc{c:boost_rsw}r/10}\Pr_\lambda\big(\cir^\star(r'',2r'')\big)\ge
    \uc{c:tiny_rsw}^{124}.
  \end{equation}
  Therefore, using the circuit argument (\ref{item:11}) of
  Lemma~\ref{lma:circuit}, we obtain
  \begin{equation}
    \label{eq:31}
    \Pr_\lambda\big(\cir^\star(60r,\tfrac{\uc{c:boost_rsw}r}5)\big)\ge
    1-f(\tfrac{300}{\uc{c:boost_rsw}}) \overset{~\eqref{e:choose_boost}}{\ge}  \frac{\uc{c:rsw}}{\uc{c:hit_middle}^2}.
  \end{equation}

  We can finally proceed as in the end of the proof of Lemma~\ref{l:tiny_rsw},
  observing that $p_{2 \uc{c:boost_rsw} \x}^\star \geq \Pr_\lambda(\cross^\star(3
  \uc{c:boost_rsw} \x, \uc{c:boost_rsw} \x))$, which can be lower bounded by
  \begin{equation}
    \label{e:glue_boost}
    \mathbb{P}_\lambda\Bigg[ \boxRings[\dual]{3 \uc{c:boost_rsw}\x}{\uc{c:boost_rsw} \x}{60\x}{\uc{c:boost_rsw} \x/5} \Bigg] \geq \uc{c:hit_middle}^2 \frac{\uc{c:rsw}}{\uc{c:hit_middle}^2} = \uc{c:rsw},
  \end{equation}
  as required.
\end{proof}

\nc{c:rsw_blows}
\begin{lma}
  \label{l:alpha_blows}
  Fix $r \geq r_0$ and suppose that for some $\xprime \in [30 \x, \uc{c:boost_rsw} \x]$ we have $\alpha_{\uc{c:box_height} \xprime} \geq \x$, then $\max\{p_{\x''}, p^\star_{\x''}\} \geq \uc{c:rsw}$ for some $\x'' \in [\xprime, \uc{c:rsw_blows} \xprime]$.
\end{lma}

\begin{proof}
  Let us choose an integer $K = K(\uc{c:box_height}, \xprime/\x) \geq 1$ such that
  \begin{equation}
    \Big( \frac{5 \uc{c:box_height}}{2} \Big)^K \x \geq \uc{c:box_height} (2 \uc{c:box_height})^K \xprime
  \end{equation}
  and then set $\uc{c:rsw_blows} := 2 (2 \uc{c:box_height})^K.$

  Observe now that if for some $k = 0, \dots, K - 1$ we have
  \begin{equation}
    \label{e:alpha_growing}
    \alpha_{\uc{c:box_height}(2 \uc{c:box_height})^k \xprime} \geq \big( \tfrac{2}{5 \uc{c:box_height}} \big) \alpha_{\uc{c:box_height}(2 \uc{c:box_height})^{k+1} \xprime},
  \end{equation}
  we can apply Lemma~\ref{l:alpha_no_grow}, replacing $\x$ with $(2\uc{c:box_height})^{k+1} \xprime$.
  This will yield the bound
  \begin{equation}
    \max\{p_{(2 \uc{c:box_height})^{k+1} \xprime}, p^\star_{(2 \uc{c:box_height})^{k+1} \xprime}\} \geq \uc{c:rsw}.
  \end{equation}
  By our choice of $\uc{c:rsw_blows}$, we conclude that $\xprime \leq (2 \uc{c:box_height})^{k+1} \xprime \leq (2 \uc{c:box_height})^K \xprime \leq \uc{c:rsw_blows} \xprime$, as desired.

  On the other hand if \eqref{e:alpha_growing} does not hold for any $k = 0, \dots, K - 1$, then
  \begin{equation}
    \label{e:alpha_blowing}
    r \leq \alpha_{\uc{c:box_height} \xprime} \leq \big( \tfrac{2}{5 \uc{c:box_height}} \big) \alpha_{\uc{c:box_height}(2 \uc{c:box_height}) \xprime} \leq \dots \leq \big( \tfrac{2}{5 \uc{c:box_height}} \big)^{K} \alpha_{\uc{c:box_height}(2 \uc{c:box_height})^{K} \xprime},
  \end{equation}
  so that by our choice of $K$ we would have $\alpha_{\uc{c:box_height}(2 \uc{c:box_height})^{K} \xprime} \geq (5 \uc{c:box_height} / 2)^{K} r \geq \uc{c:box_height}(2 \uc{c:box_height})^{K} \xprime$, which is a contradiction with \eqref{e:alpha_height}.
\end{proof}

\nc{c:rsw_max}
\begin{lma}
  \label{l:rsw_max}
  There exists $\uc{c:rsw_max} > 0$ such that
  \begin{equation}
    \label{e:hard_crossing}
    \inf_{\x \geq 1} \big(\max\{p_r, p^\star_r\}\big) \geq \uc{c:rsw_max}.
  \end{equation}
\end{lma}

\begin{proof}
  We first claim that there exists an increasing sequence $\x_1\le
  \x_2\le \cdots$ such that
  \begin{enumerate}[\quad a)]
  \item $3\le \tfrac{\x_{i+1}}{\x_i}\le \max\{2\uc{c:boost_rsw}, \uc{c:rsw_blows}\}$ for every $i \geq 1$ and
  \item $\max\{p_{\x_i}, p^\star_{\x_i}\} \geq \uc{c:rsw}$, for every $i \geq 1$.
  \end{enumerate}
  We first construct $\x_1$. Since  $\alpha_r\le r/6$ for every $r\ge1$ (see \eqref{e:alpha_height}),  there
  must exist $k\ge 1$ such that
  \begin{equation}
    \alpha_{r_0(2 \uc{c:box_height})^k} \geq \big( \tfrac{2}{5 \uc{c:box_height}} \big) \alpha_{r_0(2 \uc{c:box_height})^{k+1}},
  \end{equation}
  Applying Lemma~\ref{l:alpha_no_grow} to $r_1:=2r_0(2\uc{c:box_height})^{k}$
  yields the bound
  \begin{equation*}
    \max\{p_{r_1}, p^\star_{r_1}\} \geq \uc{c:rsw}
  \end{equation*}
  as desired.

  To finish the construction of the $(\x_i)$'s, it is enough to show that
  \begin{display}
    \label{e:reduction_max}
    for any $\x \geq 1$ such that $p_\x \geq \uc{c:rsw}$, there exists $\x'' \in [3\x, \max\{2\uc{c:boost_rsw}, \uc{c:rsw_blows}\} \x]$\\ for which $\max\{p_{\x''}, p^\star_{\x''}\} \geq \uc{c:rsw}$.
  \end{display}
  We then split the proof of this claim into two cases, depending on whether or not
  \begin{display}
    \label{e:there_is_large_alpha}
    there exists $\xprime \in [3 \x, \uc{c:boost_rsw} \x]$ such that $\alpha_{\uc{c:box_height} \xprime} \geq \x$.
  \end{display}
  In this case, Lemma~\ref{l:alpha_blows} shows that the existence of $\x'' \leq \uc{c:rsw_blows} \x$ as in \eqref{e:reduction_max}.
  On the other hand, if \eqref{e:there_is_large_alpha} does not hold, than we can use $\x'' = 2\uc{c:boost_rsw}$, according to Lemma~\ref{l:boost_rsw}.
  In any case we have proved \eqref{e:reduction_max}, which shows the existence of $\x_1,\x_2, \dots$ as above.

  To end the proof we use the standard inequalities of
  Lemma~\ref{lem:standardInequalities} to interpolate between the values $\x_i$.
\end{proof}

We now have obtained a statement similar to the RSW result stated in Theorem~\ref{thm:rsw}.
However, we only know that the crossing probability is bounded below for the primal or the dual and not both at the same time.
It could be still the case that only the dual crossings have a bounded probability.
The purpose of the following lemma is to show that this cannot be the case.

\nc{c:rsw_both_1}
\nc{c:rsw_both_2}
\begin{lma}
  \label{l:rsw_both}
  There exist constants $\uc{c:rsw_both_1} \geq 10$ and $\uc{c:rsw_both_2}>0$ such that the following holds.
  Given $\x \geq r_0$, if $p_{\xprime}^\star \geq \uc{c:rsw_max}$ for every $\xprime \in [\x, \uc{c:rsw_both_1} \x]$, then $p_{\uc{c:rsw_both_1} \x} \geq \uc{c:rsw_both_2}$.
  Moreover, the same holds true if we swap the roles of $p^\star$ and $p$ above.
\end{lma}

\begin{proof}
  We wish to apply the circuit argument of Lemma~\ref{lma:circuit}. To this
  end, let us fix a function $f=f_{\lambda,c}$ as in Item (\ref{item:11}) of
  Lemma~\ref{lma:circuit}, corresponding to the constant
  $c:=\uc{c:rsw_max}^{124}$.
  We first choose $\uc{c:rsw_both_1} = \uc{c:rsw_both_1}(\uc{c:easy_cross},
  \uc{c:rsw_max}) \geq 30$ large enough such that
  \begin{equation}
    \label{e:choose_both_1}
    f(\tfrac{10}{\uc{c:rsw_both_1}}) \leq \frac{\uc{c:easy_cross}}{3}
  \end{equation}
  and $\uc{c:rsw_both_2}>0$ small enough that
  \begin{equation}
    \label{e:choose_both_2}
    \big(1 - \uc{c:rsw_both_2}^{1/(2\uc{c:rsw_both_1})}\big) \big( 1 - \frac{\uc{c:easy_cross}}{3} \big) \geq 1 - \frac{\uc{c:easy_cross}}{2}.
  \end{equation}

  Suppose now, for a contradiction, that $p_{\uc{c:rsw_both_1} \x} \leq \uc{c:rsw_both_2}$, therefore
  \begin{equation}
    \mathbb{P}_\lambda\big[\cross^\star \big(\tfrac{4}{3} \uc{c:rsw_both_1} \x, \tfrac{3}{2} \uc{c:rsw_both_1} \x \big) \big] = 1 - p_{\uc{c:rsw_both_1} \x} \geq 1 - \uc{c:rsw_both_2}.
  \end{equation}
  Cover the interval $[0, \tfrac{3}{2} \uc{c:rsw_both_1} \x]$ by at most $H =
  \lceil \tfrac{3}{2} \uc{c:rsw_both_1} \rceil \leq 2 \uc{c:rsw_both_1}$
  intervals of length $\x$. The square-root trick \eqref{eq:32} implies that
  there exists $h \in \{0, 1, \dots, H-1\}$, such that
  \begin{equation}
    \mathbb{P}_\lambda \bigg[ \boxThreeFourMiddle[\dual]{(4/3) \uc{c:rsw_both_1} \x}{\tfrac{3}{2} \uc{c:rsw_both_1} \x}{I_h} \bigg] \geq 1 - \uc{c:rsw_both_2}^{1/H},
  \end{equation}
  where $I_h = \{\tfrac{4}{3} \uc{c:rsw_both_1} \x\} \times [h \x, (h + 1) \x]$.
  We now recall that $p_{\xprime}^\star \geq \uc{c:rsw_max}$ for every $\xprime
  \in [\x, \uc{c:rsw_both_1} \x]$. Therefore, by the standard inequalities of
  Lemma~\ref{lem:standardInequalities}, we obtain in particular
  \begin{equation}
    \label{eq:33}
    \inf_{2r\le r'\le\uc{c:rsw_both_1} r/10}
    \Pr_\lambda\big(\cir^\star(r',2r')\big)\,\ge\, \uc{c:rsw_max}^{124}.
  \end{equation}
  Therefore, by the circuit argument (\ref{item:11}) in Lemma~\ref{lma:circuit},
  \begin{equation*}
    \mathbb{P} \Big( \cir^\star(2 \x, \uc{c:rsw_both_1} \x / 5) \Big) \,\geq\, 1-
    f(\tfrac{10}{\uc{c:rsw_both_1}}) \,\overset{\eqref{e:choose_both_1}}\geq\, 1 - \frac{\uc{c:easy_cross}}{3}.
  \end{equation*}
 As in \eqref{e:glue_boost}, we obtain
  \begin{align*}
      \mathbb P_\lambda\big( \cross^\star(\tfrac{8}{3} \uc{c:rsw_both_1} \x, \tfrac{3}{2} \uc{c:rsw_both_1} \x)\big)
      \,& \geq\, \mathbb{P}_\lambda\Bigg[ \boxRings{\tfrac{8}{3} \uc{c:rsw_both_1} \x}{\tfrac{3}{2} \uc{c:rsw_both_1} \x}{2 \x}{\uc{c:rsw_both_1}r/5} \Bigg]\\& \geq\, (1 - \uc{c:rsw_both_2}^{1/H})^2 \big(1 - \frac{\uc{c:easy_cross}}{3} \big) \,\overset{\eqref{e:choose_both_2}}\geq\, 1 - \frac{\uc{c:easy_cross}}{2}
  \end{align*}
  which is a contradiction with \eqref{e:easy_cross}.
\end{proof}

\begin{remark}
  In the proof of Lemmas~\ref{l:boost_rsw} and \ref{l:rsw_both} we have focused on the dual paths instead of primal to emphasize that we only needed the weak decoupling inequality \eqref{e:rho_to_zero} that has been derived from \eqref{eq:7}.
\end{remark}

\begin{proof}[Proof of Theorem~\ref{thm:rsw}]
  We first claim that it is enough to prove that for any $\x \geq r_0$ there exists
  some $\xprime \in [\x, \uc{c:rsw_both_1} \x]$ such that $p^\star_{\xprime}
  \geq \min\{\uc{c:rsw_max}, \uc{c:rsw_both_2}\}$. Indeed, if this is the case,
  we can use the standard inequalities in Lemma~\ref{lem:standardInequalities} to
  show that $\inf_{\x \geq 1} p^\star_{\x} > 0$. (The case $p_\x$ is
  completely analogous.)
  So, let us fix some $\x \geq r_0$ and assume that
  $p^\star_{\xprime} < \uc{c:rsw_max}$ for every $\xprime \in [\x,
  \uc{c:rsw_both_1} \x]$, then by
  Lemma~\ref{l:rsw_max}, we would have $p_{\xprime} \geq \uc{c:rsw_max}$ for
  every $\xprime \in [\x, \uc{c:rsw_both_1} \x]$. But then
  Lemma~\ref{l:rsw_both} implies that $p^\star_{\uc{c:rsw_both_1} \x} \geq
  \uc{c:rsw_both_2}$ as desired.
\end{proof}

\section{Sharpness of the phase transition}
\label{s:sharp}

We have so far showed that decay of correlations is sufficient for the two thresholds $\lambda_0$ and $\lambda_1$ to be non-trivial; see Section~\ref{s:finite-size}.
Moreover, we have shown that a distinctive critical behavior, namely the box-crossing property, occurs throughout the interval $[\lambda_0,\lambda_1]$; see Section~\ref{s:critical}.
Our goal for this section is to prove that these two thresholds $\lambda_0$ and $\lambda_1$ are in fact the same.

\begin{thm}\label{thm:equality}
  For Poisson Boolean percolation satisfying~\eqref{eq:7}, we have $\lambda_0 = \lambda_1$.
\end{thm}

In order to establish an equality between the two thresholds $\lambda_0$ and $\lambda_1$, our aim will be to show that the crossing probabilities grow very fast from close to zero to close to one.
More precisely, we have seen that the probability of crossing a fixed-ratio rectangle is bounded away from zero and one for $\lambda \in [\lambda_0, \lambda_1]$.
In this section, we show that this implies that the derivative of the crossing probability is large throughout the same interval.

The phenomenon of sharp thresholds is well understood in the context of Boolean functions on the discrete cube $\{0,1\}^n$ equipped with product measure.
In order to apply the tools from this theory in our context we will following an approach introduced in~\cite{ahlbrogrimor14}.
Thus, the first step of our proof will be to enhance the probability space in which we work to suit this setting, see Subsection~\ref{ss:two_stage}.

After having defined this alternative construction, we will employ Margulis-Russo's formula to relate the derivative of crossing probabilities to influences of these auxiliary random variables.
At this point, an inequality of Talagrand will become handy in order to bound the total influence from below.
This is done in two steps, where we bound the conditional variance of the crossing probabilities from below and the individual influences from above, see Subsections~\ref{ss:controlling_variance} and \ref{ss:controlling_influences} respectively.
Subsection~\ref{ss:proof_of_sharp} combines the above bounds to prove Theorem~\ref{thm:equality}.

\begin{remark}
Of the three main theoretical components of this paper, the one presented in this section is the least general one.
Although what we will present below does apply to some other models, it strongly uses the Poissonian nature of the process, see Section~\ref{s:other_models}.
\end{remark}

\subsection{Sharp thresholds for Boolean functions}
\label{ss:sharp_boolean}

It has been known for some time that monotone events involving a large number of independent Bernoulli variables typically exhibit sharp thresholds, in the sense that the probability of the event increases sharply as the parameter of the Bernoulli variables passes a certain value.
The first evidence thereof dates back to the study of random graphs by Erd{\H o}s and R\'enyi~\cite{erdren60}, and the proof that the critical probability for bond percolation on $\Z^2$ equals $1/2$ by Kesten~\cite{kesten80}.
A more general understanding of these phenomena has been pursued since in work
of Russo~\cite{russo82}, Bollob\'as and Thomason~\cite{boltho87}, and Kahn,
Kalai and Linial~\cite{kahkallin88} as well as elsewhere.
See e.g.~\cite{odonnell14} or~\cite{garste}.

For $p \in [0, 1]$ let $\textup{P}_p$ denote product measure of intensity $p$ on $\{0,1\}^n$.
A Boolean function $f : \{0, 1\}^n \to \{0, 1\}$ is said to be \emph{monotone} when $f(\eta) \le f(\eta')$ whenever $\eta \le \eta'$ coordinate wise.
For monotone Boolean functions, the map $p \mapsto \textup{P}_p(f = 1)$ is also monotone as a function of $p$.
Moreover, the rate of change of $\textup{P}_p(f = 1)$ is related to the concept of influences of the individual variables, as made explicit by the Margulis-Russo formula (see~\cite{russo1978note}),
\begin{equation}
  \label{eq:russo}
  \frac{d}{dp} \textup{P}_p(f = 1) \,=\, \sum_{i = 1}^n \Inf_i^p(f),
\end{equation}
where $\Inf_i^p(f) := \textup{P}_p\big(f(\eta) \neq f(\sigma_i\eta)\big)$ and $\sigma_i$ is the operator that flips the value of $\eta$ at position $i$.
In case $f(\eta) \neq f(\sigma_i\eta)$, we say that the bit $i$ is \emph{pivotal} for $f(\eta)$.

Against one's initial intuition, Russo~\cite{russo82} further showed that a uniform upper bound on the individual influences implies a lower bound of their sum, thus assuring a rapid increase of $\textup{P}_p(f=1)$.
This mystique has since become better understood with the work of Kahn, Kalai and Linial~\cite{kahkallin88} and its extensions, which shows that not all influences may be too small.
The following inequality, due to Talagrand~\cite{talagrand94}, illustrates this fact well and gives a precise formulation of Russo's initial observation: There exists a constant $\uc{c:talagrand} > 0$ such that for every Boolean function $f:\{0,1\}^n\to\{0,1\}$ we have
\nc{c:talagrand}
\begin{equation}
  \label{eq:talagrand}
  \sum_{i=1}^n \Inf_i^p(f) \,\ge\, \frac{\uc{c:talagrand}}{\log \big( \tfrac{2}{p (1 - p)} \big)} \, \Var_p(f) \log \Big( \frac{1}{\max\{\Inf_i^p(f)\}} \Big).
\end{equation}
Talagrand stated this inequality for monotone functions, but this restriction is not necessary.

From \eqref{eq:russo} and \eqref{eq:talagrand}, the strategy of the proof becomes intuitive: We first give an alternative construction of the process involving independent Bernoulli random variables, then we bound the variance of the crossing probabilities from below and the influence of individual bits from above.
These are the contents of the next subsections.

\subsection{A two-stage construction}
\label{ss:two_stage}

We will for the rest of this section work on an enlarged probability space $(\overline\Omega,\overline{\mathcal{M}},\overline\Pr_\lambda)$ in order to construct our process in two stages, where the second step employs independent coin flips in order to define the final configuration.
In this way we will be able to condition on the outcome of the first step in order to arrive at a situation where tools from the discrete setting applies.

Recall that $\Omega$ denotes the space of locally finite point measures on $\R^2 \times \R_+$.
We will by $\overline\Omega$ denote the space of locally finite point measures on $\R^2 \times \R_+ \times [0, 1] \times [0, 1]$, by $\overline{\mathcal{M}}$ we denote the corresponding Borel sigma algebra, and we let $\overline\Pr_\lambda$ denote the law of a Poisson point process defined on $(\overline\Omega,\overline{\mathcal{M}})$ with intensity $\lambda\,dx\,\mu(dz)\,du\,dv$.
In this section we will use a very dense set of points and the last two coordinates $u$ and $v$ in the above construction will be used to perform a thinning of the Poisson process.

Throughout this section we will suppose that $m\ge1$ is a large integer and that $\lambda=2m\lambda_1$ is fixed. Given $\overline\omega\in\overline\Omega$ we denote by $\omega$ the projection of $\overline\omega$ onto $\Omega$. For $p\in[0,1]$ and $m\ge1$ we define
$$
\omega_p:=\sum_{\substack{(x,z,u,v)\in\overline\omega\\ u\le p,\, v\le1/m}}\delta_{(x,z)}\quad\text{and}\quad\mathcal{O}_p:=\bigcup_{(x,z)\in\omega_p} B(x,z).
$$
The usual properties of a Poisson point process imply that $\omega$ and $\omega_p$ are Poisson point processes on $\Omega$ distributed according to $\Pr_\lambda$ and $\Pr_{p\lambda/m}$, respectively.
Notice that $\omega_p$ is a $(p/m)$-thinning of $\omega$, and therefore it can be seen as a percolation process on $\{0,1\}^\omega$ with density $p/m$.
Since $\lambda=2m\lambda_1$, the above construction comes with two parameters, $p$ and $m$, of which we think of $m$ as large and fixed while $p$ will be varying. We have for this reason suppressed $m$ in the above notation, and apply a subscript $p$ to indicate that the events $\cross_p$, $\cir_p$ and $\arm_p$ refer to crossings in $\mathcal{O}_p$.

We will study the conditional probability
$$
Z_p:=\overline\Pr_\lambda\big(\cross_p(4r,r)\big|\omega\big).
$$
The product structure of the above conditional measure implies that for $\Pr_\lambda$-almost every $\omega\in\Omega$ the variable $Z_p$ can be viewed as the expectation of some Boolean function $f_\omega:\{0,1\}^\omega\to\{0,1\}$.

While $\omega$ has infinite support, it follows from~\eqref{e:finite_in_compact} that $\overline\Pr_\lambda$-almost surely, only finitely many balls touch the box $[-4r, 4r]^2$.
Therefore the function $f_\omega$ can be viewed as having only finitely many variables that we index by $i \in \{1, \dots, i_0\}$ ($i_0$ is the number of balls touching $[-4r, 4r]^2$).

In this conditional setting, the tools from the discrete analysis will apply, in particular the Margulis-Russo formula and Talagrand's inequality.
With the use of these we will be able to bound the gap between $\lambda_0$ and $\lambda_1$ by controlling the conditional variance of $f_\omega$ and its influences, defined, for $\overline\Pr_\lambda$-almost every $\overline\omega\in\Omega$, as
\begin{equation}\label{e:cond_inf}
\Inf_i^p(f_\omega):=\overline\Pr_\lambda\big(f_\omega(\omega_p)\neq f_\omega(\sigma_i\omega_p)\big|\omega\big),
\end{equation}
where $\sigma_i\omega_p=\omega_p+(1-2\ind_{\{u_i\le p,\,v_i\le1/m\}})\delta_{(x_i,z_i)}$.

\subsection{Quenched decoupling}

In order to handle dependencies in the quenched setting we will need an alternative to the decoupling inequality proved in Proposition~\ref{c:decouplef1f2}.
Given $\overline\omega\in\overline\Omega$, this can be done via an estimate on the \emph{influence region} of each point in the support of $\omega$. Since we will only consider events defined in terms of $\mathcal{O}_p$ on $[-4r,4r]^2$ it will suffice to consider points intersecting this region.
For $r\ge3n\ge3$, let $G = G(r,n)$ denote the event
\begin{equation}
  \label{e:region_of_inf_small}
  G := \Big[ \max\big\{z:(x, z)\in\supp(\omega)\text{ and } B(x, z) \cap [-4r, 4r]^2 \neq \emptyset \big\} \le r/3n \Big].
\end{equation}

\begin{lma}
  \label{l:region_of_inf_small}
  Assume the second moment condition~\eqref{eq:7}. For every $n\ge1$, $\lim_{\x \to \infty} \overline\Pr_\lambda (G) = 1$.
\end{lma}

\begin{proof}
  The proof follows the same steps as that of Lemma~\ref{l:decouple}.
\end{proof}

A first consequence of the above observation concerns quenched decouplings of events defined in well separated regions.
To state this precisely, fix integers $n, J \geq 1$ and let $D_1, \dots, D_J \subseteq [-4 \x, 4 \x]^2$ be measurable sets such that
\begin{equation}
  \label{e:A_j_separated}
  d(D_i, D_j) \geq \frac{\x}{n}, \text{ for every $i \neq j$}.
\end{equation}
Having these well separated sets, our aim is to decouple what happens inside each of them.
For this, fix some events $A_1, \dots, A_J$ that only look inside the sets $D_j$, that is $A_j \in \sigma(Y_x, x \in D_j)$.
Observe that on the event $G = G(\x, n)$ the events $A_j$ depend on disjoint subsets of $\omega$ under the conditional measure $\overline\Pr_\lambda(\,\cdot\,|\omega)$.
This implies that
\begin{equation}
  \label{e:f_decouple_in_G}
  \overline\Pr_\lambda \big( A_1 \cap \dots \cap A_J \big| \omega \big) \,\leq\, \ind_{G^c}(\omega) + \prod\nolimits_{j \leq J} \overline\Pr_\lambda(A_j | \omega).
\end{equation}

For us it is important to observe that not only the events $A_j$ decouple on $G(\x, n)$, but also their conditional probabilities $\overline\Pr_\lambda (A_j | \omega)$.
To see why this is true, let us first decompose the point measure $\omega$ as $\omega_1 + \dots + \omega_J + \tilde{\omega}$, where $\omega_j = \sum_i \delta_{(x_i, z_i)}$ with the sum ranging over $i \geq 0$ such that $d(x_i, D_j) \leq \x/(3n)$ and $\tilde{\omega}$ stands for the remainder.

Observe that on $G$, the conditional probabilities $\overline\Pr_\lambda(A_j | \omega)$ depend on $\omega$ through $\omega_j$ only.
In particular, we may exhibit a function $g_j$ such that $[\ind_G \cdot \overline\Pr_\lambda(A_j | \omega)](\omega) = \ind_G \cdot g_j(\omega_j)$, and by the Poissonian nature of $\omega$, all the $\omega_j$'s are independent.
As a consequence, for any $\delta > 0$,
\begin{equation}
  \label{e:f_given_omega_decouples}
  \begin{split}
    \overline\Pr_\lambda \Big(G \cap \bigcap_{j \leq J} \big[ \overline\Pr_\lambda (A_j | \omega) \leq \delta \big] \Big) \,& \le\, \prod_{j \leq J} \overline\Pr_\lambda \big[ g_j(\omega_j) \leq \delta \big]\\
    & \leq\, \prod_{j \leq J}\bigg[ \overline\Pr_\lambda \big( \overline\Pr_\lambda(A_j | \omega) \leq \delta \big)+\overline\Pr(G^c)\bigg].
  \end{split}
\end{equation}
Equations \eqref{e:f_decouple_in_G} and \eqref{e:f_given_omega_decouples} give a more precise version of the decoupling of Section~\ref{s:preliminaries} and will be important in bounding the conditional variance and influences in the next two sections.

\subsection{Controlling the variance}
\label{ss:controlling_variance}

In Section~\ref{s:critical} we showed that the probability of crossing a rectangle is non-degenerate throughout the critical regime $[\lambda_0, \lambda_1]$.
In the notation used in this section, that is to say that $\overline\E_\lambda[Z_p]\in(c,1-c)$ for some $c>0$ and all $r\ge1$ and $p\lambda / m \in \big[ \lambda_0, \lambda_1 \big]$.
Recall that $\lambda = 2 m \lambda_1$.
Our goal for this section is to show that this holds with probability close to one also for $Z_p$.

\nc{c:lambda}
\begin{prop}\label{prop:q-bound}
Assume the second moment condition~\eqref{eq:7}. Then, for every $\eps>0$ there exists $\delta>0$ and $\uc{c:lambda} = \uc{c:lambda}(\eps) > 0$ such that for every $m\ge\uc{c:lambda}$ we have
$$
\overline\Pr_\lambda\bigg(Z_p\in(\delta,1-\delta)\text{ for all }p\in\big[\tfrac{\lambda_0}{2\lambda_1},\tfrac{\lambda_1}{2\lambda_1}\big]\bigg)\ge1-\eps\quad\text{for all }r\ge1.
$$
The same holds with $\cross_p(4\x, \x)$ replaced by either $\cross_p^\star(4\x, \x)$, $\cir_p(\x, 2\x)$ or $\cir_p^\star(\x, 2\x)$.
\end{prop}

An immediate consequence of the above proposition is that
$$
\overline\Pr_\lambda\bigg(\Var_\lambda(\ind_{\cross_p(4r,r)}|\omega)\in(\delta,1-\delta)\text{ for all }p\in\big[\tfrac{\lambda_0}{2\lambda_1},\tfrac{\lambda_1}{2\lambda_1}\big]\bigg)\ge1/2\quad\text{for all }r\ge1.
$$

The first step in the proof of Proposition~\ref{prop:q-bound} is the following estimate on the effect of observing $\omega$ in determining whether there is a crossing of $[0,4r] \times [0,r]$ in $\mathcal{O}_p$.

\begin{lma}\label{lma:q}
For every $m\ge1$ and $p\in[0,1]$, and every measurable event $A$, we have
$$
\Var_\lambda\big(\overline\Pr_\lambda(A|\omega)\big)\le1/m.
$$
\end{lma}

\begin{proof}
Recall that $\omega_p$ is obtained as a subset of $\omega$ based on the value of the third and fourth coordinates of each point in $\overline\omega$. Alternatively, we could obtain $\omega_p$ in two stages: First, partition $\overline\omega$ into $\omega_p^1,\omega_p^2,\ldots,\omega_p^m$ and some remainder $\tilde{\omega}$, where $\omega_p^j$ contains all points with third coordinate at most $p$ and fourth coordinate in $\big[\frac{j-1}{m},\frac{j}{m}\big)$. Second, use auxiliary randomness, independent from everything else, to determine which of the $\omega_p^j$'s to choose as the set $\omega_p$. Since each of the $\omega_p^j$'s have the same distribution, the result would be the same.

Let us assume that this construction has been made, and that $J$ is the auxiliary random variable, uniformly distributed on $\{1,2,\ldots,m\}$, that determines which of the $\omega_p^j$'s that is chosen in the second step. Let $\Fc$ denote the sigma algebra containing the information of the partitioning $\omega_p^1, \ldots, \omega_p^m$. Since revealing the partitioning can only increase the variance, we have that
$$
\Var_\lambda\big(\overline\Pr_\lambda(A|\omega)\big)\le\Var_\lambda\big(\overline\Pr_\lambda(A|\Fc)\big).
$$
Moreover, by conditioning on the outcome of $J$,
$$
\overline\Pr_\lambda(A|\Fc)=\frac1m\sum_{j=1}^m\ind_{\{\omega_p^j\in A\}}.
$$
Since the configurations $\omega_p^1,\omega_p^2,\ldots,\omega_p^m$ are independent, we find that
$$
\Var_\lambda\big(\overline\Pr_\lambda(A|\omega)\big)\,\le\,\Var_\lambda\bigg(\frac1m\sum_{j=1}^m\ind_{\{\omega_p^j\in A\}}\bigg)\,\le\,\frac1m,
$$
as required.
\end{proof}

\begin{proof}[Proof of Proposition~\ref{prop:q-bound}]
Note that $Z_p$ is monotone in $p$, so it will suffice to consider the extremes of the interval $\big[\tfrac{\lambda_0}{2\lambda_1},\tfrac{\lambda_1}{2\lambda_1}\big]$.
Due to Corollary~\ref{cor:box-crossing} there exists a constant $\uc{c:q-bound} = \uc{c:q-bound}(4) > 0$ such that $\overline\E_\lambda[Z_p]\in(\uc{c:q-bound},1-\uc{c:q-bound})$ for all $p\in\big[\tfrac{\lambda_0}{2\lambda_1},\tfrac{\lambda_1}{2\lambda_1}\big]$. So, by Chebyshev's inequality and Lemma~\ref{lma:q}, for $b = 0, 1$,
$$
\overline\Pr_\lambda\Big(Z_{\lambda_b/(2\lambda_1)}\not\in(\uc{c:q-bound}/2,1-\uc{c:q-bound}/2)\Big)\,\le\,\frac{\Var_\lambda(Z_p)}{(\uc{c:q-bound}/2)^2}\,\le\,\frac{4}{m\uc{c:q-bound}^2},
$$
which is at most $\eps/2$ for $m\ge 8/(\eps\uc{c:q-bound}^2)$. Now, the result follows via a union bound.
\end{proof}

\subsection{Controlling the influences}
\label{ss:controlling_influences}

The next step is to obtain an upper bound on the individual conditional influences, defined as in~\eqref{e:cond_inf}, allowing for Talagrand's inequality to give a lower bound on the total conditional influence.
The result proved here is the following.

\begin{prop}
  \label{p:max_inf_small}
  Assume the second moment condition~\eqref{eq:7}.
  There exists a constant $\uc{c:lambda2}$ such that for every $m\ge\uc{c:lambda2}$ and $\eta>0$ we have
  \begin{equation}
    \overline\Pr_\lambda \Big( \max_{i \geq 1} \big( \Inf^p_i (f_\omega) \big) > \eta\text{ for some }p\in\big[ \tfrac{\lambda_0}{2\lambda_1}, \tfrac{\lambda_1}{2\lambda_1} \big] \Big) \leq \eta,
  \end{equation}
  for all $\x$ large enough depending only on $\eta$.
\end{prop}

This will be done in two steps.
First we relate the influence of a point to arm probabilities.
Then we bound the probability of the latter.
Before stating the first lemma, consider a partitioning of the rectangle $[0, 4\x] \times [0, \x]$ into $36 n^2$ smaller squares of side length $\x / 3n$.
We write $\arm_p^\ell(\tfrac{r}{n}, r)$, with $\ell = 1, \dots, 36 n^2$, for the arm event defined as in \eqref{e:def_arm}, but centered around each of these boxes.

These arm events help us to bound the influences $\Inf^p_i(f_\omega)$ as follows.
Assume that $i$ is pivotal, and the corresponding ball has radius smaller or equal to $r / 3n$.
Then there must exist some $\ell$ for which the event $\arm_p^\ell(r/n,r)$ must occur.
Therefore,

\begin{equation}
  \label{e:then_arm}
  \Inf_i^p(f_\omega)\le\ind_{G^c}(\omega) + \overline\Pr_\lambda\big(\arm_p^\ell(\tfrac{r}{n},r)\big|\omega\big)\quad \text{for some }\ell=1,2,\ldots, 36 n^2,
\end{equation}
for $\overline\Pr_\lambda$-almost every $\overline\omega$.

Estimating influences will now boil down to estimating conditional arm probabilities.

\nc{c:lambda2}
\begin{lma}\label{lma:q-arm}
  Assume the second moment condition~\eqref{eq:7}.
  There exist $\gamma > 0$ and $\uc{c:lambda2} > 0$ such that for every $m > \uc{c:lambda2}$ and $n \ge 16$ we have
  \begin{equation}
    \overline\Pr_\lambda \Big[ \overline\Pr_\lambda \big( \arm_p\big( \tfrac{\x}{n}, \x \big) \big| \omega \big) > n^{-\gamma} \text{ for some }p\in\big[\tfrac{\lambda_0}{2\lambda_1}, \tfrac{\lambda_1}{2\lambda_1} \big]\Big] < n^{-100},
  \end{equation}
  for all sufficiently large $\x$ depending only on $n$.
\end{lma}

\begin{proof}
  The proof follows the same structure of that of Lemma~\ref{lma:circuit}.
  For $r\ge n\ge16$ let $\ell_i=4^i\frac{r}{n}$.
  Note that if $\cir_p^\star(\frac{r}{n},r)$ fails to occur, then $A_i = \cir_p^\star(\ell_i, 2\ell_i)$ must also fail for every $i = 0, 1, \ldots, k - 1$, where $k = \lfloor \frac{1}{2} \log_4 n \rfloor$ ($k \ge 1$ as we assume $n \ge 16$).
  By~\eqref{e:f_decouple_in_G} we therefore have
  \begin{equation}\label{eq:decoup}
    \overline\Pr_\lambda\big(\arm_p(\tfrac rn, r) \big| \omega \big)\, \le\, \ind_{G^c} (\omega) + \prod_{i = 0}^{k - 1} \overline \Pr_\lambda(A_i^c | \omega).
  \end{equation}

  Fix $\eps>0$ so that $2^k\eps^{k/2}\le\frac12n^{-100}$. By Proposition~\ref{prop:q-bound} there exists $\delta>0$ and $\uc{c:lambda}(\varepsilon) > 0$ large enough so that for every $r \ge 1$ and $m\ge\uc{c:lambda}$ we have
  \begin{equation}\label{eq:quenched_cir}
    \overline\Pr_\lambda \Big( \overline\Pr_\lambda \big( \cir_p^\star(r, 2r) \big| \omega \big) > \delta \text{ for every }p\in\big[\tfrac{\lambda_0}{2\lambda_1}, \tfrac{\lambda_1}{2\lambda_1} \big]\Big) > 1 - \eps.
  \end{equation}
  Fix $\gamma > 0$ (independent of $n$) so that $n^{-\gamma} \ge (1-\delta)^{k/2}$.
  By~\eqref{eq:decoup} we have, on $G$, that if $\overline\Pr_\lambda(A_i|\omega)>\delta$ for at least half of the indices $i=0,1,\ldots,k-1$, then
  $$
  \overline\Pr_\lambda\big(\arm_p(\tfrac{r}{n},r)\big|\omega\big)\,\le\,(1-\delta)^{k/2}\,\le\, n^{-\gamma}.
  $$
In particular, via a union bound,
\begin{equation}
\overline\Pr_\lambda\Big(G\cap\big[\overline\Pr_\lambda\big(\arm_p(\tfrac rn,r)\big|\omega\big)>n^{-\gamma}\big]\Big)\,\le\,2^k\sup_{I}\overline\Pr_\lambda\Big(G \cap \bigcap_{i \in I}\big[\overline\Pr_\lambda(A_i|\omega)\le\delta\big]\Big),
\end{equation}
where the supremum is taken over all $I\subseteq\{0,1,\ldots,k-1\}$ of size at least $k/2$, of which there are at most $2^k$. By~\eqref{e:f_given_omega_decouples} and~\eqref{eq:quenched_cir}, this is bounded by
$$
2^k\sup_{I}\prod_{i\in I}\Big[\overline\Pr_\lambda\big(\overline\Pr_\lambda(A_i|\omega)\le\delta\big)+\overline\Pr_\lambda(G^c)\Big]\,\le\,2^k\big(\eps+\overline\Pr_\lambda(G^c)\big)^{k/2}.
$$
For large values of $r$, depending on $n$, this probability is at most $n^{-100}$ by the choice of $\eps>0$.
\end{proof}

We have now gathered all the pieces necessary to prove the main result of this subsection.

\begin{proof}[Proof of Proposition~\ref{p:max_inf_small}]
  We start by choosing $n$ large enough so that $36 n^{-98} \leq \eta/2$. Then,
  $$
   \overline\Pr_\lambda \Big( \max_{i \geq 1} \big( \Inf_i^p(f_\omega) \big) > n^{-\gamma} \text{ for some }p \in \big[ \tfrac{\lambda_0}{2\lambda_1}, \tfrac{\lambda_1}{2\lambda_1} \big]\Big)
  $$
  may via \eqref{e:then_arm} be bounded from above by
  \begin{equation*}
     \overline\Pr_\lambda \big( G^c \big) + \overline\Pr_\lambda \Big(\overline\Pr_\lambda \big( \arm_p^\ell\big( \tfrac{\x}{n}, \x \big) \big| \omega \big) > n^{-\gamma}\text{ for some }p \in \big[ \tfrac{\lambda_0}{2\lambda_1}, \tfrac{\lambda_1}{2\lambda_1} \big]\text{ and }\ell\le 36 n^2 \Big).
  \end{equation*}
  By Lemmas~\ref{l:region_of_inf_small} and \ref{lma:q-arm}, this is smaller than $\eta$ by our choice of $n$, once $r$ is taken large enough.
\end{proof}

\subsection{Proof of Theorem~\ref{thm:equality}}
\label{ss:proof_of_sharp}

Roughly speaking, in order to prove Theorem~\ref{thm:equality}, we are going to show that for every $\eps > 0$
\begin{equation}\label{eq:q-sharp}
  Z_{\lambda_1/(2\lambda_1)} - Z_{\lambda_0/(2\lambda_1)} \ge \frac{\lambda_1 - \lambda_0}{\eps},
\end{equation}
with positive probability.
Since $Z_p(\omega) \in [0, 1]$ almost surely, we must have that $\lambda_0=\lambda_1$.

\begin{proof}[Proof of Theorem~\ref{thm:equality}]
  We first give estimates on the derivatives of $Z_p$. Recall that, for $\overline\Pr_\lambda$ almost every $\overline\omega\in\overline\Omega$, $Z_p$ coincides with the expectation, with respect to the conditional measure $\overline\Pr_\lambda(\,\cdot\,|\omega)$, of some function $f_\omega:\{0,1\}^\omega\to\{0,1\}$. Moreover, with probability one, only finitely many points in $\omega$ will affect the outcome of events defined on $[-4r,4r]^2$. Hence, the domain of $f_\omega$ is for $\overline\Pr_\lambda$-almost every $\overline\omega$ finite dimensional, and the Margulis-Russo formula gives that
  \begin{equation}
    \label{e:dZ_dp}
    \frac{d Z_p}{d p} \,=\, \frac{d Z_p (\omega)}{d p} \,=\, \frac{1}{m}\sum_{i \geq 0} \Inf_i^p (f_\omega),
  \end{equation}
  where $\Inf_i^p (f_\omega)$ is as defined in~\eqref{e:cond_inf}, and the additional factor $1/m$ comes from the chain rule.

  Let us now assume that $\lambda_0 < \lambda_1$ and fix $m\ge \uc{c:lambda}(1/2) \vee \uc{c:lambda2}$.
  Combining~\eqref{e:dZ_dp} with the mean-value theorem and Talagrand's inequality \eqref{eq:talagrand} we obtain, for some $q\in\big[\frac{\lambda_0}{2\lambda_1},\frac{\lambda_1}{2\lambda_1}\big]$,
  \begin{equation*}
    \frac{Z_{\lambda_1/(2\lambda_1)} - Z_{\lambda_0/(2\lambda_1)}}{\lambda_1 - \lambda_0}\, =\, \frac{1}{2\lambda_1} \frac{d Z_p(\omega)}{d p} (q)\, \geq\, \frac{c}{2 \lambda_1 m \log(m)} \Var_\lambda(f_\omega | \omega) \log \bigg( \frac{c}{\max_i \big( \Inf_i^q (f_\omega) \big) } \bigg).
  \end{equation*}
  From Propositions~\ref{prop:q-bound} and \ref{p:max_inf_small}, we conclude that for $\x$ large enough (depending only on $\lambda_1 - \lambda_0$) the right hand side of the above equation is larger than $2/(\lambda_1 - \lambda_0)$ with positive probability.
  This contradicts the fact that $Z_p \in [0,1]$, thus finishing the proof of the theorem.
\end{proof}

\section{Continuity of the critical parameter}
\label{s:continuity}

We have in the previous section seen that $\lambda_0=\lambda_1$, and thus coincides with $\lambda_c$, under the assumption of finite second moment of the radii distribution.
In this section we will investigate the dependence of this critical parameter with respect to the radii distribution $\mu$, and write $\lambda_c=\lambda_c(\mu)$ throughout this section in order to emphasize this dependence.

\begin{thm}\label{thm:continuity}
  Let $(\mu_m)_{m\ge1}$ be a sequence of radii distributions, all dominated by some distribution $\nu$ with finite second moment.
  If $\mu_m\to \mu$ weakly, then $\lambdaoccupied(\mu_m)\to\lambdaoccupied(\mu)$.
\end{thm}

\begin{remark}
A weaker version of Theorem~\ref{thm:continuity}, proved under the stronger assumption of uniformly bounded support, has been known since long (see~\cite[Theorem~3.7]{meeroy96}). The argument given here is significantly simpler and shorter, but is in contrast to the aforementioned proof distinctively two-dimensional.
\end{remark}

The convergence may fail when the assumption of uniformly bounded tails is relaxed. For a first example of this, let $(\mu_m)_{m\ge1}$ be a sequence with infinite second moments converging weakly to a point mass at 1.
The critical parameter for each $\mu_m$ is degenerate (meaning that $\lambda_0=\lambda_1=0$), and cannot converge to the critical value of the point mass, which is strictly positive.

In the above example the convergence fails for rather obvious reasons. A second example, which gives further insight to what can go wrong, may be obtained from a cumulative distribution function $F$ of some distribution with finite second moment, when $F_m$ is given by
\begin{equation*}
  F_m(x):=(1-\tfrac{1}{m})F(x)+\tfrac{1}{m}F(\tfrac{x}{m}).
\end{equation*}
Each $F_m$ has finite second moment, but there is no uniform upper bound on its value. The dilation of $F$ by a factor $\frac1m$ corresponds to a scaling of the radii by $m$, and thus has an inverse quadratic effect on the critical parameter in that $\lambdaoccupied(F(\frac{x}{m}))=\frac{1}{m^2}\lambdaoccupied(F(x))$. So, although $\lambdaoccupied(\mu)>0$, we have
\begin{equation*}
  \lambdaoccupied(F_m)\,\le\,\lambdaoccupied\big((1-\tfrac{1}{m})\ind_{\{x\ge0\}}+\tfrac{1}{m}F(\tfrac{x}{m})\big)\,=\,m\,\lambdaoccupied\big(F(\tfrac{x}{m})\big)\,=\,\tfrac{1}{m}\lambdaoccupied(F).
\end{equation*}
Although $F_m\to F$ weakly, the critical parameters diverge.

\begin{proof}[Proof of Theorem~\ref{thm:continuity}]
  Fro the proof, we will write $F$, $F_m$ and $H$ for the cumulative distribution functions of $\mu$, $\mu_m$ and $\nu$. We first observe that if $F_m\ge H$ for all $m\ge1$ and $F_m\to F$, then also $F\ge H$ and $F$ has finite second moment too. To see this assume the contrary, in which case $F(x)<H(x)$ for some $x\ge0$. Either $x$ is a continuity point of $F$, or right continuity of $F$ and $H$ allows us to find a continuity point $x'$ of $F$ for which $F(x')<H(x')$ remains. But, in that case $\liminf_{m\to\infty}F_m(x')\ge H(x')>F(x')$, contradicting the assumed convergence. That $F$ has finite second moment in particular implies that $\lambdaoccupied(F)\in(0,\infty)$.

  We will base the proof on Propositions~\ref{prop:perturbation} and~\ref{prop:finite-size}.
  We therefore fix $\theta>0$ and $r_0 = r_0(\lambdaoccupied(F) + 1)$ accordingly.
  Write $\Pr_\lambda^{F_m}$ for the measure at intensity $\lambda$ and with radii distributed as $F_m$.
  To complete the proof, it will suffice to show that for every $\eps\in(0,1)$ there exists $\x\ge r_0$ and $m_0$ such that for all $m\ge m_0$
  \begin{equation}\label{eq:convergence}
  \Pr_{\lambdaoccupied(F)-\eps}^{F_m}\big(\cross(\x,3\x)\big)<\theta\quad\text{and}\quad\Pr_{\lambdaoccupied(F)+\eps}^{F_m}\big(\cross(3\x,\x)\big)>1-\theta.
  \end{equation}
  Proposition~\ref{prop:finite-size} then implies that $\lambdaoccupied(F_m)\in[\lambdaoccupied(F)-\eps,\lambdaoccupied(F)+\eps]$ for all $m\ge m_0$.

  Fix $\eps>0$.
  By Corollary~\ref{cor:nontriviality} and Theorem~\ref{thm:equality}, we may choose $\x$ large so that
  $$
  \Pr_{\lambdaoccupied(F) - \eps}^F\big(\cross(\x, 3\x)\big) < \theta/4 \quad \text{and} \quad \Pr_{\lambdaoccupied(F) + \eps}^F\big(\cross(3\x, \x)\big) > 1 - \theta/4.
  $$
  Next, based on~\eqref{e:decouple1} of Lemma~\ref{l:decouple}, fix $K\subseteq\R^2$ such that
  $$
  \Pr_{\lambdaoccupied(F)+\eps}^H\big(\mathcal{O}_{K^c}\cap[0,\x]^2=\varnothing\big)>1-\theta/4.
  $$
  And, choose $\delta>0$ such that the event in Proposition~\ref{prop:perturbation} has $\Pr_{\lambdaoccupied(F)\pm\eps}^{F}$-probability at most $\theta/4$ to occur.

  By weak convergence of $F_m$ to $F$ it follows that $F_m^{-1}(U)\to F^{-1}(U)$ almost surely, as $m\to\infty$, if $U$ is uniformly distributed on $[0,1]$ and $F^{-1}$ denotes the generalized inverse of $F$. Consequently, as the projection of $\omega$ onto $K$, denoted by $\omega(K)$, is almost surely finite, we will have $F_m^{-1}(z)\in(F^{-1}(z)-\eps,F^{-1}(z)+\eps)$ for all $(x,z)\in\omega(K)$ for large enough $m$. We may, in particular, choose $m_0$ so that for all $m\ge m_0$
  $$
  \Pr_{\lambdaoccupied(F) \pm \eps}^F \Big( \interior(\mathcal{O}_K, \delta) \subseteq \mathcal{O}_K^{F_m} \subseteq \closure(\mathcal{O}_K, \delta) \Big) > 1 - \theta/4.
  $$
  Combining the above estimates we conclude that~\eqref{eq:convergence} holds, as required.
\end{proof}

\section{Proof of main results}
\label{s:main_results}

In this section we complete the proofs of the three main theorems stated in the introduction.
We first prove Theorem~\ref{thm:sharpnessCrossing}, and then deduce Theorems~\ref{thm:sharpnessVacant} and~\ref{thm:sharpnessOccupied}.
Before starting we recall the content of Corollary~\ref{cor:nontriviality}.
Assuming~\eqref{eq:7}, we have
$$
0<\lambda_0\le\lambda_c\le\lambda_1<\infty,
$$
and the analogous relation holds for $\lambda_c$ replaced by $\lambda_c^\star$
under the stronger condition~\eqref{e:log_condition}.

\subsection{Proof of Theorem~\ref{thm:sharpnessCrossing}}

For any $\lambda>\lambda_1$ we have, by definition, that
$\lim_{r\to\infty}\Pr_\lambda(\cross(3r,r))=1$. The standard inequalities of
Lemma~\ref{lem:standardInequalities} and monotonicity imply that for any
$\lambda>\lambda_1$ and $\kappa>0$,
\begin{equation}\label{eq:proof_cross}
\lim_{r\to\infty}\Pr_\lambda(\cross(\kappa r,r))=1.
\end{equation}
Similarly, for any $\lambda<\lambda_0$ and $\kappa>0$ the limit
in~\eqref{eq:proof_cross} equals 0. Since $\lambda_c\in[\lambda_0,\lambda_1]$,
then Corollary~\ref{cor:box-crossing} shows that for any $\kappa>0$ there exists
a constant $c=c(\kappa)>0$ such that
$$
c<\Pr_{\lambda_c}(\cross(\kappa r,r))<1-c\quad\text{for all }r\ge1.
$$
Moreover, by Corollary~\ref{cor:no_cluster}, there is at $\lambda_c$ almost surely no unbounded component of either kind.
Finally, the proof is completed by Theorem~\ref{thm:equality}, which shows that indeed $\lambda_c=\lambda_0=\lambda_1$.

\subsection{Proof of Theorem~\ref{thm:sharpnessVacant}}

Assume~\eqref{e:log_condition}, in which case we by now know that $\lambda_c^\star=\lambda_0=\lambda_1=\lambda_c$.
Part~\emph{(i)} is then a consequence of Corollary~\ref{cor:arm_events}.
The proof of part~\emph{(ii)} will involve a truncation of the radii distribution, allowing for a comparison with a highly supercritical 1-dependent percolation process on $\Z^2$, for which a standard Peierls argument gives the required exponential decay.

Fix $\lambda>\lambdacritical$. For $m\ge1$, let $\mu_m$ denote the truncation of $\mu$ at $m$, meaning that $\mu_m$ coincides with $\mu$ on $[0,m)$ but that $\mu_m([0,m])=1$. We shall write $\lambda^\star(\mu)$ for the critical parameter associated with radii distribution $\mu$. By continuity of the critical parameter, Theorem~\ref{thm:continuity}, we have for large $m$ that
$$
\lambdacritical(\mu_m)\in[\lambdacritical(\mu),\lambda),
$$
and assume for the rest of this proof that $m$ is chosen accordingly. Given $\gamma>0$, we also assume that $\x>2m$ is fixed so that
$$
\Pr_\lambda^{\mu_m}(\cross(3\x,\x))>1-\gamma.
$$

Now, tile the plane with overlapping $3\x\times\x$ and $\x\times3\x$ rectangles as follows: For each $z\in\Z^2$ consider both a $3\x\times\x$ and a $\x\times3\x$ rectangle positioned so their lower left corners coincide with $2\x z$.
We want to register which of these rectangles that are crossed in the hard direction by an occupied path.
Define thus, for each $z\in\Z^2$, $\eta(z,z + (1,0))$ to be the indicator of the event $\cross(3\x,\x)$ translated by the vector $2\x z$, and $\eta(z,z + (0, 1))$ as the indicator of $\cross(\x,3\x)$ translated by the same vector.

By the choice of $m$ and $\x$, the process $\eta$ is a $1$-dependent bond percolation process on $\Z^2$ with marginal edge probability exceeding $1 - \gamma$.
Therefore, by a standard Peierls argument, there is a constant $\uc{c:exp_decay}>0$ such that the event that $\eta$ contains an open circuit which is contained in $B^\infty(n)$ and surrounds the origin has probability at least $1-\uc{c:exp_decay}\exp\{n/\uc{c:exp_decay}\}$.
However, by construction, an open path in $\eta$ corresponds to an occupied crossing in $\mathcal{O}$, and thus that for all $n\ge1$
$$
\mathbb{P}_\lambda^{\mu_m}\big[ 0\lrv \partial B^\infty(\x n) \big] \leq \uc{c:exp_decay}\exp\{-n/\uc{c:exp_decay}\},
$$
from which part~\emph{(ii)} of the theorem follows.

\subsection{Proof of Theorem~\ref{thm:sharpnessOccupied}}

Under assumption~\eqref{e:alpha_condition} Item~\emph{(i)} is a consequence of part~\emph{(iii)} of Corollary~\ref{cor:arm_events} (see also Remark~\ref{r:D_is_bad}).
So, it remains to argue for Item~\emph{(ii)}. However, if~\eqref{e:alpha_condition} holds, then by Proposition~\ref{c:decouplef1f2} there exists $c > 0$ such that for all $\lambda\le\lambda_c$ and $\x\ge1$
$$
\rho_\lambda(5\x,\x)\le c\x^{-\alpha}.
$$
Now, fix $\lambda<\lambda_c$, so that $\lim_{r\to\infty}\Pr_\lambda(\cross^\star(3\x,\x))=1$. Proposition~\ref{prop:finite_quant} then gives that
$$
\Pr_\lambda(\cross^\star(3\x,\x))\ge 1 - c'r^{-\alpha}
$$
for all $r\ge1$ and some constant $c'$, possibly depending on $\lambda$.
Consequently, for all $r\ge1$,
$$
\Pr_\lambda\big[0 \lro \partial B(r)\big]\,\le\,4\,\Pr_\lambda(\cross(\x,3\x))\,\le\,4c'r^{-\alpha},
$$
as required. This ends the proof of Theorem~\ref{thm:sharpnessOccupied}.

\section{Other models}
\label{s:other_models}

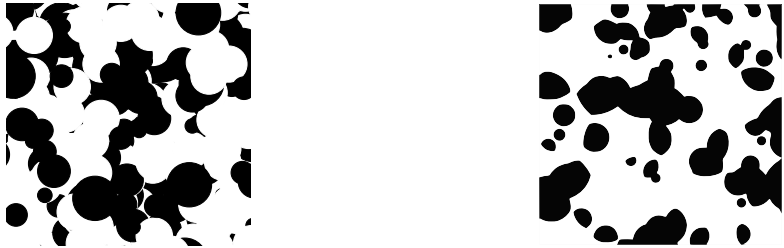
\begin{figure}[ht]
  \centering
  \begin{tikzpicture}[scale=1]
    \begin{scope}
      \clip (-.6, -.6) rectangle (2.6, 2.6);
      \iffinal
      \foreach \z in {1,...,300}
      { \pgfmathrandominteger{\x}{0}{200};
        \pgfmathrandominteger{\y}{0}{200};
        \pgfmathrandominteger{\c}{4}{7};
        \draw[color=black,fill=black] ({.02*\x-1}, {.02*\y-1}) circle (.006 * \c * \c);
        \pgfmathrandominteger{\x}{0}{200};
        \pgfmathrandominteger{\y}{0}{200};
        \draw[color=white,fill=white] ({.02*\x-1}, {.02*\y-1}) circle (.006 * 40);
      }
      \fi
    \end{scope}
    \node at (8, 1) {\includegraphics[width=.20\textwidth]{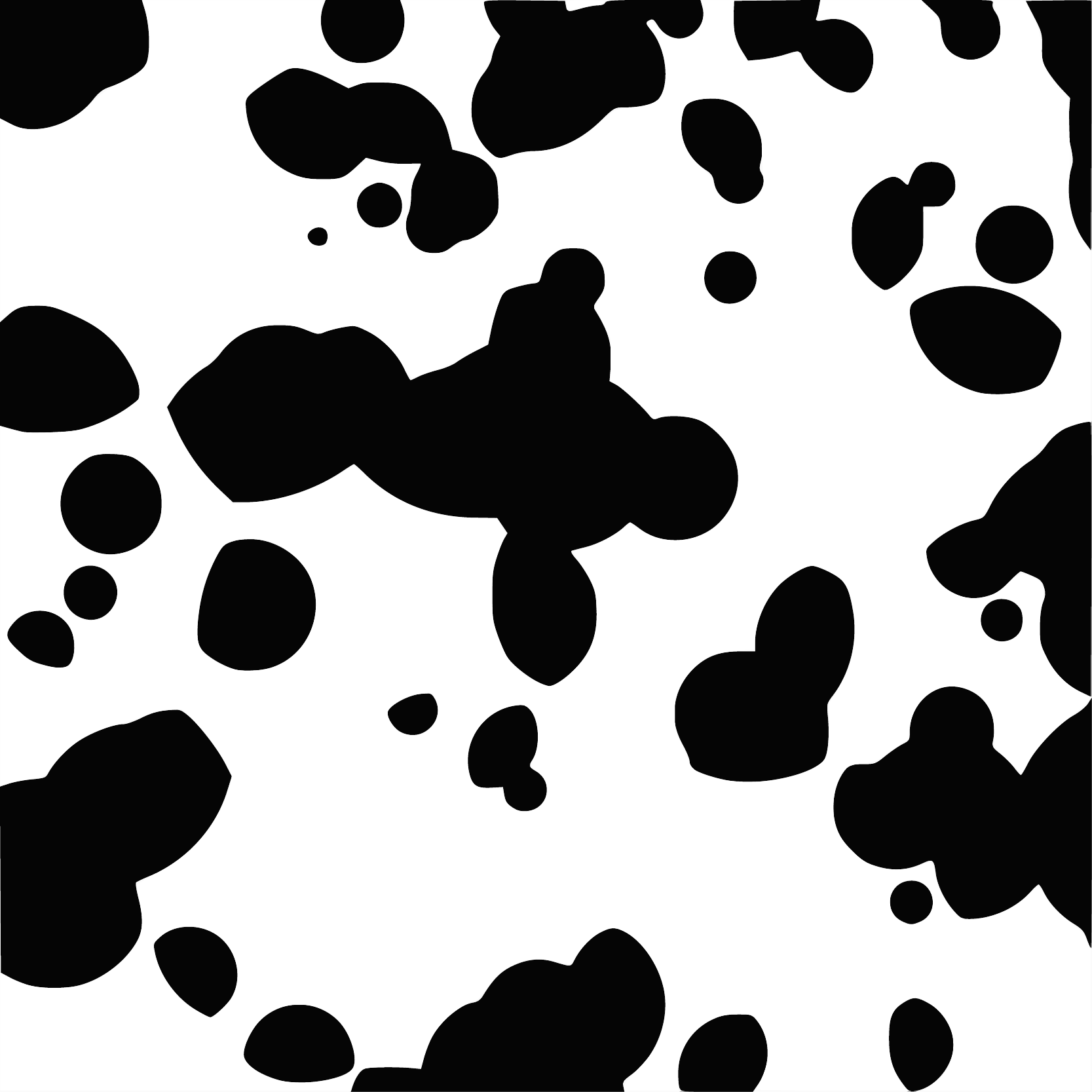}};
  \end{tikzpicture}
  \caption{Simulations of Poisson confetti (left) and Voronoi (right) percolation.
  In the left picture, black confetti have random radii, while the radii of white ones are fixed.
  The Voronoi picture to the right uses $G_0 = 1$, $G_1 = 2$ and $q = 0.8$.
  Note that the Voronoi cells above are determined by segments of ellipses.}
   \label{f:confetti}
\end{figure}

We have in Sections~\ref{s:finite-size}-\ref{s:sharp} provided a framework to prove the existence of a sharp phase transition and some interesting properties of the critical behavior for Poisson Boolean percolation.
In this section we will show how our techniques can be used to obtain alternative proofs for the facts that the critical probability equals $1/2$ for Poisson Voronoi percolation, earlier proved by Bollob\'as and Riordan~\cite{bolrio06-2}, and Poisson confetti percolation, first proved by Hirsch~\cite{hirsch15} (for unit squares) and later M\"uller~\cite{muller} (for unit discs).
Our results will in fact apply to some generalized versions of these models that are not self-dual.
See Figure~\ref{f:confetti} for a simulation of these processes.

The finite-size criterion of Section~\ref{s:finite-size} and the RSW techniques of Section~\ref{s:critical} solely rely on three basic facts: That the probability measure in question is invariant (with respect to translations, right angle rotations and reflections in coordinate axis), positively associated, and features a decay of spatial correlations.
However, the argument used to prove the sharp threshold behavior in Section~\ref{s:sharp} used the further assumption of a Poissonian structure.
All of these properties are satisfied by Poisson Voronoi and confetti percolation.
We shall below illustrate what our techniques give for these models, but emphasize that we make no attempt to sharpen the hypothesises in the results we present.

\subsection{Poisson Voronoi percolation}

The study of random Voronoi tessellations goes back several decades in time, yet it was only about a decade ago that Bollob\`as and Riordan~\cite{bolrio06-2} gave the first proof for the fact that the critical probability in Poisson Voronoi percolation in the plane equals $1/2$; see also \cite{ahlberg2015quenched}.
One of the main difficulties faced in studying the phase transition is to derive RSW techniques that apply in this setting.
In~\cite{bolrio06-2} a weak RSW result was provided, whereas a strong version of RSW was established later by Tassion~\cite{tassion14+}, and was used to prove polynomial decay of the one-arm probability at criticality.

The RSW techniques developed here allows us to consider more general variants of the standard Voronoi model, where black and white cells do not necessarily follow the same law.
Informally speaking, to each point we associate a random `gravitational pull', which gives a bias to the size of the associated tile.
The symmetry between black and white points may be broken by considering different laws for the gravitational pull associated to them.

Consider a Poisson point process of unit density on the space of locally finite counting measures on $\R^2 \times [0, 1] \times [0, 1]$, endowed with Lebesgue measure.
The first two coordinates of a point in the support of a realization $\omega$ will mark the location of the seed in the plane, the third and fourth coordinates will determine the gravitational pull of a given point and its color, respectively.

More precisely, given a parameter $q \in [0, 1]$, we say that the seed $(x, z, s) \in \supp(\omega)$ is \emph{black} if $s \leq q$ and \emph{white} otherwise.
We also fix two non-decreasing functions $G_0, G_1: [0,1] \to (0, \infty)$, which will be applied to the $z$ coordinate of $(x, z, s)$ to determine the gravitational pull of a given seed ($G_0$ will be used for black seeds and $G_1$ for white).
Given $q \in [0, 1]$, define the occupied (black) set as
\begin{equation}
  \label{e:O_Voronoi}
  \mathcal{O}_q := \bigcup_{\substack{(x, z, s) \in \omega;\\s \leq q}} \Big\{ y \in \R^2 : \frac{|y - x|}{G_0(z)} \leq \frac{|y - x'|}{G_1(z')} \text{ for all $(x', z', s') \in \omega$ with $s' > q$} \Big\}.
\end{equation}
As before, let $\mathcal{V}_q := \mathbb{R}^2 \setminus \mathcal{O}_q$ denote the corresponding vacant set.
Notice that the cell corresponding to a given seed may be disconnected.
Finally, we define the critical values of the parameter $q$ as
\begin{equation}\label{eq:q_c}
\begin{aligned}
q_c:=&\sup\big\{q\in[0,1]:\Pr\big[0\overset{\tiny \mathcal{O}_q}\longleftrightarrow\infty\big]=0\big\},\\
q_c^\star:=&\inf\big\{q\in[0,1]:\Pr\big[0\overset{\tiny \mathcal{V}_q}\longleftrightarrow\infty\big]=0\big\}.
\end{aligned}
\end{equation}

We will here only consider the case when $G_0$ and $G_1$ take values in some interval $[a, b] \subset (0, \infty)$.
It would be interesting to attempt to relax these assumptions.
When $G_0 = G_1$, then the gravitational pull of all black and white seeds are equally distributed.
In this case the model is self-dual, and we recover the equation $q_c+q_c^\star=1$; see Section~\ref{ss:remarks} below.
The techniques developed in this paper may easily be adapted to prove the following result for Poisson Voronoi percolation.

\begin{thm}\label{thm:Voronoi}
Assume that $G_0$ and $G_1$ take values in some interval $[a,b]\subset(0,\infty)$. Then, $q_c=q_c^\star$ and their common value is strictly between zero and one.
Moreover,
\begin{enumerate}[(i)]
\item for all $q<q_c$, there exists $c = c(q)>0$ such that
$$
\Pr \big[ 0\overset{\tiny \mathcal{O}_q}\longleftrightarrow\partial B(r) \big] \le \exp\{-c r\}.
$$
\item for $q=q_c$ and all $\kappa>0$, there exist $c=c(\kappa)>0$ and $\alpha>0$ such that for every $r \geq 1$
\begin{equation*}
  c < \Pr(\cross(\kappa r,r)) < 1 - c\quad\text{and}\quad  \mathbb{P} \big[ 0 \overset{\tiny \mathcal{O}_q}\longleftrightarrow \partial B(r) \big] \leq r^{-\alpha}.
\end{equation*}
\item for all $q>q_c$, there exists $c=c(q)>0$ such that
$$
\Pr \big[ 0\overset{\tiny \mathcal{V}_q}\longleftrightarrow\partial B(r) \big] \le \exp\{-c r\}.
$$
\end{enumerate}
Moreover, at $q_c$ there is almost surely no unbounded cluster of either kind.
\end{thm}

We will for the remainder of this subsection outline the proof of this result, based on the techniques developed over the previous sections.
First of all, we note that the law of the above model is invariant with respect to translations, rotations and reflections.
The measure is further positively associated, meaning that is satisfies an FKG inequality similar to~\eqref{eq:3}, for events that are increasing in the addition of black points and removal of white points. This can be seen by adapting the proof of~\cite[Lemma~8.14]{bolrio06} for the case when $G_0=G_1$ are constant.
What remains in order for the techniques of Sections~\ref{s:finite-size} and~\ref{s:critical} to apply is an estimate on the spatial correlations.

\begin{lma}
  \label{l:decouple_voronoi}
  Fix a bounded measurable set $D \subseteq \mathbb{R}^2$ and an arbitrary function $f(\mathcal{O}_q)$ satisfying $f \in \sigma(Y_v; v \in D)$.
  Then, defining the event
  \begin{equation}
    G_{D, \x} = \Big\{ \text{for every $y \in D$, there is a point $(x, z, s)$ in $\omega$ with $d(y, x) \leq \frac{a}{b} \x$} \Big\},
  \end{equation}
  we have that $\ind_{G_{D, \x}} \cdot f$ only depends on the restriction of $\omega$ to $B(D, \x) \times [0, 1] \times [0, 1]$.
  Moreover, for every $\eps>0$ there exists a constant $c=c(\eps)>0$ such that for any $D\subseteq B^\infty(\x)$ we have $\mathbb{P}(G_{D, \eps\x}) \geq 1 - \exp\{- c \x^2 \}$.
\end{lma}

\begin{proof}
  The proof of the first claim follows the same steps as those of Lemma~1.1 in \cite{tassion14+}.
  Roughly speaking, the argument goes like this.
  Given $x \in \mathbb{R}^2$, if a seed in $\omega$ is located within distance $\ell$ from $x$, then no other seed outside $B(x, (b/a)\ell)$ can influence the state of $x$.
  The second claim follows from a simple large deviations bound for Poisson random variables.
\end{proof}

As a first consequence of this fact we obtain an estimate on the decay of spatial correlations.
As in Definition~\ref{def:rho}, let $f_1, f_2:\mathcal{P}(\R^2)\to[-1,1]$ be two functions satisfying $f_1(\mathcal{O}_q)\in\sigma\big(Y_v;v\in B^\infty(r)\big)$ and $f_2\in\sigma\big(Y_v;v\in B^\infty(r+t)\setminus B^\infty(r+s)\big)$, for some $t>s$. Then, if $G=G_{B^\infty(r),s/2}\cap G_{B^\infty(r+t)\setminus B^\infty(r+s),s/2}$, we obtain
\begin{equation}
  \label{e:rho_voronoi}
  \begin{split}
    \mathbb{E}[f_1 f_2] \,& \leq\, \mathbb{P}(G^c) + \mathbb{E}[\ind_G \cdot f_1 f_2]\\
    &\leq \,2\mathbb{P}(G^c)+\E[f_1']\E[f_2']\\
    &\leq \,6\Pr(G^c)+\E[f_1]\E[f_2],
  \end{split}
\end{equation}
where $f_1'$ and $f_2'$ denote $f_1$ and $f_2$ evaluated at the restriction of $\omega$ to $B^\infty(r+s/2)$ and $B^\infty(r+t)\setminus B^\infty(r+s/2)$ and are hence independent. In particular, when $s=\eps r$ and $t=\frac1\eps r$, which has been the case throughout this paper, the correlation between $f_1$ and $f_2$ decays as $\exp\{- c r^2\}$ as in Lemma~\ref{l:decouple_voronoi}.
Consequently, the techniques of Sections~\ref{s:finite-size} and~\ref{s:critical} apply, and we obtain the existence of $q_0$ and $q_1$, defined analogously as $\lambda_0$ and $\lambda_1$, satisfying
$$
0<q_0\le q_c^\star, \, q_c\le q_1<1,
$$
and such that part~\emph{(i)} of Theorem~\ref{thm:Voronoi} holds for $q<q_0$, part~\emph{(ii)} for $q\in[q_0,q_1]$, and part~\emph{(iii)} for $q>q_1$. To complete the proof of Theorem~\ref{thm:Voronoi}, it remains to prove that $q_0=q_1$.

In order to show that $q_0=q_1$ we again enlarge our probability space as in Section~\ref{ss:two_stage}, to enable a two-stage construction of our process.
We thus assume that $\overline\omega$ is a Poisson point process on the space of locally finite counting measures on $\big(\R^2\times[0,1]\times[0,1]\big)\times[0,1]$.
The law of $\overline\omega$ (denoted $\overline{\mathbb{P}}_\lambda$) is chosen with a large density $\lambda=m\ge1$ and the fifth coordinate will be used to thin the process down to density one, which we then use to define $\mathcal{O}_q$ as in~\eqref{e:O_Voronoi}.

We denote by $\omega$ the projection of $\overline\omega$ onto $\mathbb{R}^2 \times [0, 1] \times [0, 1]$ (the first four coordinates) and by $\underline\omega$ the further projection onto $\mathbb{R}^2 \times [0, 1]$.
Define
$$
W_q:=\overline\Pr_\lambda(\cross(4\x,\x)|\underline\omega).
$$
As in Section~\ref{ss:two_stage}, we note that $W_q$ for  $\overline\Pr_\lambda$-almost every $\overline\omega$ coincides with the expectation of some Boolean function $f_{\underline\omega}:\{0,1\}^{\underline\omega}\times\{0,1\}^{\underline\omega}\to\{0,1\}$ evaluated with respect to product measure with density $q$ and $1/m$ for choosing `color' and `presence', respectively, of each point in $\underline\omega$. This function is increasing in the choice of color, so we obtain via the Margulis-Russo formula that
\begin{equation}\label{eq:russo_voronoi}
\frac{dW_q}{dq}=\sum_{i}\overline\Pr_\lambda \big( \,\text{point $i$ is present and its color is pivotal for }f_{\underline\omega}\, \big| \,\underline\omega\, \big).
\end{equation}

Note that the concept of pivotality arising here is with respect to \emph{color} and not \emph{presence} as in Section~\ref{s:sharp}.
The crucial observation is that is that switching the color of a point has a larger potential change than switching its presence.
That is, for the existence of a black crossing, a black point is better than no point, and no point is better than a white point.
Hence, from~\eqref{eq:russo_voronoi}, we obtain that
\begin{equation*}
  \begin{aligned}
    \frac{dW_q}{dq} \,&=\, \overline\E_\lambda \Big[ \sum_i \overline\Pr_\lambda \big( \, \text{point $i$ is present and its color is pivotal} \, \big| \, \omega \, \big) \Big| \, \underline\omega \, \Big]\\
    &\ge\,\overline\E_\lambda \Big[ \sum_i\overline\Pr_\lambda \big( \, \text{point $i$ is present and its presence is pivotal} \, \big| \,\omega \, \big) \Big|\,\underline\omega\,\Big]\\
    &=\, \frac{1}{m} \overline\E_\lambda \Big[ \sum_i \Inf^q_i(f_{\underline\omega}) \Big|\,\underline\omega\,\Big]\\
    &\ge \, \frac{1}{m} \overline\E_\lambda \Big[ \frac{c}{\log(m)} \Var_\lambda(f_{\underline\omega}|\omega)\log\frac{1}{\max_i \big( \Inf_i^q(f_{\underline\omega}) \big) } \,\Big|\,\underline\omega\,\Big],
  \end{aligned}
\end{equation*}
where in the last step we applied~\eqref{eq:talagrand}, and $\Inf_i^q(f_{\underline\omega})$ denotes the conditional probability, given $\omega$, that changing point $i$ from present to absent changes the outcome of $f_{\underline\omega}$.
Controlling the variance $\Var_\lambda(f_{\underline\omega}|\omega)$ can now be done exactly as in Section~\ref{ss:controlling_variance}.
Also the influences $\Inf_i^q(f_{\underline\omega})$ with respect to presence of a point can be estimated similarly as it was done in Section~\ref{ss:controlling_influences}, but requires a slight modification of the decoupling inequalities~\eqref{e:f_decouple_in_G} and~\eqref{e:f_given_omega_decouples}.

Let $D_1,D_2,\ldots, D_J$ be well separated sets as in~\eqref{e:A_j_separated}, and let $A_1,A_2,\ldots,A_J$ be events where $A_j$ is determined by the restriction of $\mathcal{O}_q$ to $D_j$.
On the event $G=\mcap_{j\le J}G_{D_j,r/(2n)}$ the events $A_1,A_2,\ldots,A_J$ are, by Lemma~\ref{l:decouple_voronoi}, determined by restrictions of $\overline\omega$ to disjoint subsets of the plane, and are hence independent.
We obtain, therefore, the following variant of~\eqref{e:f_decouple_in_G} for the Voronoi model
\begin{equation}\label{e:V_f_decouple_in_G}
\overline\Pr_\lambda\big(A_1\cap\cdots\cap A_J\big|\omega\big)\,\le\,\overline\Pr_\lambda(G^c|\omega)+\prod_{j\le J}\Big[\overline\Pr_\lambda(A_j|\omega)+\overline\Pr_\lambda(G^c|\omega)\Big].
\end{equation}
Since $\overline\Pr_\lambda(G^c)\to0$ fast, also $\overline\Pr_\lambda(G^c|\omega)\to0$ fast with probability tending to one. This suffices for the application of~\eqref{e:f_decouple_in_G} in the proof of Lemma~\ref{lma:q-arm}.

In order to obtain a variant for~\eqref{e:f_given_omega_decouples} we observe that the difference between $\overline\Pr_\lambda(A_j|\omega)$ and the conditional probability (given $\omega$) that $A_j$ occurs with respect to the restriction of $\overline\omega$ restricted to the set $B(D_j,r/(2n))$ is bounded by $\overline\Pr_\lambda(G^c|\omega)$. Let $H=[\overline\Pr_\lambda(G^c|\omega)\le\delta]$. Hence, arguing similarly as for~\eqref{e:f_given_omega_decouples}, we obtain that
\begin{equation}\label{e:V_f_given_omega_decouples}
\overline\Pr_\lambda\bigg(H \cap \bigcap_{j\le J}\big[\overline\Pr_\lambda(A_j|\omega)\le\delta\big]\bigg)\,\le\,\prod_{j\le J}\Big[\overline\Pr_\lambda\big(\overline\Pr_\lambda(A_j|\omega)\le3\delta\big)+\overline\Pr_\lambda(H^c)\Big].
\end{equation}
Since $\overline\Pr_\lambda(H^c)\to0$ fast, as a consequence of Lemma~\ref{l:decouple_voronoi}, the two expressions~\eqref{e:V_f_decouple_in_G} and~\eqref{e:V_f_given_omega_decouples} replace~\eqref{e:f_decouple_in_G} and~\eqref{e:f_given_omega_decouples} in the proof of Proposition~\ref{p:max_inf_small}. The proof that $q_0=q_1$ is now a straightforward adaptation of the arguments presented in Section~\ref{s:sharp}, and this ends the outline of the proof of Theorem~\ref{thm:Voronoi}.

\subsection{Poisson confetti percolation}

Confetti percolation, or the `dead leaves' model, was introduced by Jeulin in~\cite{jeulin97}. In this model, black and white confetti `rain down' on the plane according to a Poisson point process, and each point in the plane is colored according to the color of the first confetti to cover it. In the case of circular confetti with fixed diameter, Benjamini and Schramm~\cite{bensch98} conjectured that the critical probability for this model equals $1/2$. This was later confirmed by Hirsch~\cite{hirsch15} for square shaped confetti and by M\"uller~\cite{muller} for circular confetti. Just as in the settings of Poisson Boolean and Voronoi percolation, our techniques allow us to handle confetti of random radii, whose laws may differ between black and white.

Consider a Poisson point process on the space of locally finite counting measures on $\R^2 \times \R_+ \times [0, 1]^2$.
Here again the first two coordinates will assign a location in $\R^2$ to each confetti, while the third coordinate will denote the fall-time of a confetti and will be used to order overlapping confetti.
Finally, the fourth and fifth coordinates will help us determine the radius and color of the confetti, respectively.

As before, we fix $q \in [0, 1]$ and, given a realization $\omega$ with density one, declare a point $(x, t, z, s) \in \omega$ \emph{black} if $s \leq q$ and \emph{white} otherwise.
We also fix non-decreasing functions $G_0, G_1 : [0, 1] \to (0, \infty)$ that will be applied to the $z$ coordinate, in order to determine the radius of the confetti to be placed at $x\in\R^2$ ($G_0$ for black points and $G_1$ for white).
More precisely, given a realization $\omega$ of the Poisson point process, let
\begin{equation*}
\mathcal{O}_q := \bigcup_{\substack{(x, t, z, s) \in \omega\\s \leq q}} \left\{ y \in \R^2:
\begin{split}
& |y - x| \le G_0(z), \text{ and } t < t' \text{ for all $(x', t', z', s') \in \omega$}\\
& \text{such that $|y - x'| \le G_1(z')$ and $s' > q$}
\end{split}
\right\}.
\end{equation*}
Finally, let $\mathcal{V}_q := \mathbb{R}^2 \setminus \mathcal{O}_q$ and define the critical parameters $q_c$ and $q_c^\star$ as in~\eqref{eq:q_c}.

For the sake of simplicity, we assume here rather light tails for distribution of the radii induced by $G_0$ and $G_1$.
If these two functions are the same, then black and white confetti have the same radii distributions, which yields $q_c+q_c^\star=1$ by self-duality, and together with the following theorem that $q_c=1/2$; see Section~\ref{ss:remarks} below.

\begin{thm}\label{thm:confetti}
Assume that $G_0$ and $G_1$ satisfy $\Pr[G_i(Z) \geq r] \leq r^{-100}$, for $i = 0, 1$ and large $r$, where $Z$ is uniformly distributed on $[0,1]$. Then, $q_c=q_c^\star$ are strictly between zero and one, and
\begin{enumerate}[(i)]
\item for all $q<q_c$, there exists $c=c(q) > 0$ such that
$$
\Pr \big[ 0\overset{\tiny \mathcal{O}_q}\longleftrightarrow\partial B(r) \big] \le\frac{1}{c}r^{-10}.
$$
\item for $q=q_c$ and all $\kappa>0$, there exist $c=c(\kappa)>0$ and $\alpha>0$ such that, for every $r \geq 1$
$$
c<\Pr(\cross(\kappa r,r))<1-c\quad\text{and}\quad
  \mathbb{P}\big[ 0\overset{\tiny \mathcal{O}_q}\longleftrightarrow\partial B(r) \big] \leq r^{-\alpha}.
$$
\item for all $q>q_c$, there exists $c=c(q)>0$ such that
$$
\Pr \big[ 0\overset{\tiny \mathcal{V}_q}\longleftrightarrow\partial B(r) \big] \le\frac{1}{c}r^{-10}.
$$
\end{enumerate}
Moreover, at $q_c$ there is almost surely no unbounded cluster of either kind.
\end{thm}

The proof of this theorem is analogous to that of Theorem~\ref{thm:Voronoi} above.
The only distinction lies in the treatment of the spatial dependence. The following lemma will replace Lemma~\ref{l:decouple_voronoi}.

\begin{lma}
  Fix a bounded measurable set $D \subseteq \mathbb{R}^2$ and $r \geq 1$, $M \geq 1$ and define
  \begin{equation*}
    G_{D, r, M} = \left\{
    \begin{array}{c}
      \text{for every $(x, t, z, s) \in \omega \cap B(D, Mr)$ with $t \leq r$, $G_0(z) \vee G_1(z) \leq r/M$}\\
      \text{for every $(x,t,z,s) \in \omega\cap B(D,Mr)^c$ with $t\le r$, $G_0(z) \vee G_1(z) \leq d(x,D)$}\\
      \text{for every $y \in D$, there is $(x, t, z, s) \in \omega$ with $t \leq r$, $G_0(z) \wedge G_1(z) \geq |y - z|$}
    \end{array}
\right\}.
  \end{equation*}
  Then, $\lim_{r \to \infty} r^{10} \mathbb{P}(G_{D, r, M}) = 1$ and for an arbitrary function $f(\mathcal{O}_q) \in \sigma(Y_x; x \in D)$, we have that $\ind_G \cdot f$ is measurable with respect to $\omega$ restricted to $B(D, r/M) \times [0, r] \times [0, 1] \times [0, 1]$.
\end{lma}

The proof of this lemma follows from a simple large deviations bound and is omitted.
As a consequence of the lemma, spatial correlations will decay as $r^{-10}$.
This rate of decay is more than enough to follow the outlined proof of Theorem~\ref{thm:Voronoi} and obtain a proof also for Theorem~\ref{thm:confetti}.

\subsection{The consequence of self-duality}\label{ss:remarks}

For the two models considered in this section, the case when the two functions $G_0$ and $G_1$ are equal stands out due to the self-duality retained in the model. Self-duality refers to the property that the occupied set at parameter $q$ has the same law as the vacant set has at parameter $1-q$. In particular, there is no unbounded occupied component at parameter $q$ if and only if there is no unbounded vacant component at parameter $1-q$ (almost surely). As the critical parameters are defined as the supremum for which this holds it follows that $q_c+q_c^\star=1$, in analogy to the Bernoulli case in~\eqref{eq:2}. Together with the equality $q_c=q_c^\star$ of the two parameters, obtained in Theorems~\ref{thm:Voronoi} and~\ref{thm:confetti}, it follows that $q_c=q_c^\star=1/2$, and we recover the results of~\cite{bolrio06-2},~\cite{hirsch15} and~\cite{muller}.

\begin{remark}
  There are other Poisson percolation processes that could in principle be studied with the techniques we have employed.
  These include for instance Brownian interlacements (see~\cite{sznitman13}) in $\R^d$ intersected with a plane, and Poissonian cylinders (see~\cite{TW10b}) in $\R^d$ intersected with a plane.
  A rather general setting for Poisson Boolean percolation but with dependence between the radii assigned to different points has been considered by Ahlberg and Tykesson~\cite{ahltyk}, however that study does not address the sharpness of the threshold.
  It would be interesting to investigate to what extent the techniques developed here applies in that setting.
\end{remark}

\section{Open problems}\label{s:open}

We have in this paper developed techniques for the study of the sharpness of the threshold for a wide range of models.
The weakest link in this scheme is the Poissonian assumption necessary to deduce that the critical regime is indeed a single point and not an interval.
It would be interesting to develop an alternative argument for the existence of a sharp threshold beyond the Poissonian setting considered here.

It would further be interesting to pursue optimal conditions to study Poisson Voronoi and confetti percolation.
Here we only aimed at illustrating how our techniques may be applied to other settings, but we have not pursued to make these applications optimal.
Note also that the techniques we have used do not apply to Poisson Boolean percolation models where the discs are replaced by some other deterministic shape which is not sufficiently symmetric.

Besides extending the current techniques to other models, one could be interested in answering several questions that are known for Bernoulli percolation.
Examples of such problems would be to define the incipient infinite cluster, study noise sensitivity for $\cross(r, r)$ and investigate scaling relations, just to name a few.

\bibliographystyle{alpha}
\bibliography{bib}

\end{document}